\documentclass{article}
\usepackage[utf8]{inputenc}
\usepackage[T1]{fontenc}
\usepackage{amsmath,amssymb,amsthm,amscd}
\usepackage{mathrsfs}
\usepackage{lipsum}
\usepackage{tikz-cd}
\usetikzlibrary{matrix,arrows,decorations.pathmorphing}
\usepackage{enumitem}
\setlist[description]{leftmargin=\parindent,labelindent=\parindent}
\usepackage{xcolor}
\usepackage{fullpage}
\usepackage{url}
\usepackage{cite}
\usepackage[colorlinks=true,linkcolor=blue!60!black,citecolor=green!45!black,urlcolor=blue!60!black]{hyperref}
\usepackage{sample}
\theoremstyle{plain}
\newtheorem{theorem}{Theorem}[section]
\newtheorem{proposition}[theorem]{Proposition}
\newtheorem{lemma}[theorem]{Lemma}
\newtheorem{corollary}[theorem]{Corollary}
\theoremstyle{definition}
\newtheorem{definition}[theorem]{Definition}
\newtheorem{convention}[theorem]{Convention}
\theoremstyle{remark}
\newtheorem{remark}[theorem]{Remark}

\numberwithin{equation}{section}

\title{Equivariant Levi restriction}
\author{Mark Macerato
  \and
  Wenqing Wei
  \and
  Peisheng Yu
}

\newcommand{\Addresses}{{
  \bigskip
  \footnotesize

  Mark Macerato, \textsc{University of Minnesota, Twin Cities, Department of Mathematics, Minneapolis, MN, 55455}\par\nopagebreak
  \textit{E-mail address}, Mark Macerato: \texttt{macer008@umn.edu}

  \medskip

  Wenqing Wei, \textsc{University of California Berkeley, Department of Mathematics, Berkeley, CA, 94720, USA}\par\nopagebreak
  \textit{E-mail address}, Wenqing Wei: \texttt{weiwenqing@berkeley.edu}

  \medskip

  Peisheng Yu, \textsc{University of California Berkeley, Department of Mathematics, Berkeley, CA, 94720, USA}\par\nopagebreak
  \textit{E-mail address}, Peisheng Yu: \texttt{peishengyu@berkeley.edu}

}}

\date{}

\begin{document}
\maketitle

\begin{abstract}
Let $G$ be a connected reductive group over $\bbC$, $L\subset G$ a Levi subgroup, and $\cL\subset\cG$ their Langlands duals. Macerato \cite{macerato2023leviequivariantrestrictionsphericalperverse} computed the $\cL(\calO)$-equivariant cohomology of the $!$-restriction to $\GrL$ of the regular spherical perverse sheaf on $\GrG$, and identified it with the functions on the partial Kostant--Whittaker reduction $T^*(G/U,\psiL)$, where $\psiL$ is the possibly degenerate Whittaker character of $U$ determined by $L$. We add loop rotation. Our main result is: the $\LO$-equivariant cohomology of the same object is the ring $\Dh(G/U,\psiL)$ of asymptotic $\psiL$-twisted differential operators on $G/U$. For a single $\IC_V$ the answer is the space of intertwiners $\bM_\hbar(\psiL)\to V\otimes\bM_\hbar(\psiL)$ of the universal Verma--Whittaker module. This interpolates between results of Bezrukavnikov--Finkelberg ($L=G$) and Ginzburg--Riche ($L=T$).
\end{abstract}

\tableofcontents

\section{Introduction}
\subsection{Equivariant geometric Satake}\label{sec:intro-satake}
Let $G$ be a connected reductive group over $\bbC$ and $\cG$ its Langlands dual. Put $\calK=\bbC((t))$, $\calO=\bbC[[t]]$, and let $\GrG=\cG(\calK)/\cG(\calO)$ be the affine Grassmannian, with the action of $\GO$ ,where $\bbG_m$ acting by loop rotation. Geometric Satake is a symmetric monoidal equivalence
\begin{equation}\label{eq:satake}
\bbS_G\colon \Perv_{\cG(\calO)}(\GrG)\xrightarrow{\ \sim\ }\Rep(G)
\end{equation}
under which the fiber functor to $\Vect$ corresponds to $\H^*(\GrG,-)$. Every spherical perverse sheaf is  $\GO$-equivariant. We write $\IC_V:=\bbS_G^{-1}(V)$ and $\IC_\lambda=\IC_{V^\lambda}$.

After turning on the loop rotation, $\eqref{eq:satake}$ is quantized by
Bezrukavnikov and Finkelberg \cite{bezrukavnikov2026equivariantsatakecategorykostantwhittaker}, where they computed the loop-rotation equivariant cohomology, and expressed in terms of Kostant--Whittaker reduction:
\begin{equation}\label{eq:BF}
\H_{\GO}(\GrG,\calF)\ \cong\ (\bbS_G(\calF)\otimes\calW_\hbar)^{\fraku,\psi},
\end{equation}
where $\calW_\hbar$ is the universal Whittaker module of the asymptotic enveloping algebra $\Uh(\frakg)$ for a non-degenerate character $\psi$ of $\fraku$, and $(-)^{\fraku,\psi}$ denotes $\psi$-Whittaker vectors (\S\ref{sec:universal}). For the regular sheaf $\IC_\reg=\bigoplus_{\lambda\in\Lambda^+}(V^\lambda)^*\otimes\IC_\lambda$ this becomes an isomorphism between algebras
\begin{equation}\label{eq:BF-reg}
\H_{\GO}(\GrG,\IC_\reg)\ \cong\ \Dh(G/U,\psi):=\bigl(\calO(G)\otimes\calW_\hbar\bigr)^{U},
\end{equation}
where $\Dh(G/U,\psi)$ is the quantum Hamiltonian reduction of $\Dh(G)=\calO(G)\otimes\Uh(\frakg)$ \cite[\S3.1]{ginzburg2014differentialoperatorsguaffine} at $\psi$ with respect to right translations by $U$. For the extended torus $A=\cT\times\bbG_m$, equivariant formality and Lemma~\ref{lem:GK} give
\begin{equation}\label{eq:BF-torus}
\H_A(\GrG,\calF)\ \cong\ \bar\fraku\backslash\bigl(\bbS_G(\calF)\otimes\calW_\hbar\bigr).
\end{equation}
This is also implicitly used in \cite{bezrukavnikov2026equivariantsatakecategorykostantwhittaker}.

\subsection{Equivariant corestriction to the torus}\label{sec:intro-torus}
Let $i_T\colon\GrT\hookrightarrow\GrG$ be the inclusion. Ginzburg and Riche \cite[Theorem 2.2.4, Lemma 2.4.1]{ginzburg2014differentialoperatorsguaffine} computed
\begin{equation}\label{eq:GR}
\H_A(\GrT,i_T^!\calF)\ \cong\ (\bbS_G(\calF)\otimes\bM_\hbar)^{\fraku}
=\Hom_{\Uh(\frakg)}\bigl(\bM_\hbar,\bbS_G(\calF)\otimes\bM_\hbar\bigr),
\end{equation}
with $\bM_\hbar=\Uh(\frakg)/\Uh(\frakg)\fraku$ the universal Verma module. For $\calF=\IC_\reg$:
\begin{equation}\label{eq:GR-reg}
\H_A(\GrT,i_T^!\IC_\reg)\ \cong\ \Dh(G/U)=(\calO(G)\otimes\bM_\hbar)^U .
\end{equation}

\subsection{Equivariant corestriction to a Levi}\label{sec:intro-levi}
Let $T\subset L\subset G$ be a Levi subgroup, $\cL\subset\cG$ the dual Levi, and let $i_L\colon\GrL\hookrightarrow\GrG$ be the inclusion between affine grassmannians. Let $\psiL$ be the partial Whittaker character of $U$ determined by $L$ ($\psiL$ is non-zero on the simple roots of $L$, and zero on other simple roots). In \cite[Theorem 1.5.2, Corollary 1.5.4]{macerato2023leviequivariantrestrictionsphericalperverse}, Macerato proved that
\begin{equation}\label{eq:Mac}
\H_{\cL(\calO)}(\GrL,i_L^!\IC_\reg)\ \cong\ \calO\bigl(T^*(G/U,\psiL)\bigr),\qquad T^*(G/U,\psiL)=G\times_U(\psiL+\fraku^\perp),
\end{equation}
and, on each $V$-isotypic components,
\begin{equation}\label{eq:Mac-V}
\H_{\cL(\calO)}(\GrL,i_L^!\IC_V)\ \cong\ \bigl(V\otimes\calO(\psiL+\fraku^\perp)\bigr)^U .
\end{equation}
However, loop rotational version is not considered there.

\subsection{Main result}\label{sec:intro-main}
We quantized the main result \eqref{eq:Mac} in \cite{macerato2023leviequivariantrestrictionsphericalperverse}. Our main theorem is:
\begin{theorem}\label{thm:main}
There is an isomorphism of $\Zh(\frakl)$-algebras
\begin{equation}\label{eq:main}
\Psi\colon\ \Dh(G/U,\psiL):=\bigl(\calO(G)\otimes \Uh(\frakg)/\Uh(\frakg)\fraku^{\psiL}\bigr)^{U}\ \xrightarrow{\ \sim\ }\ \H_{\LO}(\GrL,i_L^!\IC_\reg).
\end{equation}
For every $V\in\Rep(G)$ it restricts to
\begin{equation}\label{eq:main-V}
\Psi_V\colon\ \bigl(V\otimes\bM_\hbar(\psiL)\bigr)^{\fraku,\psiL}=\Hom_{\Uh(\frakg)}\bigl(\bM_\hbar(\psiL),V\otimes\bM_\hbar(\psiL)\bigr)\ \xrightarrow{\ \sim\ }\ \H_{\LO}(\GrL,i_L^!\IC_V),
\end{equation}
where $\bM_\hbar(\psiL)=\Uh(\frakg)/\Uh(\frakg)\fraku^{\psiL}$ is the universal Verma--Whittaker module (\S\ref{sec:verma-whittaker}).
\end{theorem}

Here $\Zh(\frakl)$ is the center of $\Uh(\frakl)$; it acts on the source of \eqref{eq:main-V} by the right action (postcomposition, \S\ref{sec:intertwiners}), and on the target through $\Zh(\frakl)\cong\H_{\LO}(\pt)$. The ring structure on the right of \eqref{eq:main} comes from the lax monoidal structure of $i_L^!$ \cite[\S2.2]{macerato2023leviequivariantrestrictionsphericalperverse}. For $L=G$ we recover \eqref{eq:BF}, for $L=T$ we recover \eqref{eq:GR}, and at $\hbar=0$ we recover \eqref{eq:Mac}, which is the isomorphism constructed in \cite{macerato2023leviequivariantrestrictionsphericalperverse}.

\begin{remark}[Normalizations]\label{rem:gradings}
Both sides of \eqref{eq:main} are graded and decompose over $\pi_0(\GrL)=\Lambda/\bbZ R_L$, resp.\ over the characters of $Z(L)$. To make $\Psi$ homogeneous one must shift the grading on each component of $\GrL$; these are the shifts of the normalized corestriction $i_L^{!,\natural}$ of \cite[\S1.3]{macerato2023leviequivariantrestrictionsphericalperverse}, and we suppress them. Similarly, we suppress the normalizations needed for matching the equivariant parameters and the algebraic parameters.
\end{remark}

\subsection{Strategy of the proof}\label{sec:intro-strategy}
Following  \cite{ginzburg2014differentialoperatorsguaffine,macerato2023leviequivariantrestrictionsphericalperverse}, we will use equivariant localization. Write $\alpha$ for roots in $R^+\setminus R^+_L$ and $\beta$ for roots in $R^+_L$. $\H_A(\GrL,i_L^!\calF)$ is free over $\Sh=\H_A(\pt)$, and we study it after two localizations:
\begin{itemize}
\item[(1)] inverting the $\alpha$-walls $\check\alpha+n\hbar$, $i_L^!\calF$ becomes the hyperbolic restriction $\Res^{\cG}_{\cL}\calF$, computed by equivariant Satake for $\cL$ (\S\ref{sec:case1});
\item[(2)] inverting the $\beta$-walls $\check\beta+m\hbar$, one localizes further to $\GrT$ and lands in case of \cite{ginzburg2014differentialoperatorsguaffine} (\S\ref{sec:case2}).
\end{itemize}
On the algebraic side, Case~1 uses Kostant's theory and the simplicity of standard modules off the $\alpha$-walls; Case~2 uses Backelin's functor $\widehat\Gamma_{\psiL}$ and the existence of the universal Whittaker vector in the completed Verma module after localization. The desired existence is controlled by the Shapovalov determinant. The two cases are glued along the regular semisimple locus, using a Hartogs property of $\Dh(G/U,\psiL)$ (\S\ref{sec:hartogs}).

However, the gluing is quite non-trivial, and that's what makes the problem more complicated. For $L\neq T$ the Case~2 isomorphism, given by evaluating intertwiners on the universal Whittaker vector, differs on the regular semisimple locus from the isomorphism in Case~1. On the $\nu$-weight piece the two differ, up to a constant in $T$, by a ratio of $\Gamma$-functions:
\[
\delta^L_\nu(y)=\prod_{\beta\in R_L^+}P_{\langle\nu,\check\beta\rangle}\bigl(\langle y+\hbar\rho,\check\beta\rangle\bigr)^{-1}
=\frac{\Gamma_L(y+\hbar\rho)}{\Gamma_L(y+\hbar(\nu+\rho))},\qquad
\Gamma_L(Z)=\prod_{\beta\in R^+_L}\hbar^{\langle Z,\check\beta\rangle/\hbar}\,\Gamma\bigl(1+\langle Z,\check\beta\rangle/\hbar\bigr)
\]
See theorem~\ref{thm:DF}, where $P_k$ is the polynomial in \S\ref{sec:disc-rank-one}, $\Gamma_L$ is the Gamma function for the Levi subgroup $L$. Correcting Case~2 by this \emph{$\Gamma$-gauge}\footnote{We believe it is possible to explain this $\Gamma$-gauge via more geometric approach. The geometric Gamma class for a complex manifold $X$ is defined via the Chern roots of $T_X$ and Euler's $\Gamma$-function, which can be viewed as a square root of the Todd class of $T_X$. In the quantized world, things are often normalized by $K_X^{1/2}$, and by $\Gamma$-classes. Such $\Gamma$-functions are ubiquitous in quantum cohomology, for instance by Iritani \cite{iritani2023gammaclassesquantumcohomology}, and its $q$-analog by Okounkov \cite{okounkov2017lecturesktheoreticcomputationsenumerative}. We plan to investigate this in the future work. } makes the two isomorphisms compatible. In order to compute the discrepancy, we use localization technique, localizing to simple roots first, reducing the computation to rank one case, and then use Weyl reflection to cover all roots. Such strategy is also used in \cite{bezrukavnikov2026equivariantsatakecategorykostantwhittaker} and in many other places in geometric representation theory. Finally we use multiplicativity property, combined with integrality property, together to show that, after correction, the remaining ambiguity is given by a torus element. At $\hbar=0$, this correction is given by a  $T$-valued function, which is invisible in \cite{macerato2023leviequivariantrestrictionsphericalperverse} (Remark~\ref{rem:hbar0}), essentially because the splitting of the Whittaker cotangent bundle $T^*(G/U,\psi)$ into $G\times \mathfrak c_G$ can differ by such a function.

\subsubsection*{Organization}
\S\ref{sec:prelim}: notation, universal modules, their bimodule structures, intertwiners and their composition, normalizations (\S\ref{sec:bimodules}), the category $\calN(\psiL)$ and Backelin's functor. \S\ref{sec:localization}: the geometric side and the reduction of Theorem~\ref{thm:main} to Cases~1 and~2 (\S\S\ref{sec:case1}--\ref{sec:case2}) and the Hartogs property (\S\ref{sec:hartogs}). \S\ref{sec:rs-compat} reduces the comparison on the regular semisimple locus to the group $L$; \S\ref{sec:discrepancy} solves it for any reductive group; \S\ref{sec:gluing} concludes. Appendix~\ref{app:hbar0} treats the fiber $\hbar=0$: Case~2 is carried out after inverting $\hbar$, and the result is extended across $\hbar=0$ there. Appendix~\ref{sec:gauge-why} explains how the $\Gamma$-gauge is determined.

\subsection{Conventions}\label{sec:conventions}
We follow the notations in \cite{ginzburg2014differentialoperatorsguaffine}: On automorphic side everything is with a "check", e.g. the affine Grassmannian is that of $\cG$, $\cL$, and $\cT$; On spectral side everything is without a "check", e.g. the representations, Verma modules, and Whittaker modules are that of $G$,$L$,and $T$. \cite{bezrukavnikov2026equivariantsatakecategorykostantwhittaker,macerato2023leviequivariantrestrictionsphericalperverse} use the opposite convention, and we interchange $G$ and $\cG$ when quoting them. Moreover:
\begin{itemize}
\item All modules are left modules, and Verma modules, Whittaker modules, Verma-Whittaker modules are all quotients by left ideals. \cite{ginzburg2014differentialoperatorsguaffine} uses right modules; the two are interchanged by the antipode.
\item Whittaker characters live on $\fraku$, as in \cite{Milicic1997,MCDOWELL1985161,Backelinfunctor,eberhardt2024proofsjantzenconjecturewhittaker,macerato2023leviequivariantrestrictionsphericalperverse}; \cite{bezrukavnikov2026equivariantsatakecategorykostantwhittaker,ginzburg2022differentialoperatorsgugelfandgraev} use $\bar\fraku$. The two are interchanged by the Chevalley involution.
\item Verma modules are indexed by their highest weight, so central characters are the same on dot-orbits $w\cdot\lambda=w(\lambda+\rho)-\rho$. \cite{Milicic1997,MCDOWELL1985161,eberhardt2024proofsjantzenconjecturewhittaker} index by $\lambda+\rho$ instead.
\item \eqref{eq:BF}, transported to our conventions, is normalized as in Convention~\ref{conv:BF}, for $G$ and for its Levi subgroups.
\end{itemize}

\subsection*{Acknowledgements}
The third author is grateful to his advisor, David Nadler, for valuable conversations about Plancherel algebra that inspired this project. We thank him for many helpful discussions and suggestions. In the early stage of this work, David pointed out the necessity of showing compatibility over regular semisimple locus, which turns out to be occupying the latter half of this paper. We also thank Jeremy Taylor for sharing his ideas and his work with the third author. From conception to completion, this paper is the work of the authors, including all calculations and possible errors. The authors used open-source AI for fact-checking, recalling existing theory, and reviewing the paper.

\section{Preliminaries}\label{sec:prelim}

\subsection{Notation}\label{sec:notation}

\subsubsection*{Groups and Lie algebras}
Fix $T\subset B\subset G$ with $U$ the unipotent radical of the Borel $B = T\cdot U$. Let $L\supset T$ be a Levi subgroup, $P_L\supset B$ the standard parabolic with Levi factor $L$, $U_P$ its unipotent radical, and $U_L:=U\cap L$; then $U=U_L\cdot U_P$ with $U_P$ normal. Lie algebras are $\frakg,\frakl,\frakb,\frakt,\fraku,\frakp_L,\fraku_L,\fraku_P$, and $\bar\fraku,\bar\fraku_L,\bar\fraku_P$ are the opposite nilpotent subalgebras. Let $Z(L)$ be the center of $L$, $S = L_{\mathrm{der}}=[L,L]$, $T_S=S\cap T$, with Lie algebras $\frakz,\fraks,\tS$; so $\frakt=\tS\oplus\frakz$ and $\frakl=\fraks\oplus\frakz$.

\subsubsection*{Roots}
$R\supset R^+\supset\Delta$ are the sets of roots, positive roots and simple roots, $\check\gamma\in\frakt$ the coroot of $\gamma$, and $R_L,R_L^+,\Delta_L$ those of $L$. Throughout this paper, our convention is:
\begin{center}
\emph{$\alpha$ denotes an element of $R^+\setminus R^+_L$, and $\beta$ an element of $R^+_L$.}
\end{center}
$\rho,\rho_L$ are the half sums of $R^+,R^+_L$, and $W\supset W_L$ the Weyl groups. The dot action is $w\cdot\lambda=w(\lambda+\rho)-\rho$. Both the linear and the dot action of $W_L$ act trivially on the $\frakz^*$-component of $\frakt^*=\tS^*\oplus\frakz^*$.

\subsubsection*{Whittaker characters}
Fix a pinning $\{e_\gamma\}_{\gamma\in\Delta}$. For $\Delta'\subset\Delta$ let $\psi_{\Delta'}$ be the composite
\begin{equation}\label{eq:psi-def}
U\twoheadrightarrow U/[U,U]\cong\bigoplus_{\gamma\in\Delta}\bbG_a\twoheadrightarrow\bigoplus_{\gamma\in\Delta'}\bbG_a\xrightarrow{\ +\ }\bbG_a .
\end{equation}
Put $\psi:=\psi_\Delta$ and $\psiL:=\psi_{\Delta_L}$, and use the same letters for their differentials. Then $\psiL$ vanishes on $\fraku_P$, its restriction $\psi_{\mathrm{der}}$ to $\fraku_L$ is non-degenerate, and $\psi,\psiL$ are $\mathrm{Ad}(U)$-invariant. For $\chi\in\{\psi,\psiL,0\}$ put $\fraku^{\chi}:=\{x-\chi(x)\mid x\in\fraku\}$, and for a $\frakg$-module $E$ let $E^{\fraku,\chi}$ be the space of $\chi$-Whittaker vectors, i.e. $xv=\chi(x)v$ for $x\in\fraku$. Similarly for $\fraku_L$ and $\fraku_P$.

\subsubsection*{Centers and asymptotic algebras}
By PBW, we have:
\begin{equation}\label{eq:HC-decomp}
U(\frakg)=U(\frakt)\oplus\bigl(\bar\fraku\,U(\frakg)+U(\frakg)\,\fraku\bigr),
\end{equation}
Under this, the projection onto $U(\frakt)=\SymAlg(\frakt)$ restricts to the Harish-Chandra isomorphism $\hc\colon Z(\frakg)\xrightarrow{\sim}\SymAlg(\frakt)^{W\cdot}$; $z$ acts on the Verma module $M(\lambda)$ by $\lambda(\hc(z))$. Similarly, $Z(\frakl)=Z(\fraks)\otimes U(\frakz)\cong\SymAlg(\frakt)^{W_L\cdot}$, compatibly with the parabolic Harish-Chandra map $Z(\frakg)\to Z(\frakl)$ , which comes from the projection along $\bar\fraku_PU(\frakg)+U(\frakg)\fraku_P$. Put $\frakc_L:=\Spec Z(\frakl)$.

For a Lie algebra $\frakk$, $\Uh(\frakk)$ is the Rees algebra of $U(\frakk)$ for the PBW filtration: generated over $\bbC[\hbar]$ by $\frakk$ with $xy-yx=\hbar[x,y]$, with $\frakk$ and $\hbar$ in degree $2$ \cite[\S2.1]{ginzburg2014differentialoperatorsguaffine}. Then $\Uh(\frakk)/(\hbar-a)\cong U(\frakk)$ for $a\neq0$ and $\Uh(\frakk)/(\hbar)=\SymAlg(\frakk)$. Write $\Zh(\frakg),\Zh(\frakl)$ for the centres and $\Sh:=\Uh(\frakt)=\SymAlg(\frakt)[\hbar]$. The Harish-Chandra isomorphisms become $\Zh(\frakg)\cong\Sh^{W\cdot}$, $\Zh(\frakl)\cong\Sh^{W_L\cdot}$ for the $\hbar$-dot action $w\cdot(y,a)=(w(y+a\rho)-a\rho,a)$. All these algebras are free over $\bbC[\hbar]$. 

Following \cite{bezrukavnikov2026equivariantsatakecategorykostantwhittaker} and \cite{ginzburg2014differentialoperatorsguaffine},
for $V\in\Rep(G)$, $\Uh(\frakg)$ acts on $V$ by $\hbar$ times the action of $\frakg$.

\subsubsection*{The extended torus and the localization loci}
Let $A:=\cT\times\bbG_m$ be the extended torus. Then $\H_A(\pt)=\SymAlg(\frakt)[\hbar]=\Sh$, and we put $\fraka^*:=\Spec\Sh=\frakt^*\times\bbA^1$, with points $(y,a)$.

\begin{definition}\label{def:loci}
For $\gamma\in R$, $n\in\bbZ$ let $H_{\gamma,n}:=\{\check\gamma+n\hbar=0\}\subset\fraka^*$ be a \emph{wall}, and put
\[
\aal:=\fraka^*\setminus\bigcup_{\alpha,n}H_{\alpha,n},\qquad
\abe:=\fraka^*\setminus\bigcup_{\beta,m}H_{\beta,m},\qquad
\ars:=\aal\cap\abe ,
\]
with coordinate rings $\Cal,\Cbe,\Crs$. For an $\Sh$-module $E$ put $E_\alpha:=E\otimes_{\Sh}\Cal$, and similarly $E_\beta,E_\rs$; for a $\Zh(\frakl)$-module use the $W_L$-invariant subrings. Let $\frakt^*_\alpha,\frakt^*_\beta,\frakt^*_\rs$ be the fibers over $a=1$.
\end{definition}

These sets are $W_L$-stable. It is worth noting that all walls of \emph{all} roots of the relevant type must be removed, not only just the simple ones; this is what the localization theorem requires, compare \cite[Proposition 2.1.6]{macerato2023leviequivariantrestrictionsphericalperverse}.

\begin{remark}\label{rem:specialization}
For $a\neq0$, specializing a $(\Uh(\frakg),\Sh)$-bimodule at $(y,a)$ gives a $U(\frakg)$-module on which $\frakt$ acts through $\lambda=y/a$; e.g.\ $\bM_\hbar$ specializes to $M(y/a)$. Since $\langle\rho,\check\gamma\rangle\in\bbZ$, the walls $H_{\gamma,n}$ are the asymptotic versions of $\{\langle\lambda,\check\gamma\rangle\in\bbZ\}=\{\langle\lambda+\rho,\check\gamma\rangle\in\bbZ\}$. More generally, the set of walls is stable under $y\mapsto-y$ and under translation by $\hbar\mu$ whenever $\langle\mu,\check\gamma\rangle\in\bbZ$ for all $\gamma$, e.g.\ $\mu\in\Lambda$ or $\mu=\rho$. At $a=0$ all $H_{\gamma,n}$ become $\{\check\gamma=0\}$.
\end{remark}

\subsection{Universal modules}\label{sec:universal}
For simplicity, we work with $U(\frakg)$; the asymptotic versions, with $\Uh$ in place of $U$, are denoted by a subscript $\hbar$. Their bimodule structures and the intertwiners between them are discussed in \S\ref{sec:bimodules}.

\subsubsection{Universal Whittaker modules}\label{sec:whittaker}
The universal Whittaker module $\calW:=U(\frakg)/U(\frakg)\fraku^{\psi}$ is a $(U(\frakg),Z(\frakg))$-bimodule. By Kostant \cite{Kostant1978} (see also \cite[\S4.2]{bezrukavnikov2026equivariantsatakecategorykostantwhittaker}):
\begin{enumerate}[label=(\roman*)]
\item $\fraku^\psi$ acts locally nilpotently on $\calW$;
\item $Z(\frakg)\xrightarrow{\sim}\calW^{\fraku,\psi}$, $z\mapsto z\cdot1$; equivalently $\End_{U(\frakg)}(\calW)=Z(\frakg)$;
\item $\calW$ is free over $Z(\frakg)$;
\item each specialization $\calW(\xi)=\calW\otimes_{Z(\frakg)}\bbC_\xi$ is irreducible with $\calW(\xi)^{\fraku,\psi}=\bbC\cdot1$, and every irreducible module generated by a $\psi$-Whittaker vector is of this form.
\end{enumerate}
The same holds for $\calW_\hbar$ over $\Zh(\frakg)$. Therefore, \eqref{eq:BF} can also be read as: $\H_{\GO}(\GrG,\IC_V)\cong\Hom_{\Uh(\frakg)}(\calW_\hbar,V\otimes\calW_\hbar)$.

\begin{remark}[Group and Lie algebra invariants]\label{rem:U-inv}
If $\fraku^\chi$ acts locally nilpotently on $E$, the action $x\mapsto x-\chi(x)$ integrates to an algebraic $U$-action, such that $E^U=E^{\fraku,\chi}$. On $\calW$ and on $\bM(\psiL)$ below this is the adjoint action. Hence $(V\otimes\calW_\hbar)^{U}=(V\otimes\calW_\hbar)^{\fraku,\psi}$ and $(V\otimes\bM_\hbar(\psiL))^U=(V\otimes\bM_\hbar(\psiL))^{\fraku,\psiL}$, which is the sense of the invariants in \eqref{eq:BF-reg} and \eqref{eq:main}.
\end{remark}

\subsubsection{Universal Verma modules}\label{sec:verma}
The universal Verma module $\bM:=U(\frakg)/U(\frakg)\fraku=U(\frakg)\otimes_{U(\frakb)}\SymAlg(\frakt)$ is a $(U(\frakg),\SymAlg(\frakt))$-bimodule, which is free over $U(\bar\fraku)\otimes\SymAlg(\frakt)$. Its specialization at $\lambda\in\frakt^*$ is the usual Verma module $M(\lambda)$, which is simple iff $\langle\lambda+\rho,\check\gamma\rangle\notin\bbZ_{>0}$ for all $\gamma\in R^+$ (\S\ref{sec:shapovalov}). The right action of $\Sh$ on $\bM_\hbar$ is called the \emph{Verma parameter}, with coordinate $y$ (\S\ref{sec:right-actions}). Let $\one$ denote the generator. $\bM_L,\bM_{L,\hbar}$ denote the universal Verma modules of $\frakl$.

\subsubsection{Universal Verma--Whittaker modules}\label{sec:verma-whittaker}
The \emph{universal Verma--Whittaker module} is
\begin{equation}\label{eq:MpsiL-def}
\bM(\psiL):=U(\frakg)/U(\frakg)\fraku^{\psiL}=U(\frakg)\otimes_{U(\frakp_L)}\calW_L,\qquad \calW_L:=U(\frakl)/U(\frakl)\fraku_L^{\psiL},
\end{equation}
The second equality above is because $\psiL|_{\fraku_P}=0$. Thus it is parabolically induced from the universal Whittaker module of $\frakl$. This module has already appeared in \cite[\S4.4]{bezrukavnikov2026equivariantsatakecategorykostantwhittaker}, and was called the universal Verma-Whittaker module there.\footnote{Under the Chevalley involution, $\bM_\hbar(\psiL)$ corresponds to $\mathcal{MW}_L(-\rho+\rho_L)$ of loc.\ cit.} The right action of $Z(\frakl)$ on $\calW_L$ makes $\bM(\psiL)$ a $(U(\frakg),Z(\frakl))$-bimodule, which is free over $Z(\frakl)$ since $\bM(\psiL)\cong U(\bar\fraku_P)\otimes\calW_L$. For $L=G$ it is $\calW$, for $L=T$ it is $\bM$.

For $\lambda\in\frakt^*$ let $\Omega_\lambda$ be the central character of $M_L(\lambda)$, i.e.\ $\lambda\circ\hc$ on $Z(\frakl)$, and put
\begin{equation}\label{eq:M-psiL-lambda}
\calW_L(\lambda):=\calW_L\otimes_{Z(\frakl)}\bbC_{\Omega_\lambda},\qquad M(\psiL,\lambda):=\bM(\psiL)\otimes_{Z(\frakl)}\bbC_{\Omega_\lambda}=U(\frakg)\otimes_{U(\frakp_L)}\calW_L(\lambda).
\end{equation}
Every character of $Z(\frakl)$ is some $\Omega_\lambda$, and $\Omega_\lambda=\Omega_\mu$ iff $W_L\cdot\lambda=W_L\cdot\mu$. The asymptotic version $\bM_\hbar(\psiL)=\Uh(\frakg)\otimes_{\Uh(\frakp_L)}\calW_{L,\hbar}$ is a $(\Uh(\frakg),\Zh(\frakl))$-bimodule, which is free over $\Zh(\frakl)$.

\begin{lemma}[Kostant theory for $\fraks$]\label{lem:Kostant-for-s}
Let $\calC$ be the category of $\fraks$-modules on which $\fraku_L^{\psiL}$ acts locally nilpotently, and $\mathrm{Wh}(N):=N^{\fraku_L,\psiL}$.
\begin{enumerate}
\item[(K1)] $\mathrm{Wh}$ is an equivalence between $\calC$ and $Z(\fraks)$-modules \cite[\S\S3--4]{Kostant1978}, \cite{skryabin2002}. In particular $\mathrm{Wh}$ is exact, and $\mathrm{Wh}(N)=0$ implies $N=0$.
\item[(K2)] For $F$ a finite-dimensional $\frakl$-module, $F\otimes\calW_L(\lambda)\in\calC$ and $\dim\mathrm{Wh}(F\otimes\calW_L(\lambda))=\dim F$. Indeed it has a filtration with $\dim F$ subquotients $\calW_L(\lambda_i)$ \cite[Theorem 4.6]{Kostant1978}, each with one-dimensional Whittaker space.
\end{enumerate}
The same holds for $(\frakg,\psi)$ and for the asymptotic versions.
\end{lemma}

\subsubsection{A lemma of Ginzburg--Kazhdan}\label{sec:GK}

The following lemma will be used frequently throught this paper:
\begin{lemma}[{\cite[Lemma 6.2.1]{ginzburg2022differentialoperatorsgugelfandgraev}}]\label{lem:GK}
Let $\psi$ be non-degenerate and $E$ a $\Uh(\frakg)$-module on which $\fraku^\psi$ acts locally nilpotently. Then
\[
\Sh\otimes_{\Zh(\frakg)}E^{\fraku,\psi}\longrightarrow\bar\fraku\backslash E:=E/\bar\fraku E,\qquad s\otimes v\mapsto s\cdot[v],
\]
is an isomorphism of $\Sh$-modules, natural in $E$.
\end{lemma}

This is loc.\ cit.\ transported by the Chevalley involution; the asymptotic version follows by flatness over $\bbC[\hbar]$. Note that the coinvariants are for the \emph{untwisted} $\bar\fraku$. With $E=V\otimes\calW_\hbar$, \eqref{eq:BF} and equivariant formality it gives \eqref{eq:BF-torus}.

\subsection{Bimodule structures and intertwiners}\label{sec:bimodules}
The universal modules carry two commuting structures, and the intertwiners between them can be composed. Both are used throughout the paper, so we fix the terminology here. All modules are left modules; a \emph{right action} is an action by endomorphisms of the left module. We work in the asymptotic setting; setting $\hbar=1$ gives the non-asymptotic versions.

\subsubsection{Right actions}\label{sec:right-actions}
Each universal module $X$ is a left $\Uh(\frakg)$-module with a commuting right action of a commutative ring $\calR_X\subset\Uh(\frakg)$, such that $\one\cdot s=s\cdot\one$ for $s\in\calR_X$:
\begin{enumerate}[label=(\alph*)]
\item $X=\calW_\hbar$, $\calR_X=\Zh(\frakg)$, by right multiplication; well defined since $\Zh(\frakg)$ is central. On $\calW_\hbar$ itself the left and right actions of $\Zh(\frakg)$ coincide.
\item $X=\bM_\hbar$, $\calR_X=\Sh=\Uh(\frakt)$, by right multiplication; well defined since $\frakt$ normalizes $\fraku$. This right action is the \emph{Verma parameter}, and $y$ always denotes its coordinate: a function $s(y)$ acts by $m\mapsto m\cdot s$. The left action of $\Zh(\frakg)$ factors through it: $z\cdot m=m\cdot\hc(z)$, because $z$ has weight $0$, hence $z\in\Uh(\frakt)\oplus\Uh(\frakg)\fraku$ by PBW.
\item $X=\bM_\hbar(\psiL)$, $\calR_X=\Zh(\frakl)\subset\Uh(\frakl)$, induced from the right action on $\calW_{L,\hbar}$, which commutes with $\frakp_L$. As in (b), $z\cdot m=m\cdot\hc_P(z)$, with $\hc_P\colon\Zh(\frakg)\to\Zh(\frakl)$ the parabolic Harish-Chandra map.
\end{enumerate}

\subsubsection{Left and right actions on tensor products}\label{sec:left-action}
Let $V\in\Rep(G)$ and let $X$ be as above. On $V\otimes X$ we use exactly two actions of $\calR_X$:
\begin{itemize}
\item the \emph{right action} $(v\otimes m)\cdot s:=v\otimes ms$, written on the right;
\item the \emph{left action} $s\cdot(v\otimes m)$, $s\in\calR_X\subset\Uh(\frakg)$, the tensor product module structure, written on the left.
\end{itemize}
Recall (\S\ref{sec:notation}) that $\Uh(\frakg)$ acts on $V$ by $\hbar$ times the action of $\frakg$. So the tensor product module structure on $V\otimes X$ is
\[
x\cdot(v\otimes m)=\hbar\,xv\otimes m+v\otimes xm\qquad(x\in\frakg),
\]
where $xv$ is the classical action, and it is extended multiplicatively to $\Uh(\frakg)$. Equivalently, it is the coproduct $\Delta(x)=x\otimes1+1\otimes x$ of $\Uh(\frakg)$ followed by the action on $V$ just recalled; with the classical action on $V$ it would read $\Delta_\hbar(x)=\hbar\,x\otimes1+1\otimes x$.
They commute. The left action is often called the \emph{diagonal action}. On intertwiners the left action becomes precomposition, and the right action becomes postcomposition (\S\ref{sec:intertwiners}). For $s\notin\frakg$ the left action is not $s\otimes1+1\otimes s$; for instance, for $c\in\Sh$ and $v\in V_\mu$ we have $c\cdot(v\otimes m)=v\otimes c(\,\cdot+\hbar\mu)\,m$.

\paragraph{Relative weights.}
On $V\otimes\bM_\hbar$ the difference of the two actions of $h\in\frakt$,
\[
\ad(h)(v\otimes u\one):=h\cdot(v\otimes u\one)-(v\otimes u\one)\cdot h=hv\otimes u\one+v\otimes[h,u]\one ,
\]
is the tensor product of the action on $V$ and the adjoint action on $\bM_\hbar$. It is semisimple, commutes with both actions of $\Sh$, and preserves $(V\otimes\bM_\hbar)^{\fraku}$. Its eigenvalue on $v\otimes u\one$, $v\in V_\mu$, $u$ of weight $-\gamma$, is $\hbar(\mu-\gamma)$, and we call $\mu-\gamma$ the \emph{relative weight}. Hence $(V\otimes\bM_\hbar)^{\fraku}=\bigoplus_\nu(V\otimes\bM_\hbar)^{\fraku}_\nu$, and on the $\nu$-piece
\begin{equation}\label{eq:left-right}
s\cdot\xi=\xi\cdot s(y+\hbar\nu),\qquad s\in\Sh,\ \xi\in(V\otimes\bM_\hbar)^{\fraku}_\nu ,
\end{equation}
i.e.\ \emph{the left action is the right action shifted by $\hbar\nu$}. The same holds on the completion $V\otimes\widehat\bM_\hbar$. Both actions descend to $\bar\fraku\backslash(V\otimes\bM_\hbar)$ because $\frakt$ normalizes $\bar\fraku$, and there $s\cdot[v\otimes\one]=[v\otimes\one]\cdot s(y+\hbar\nu)$ for $v\in V_\nu$. On $V\otimes\bM_\hbar(\psiL)$ one uses $\frakz$ in place of $\frakt$.

\subsubsection{Intertwiners and their composition}\label{sec:intertwiners}
Let $\chi\in\{\psi,0,\psiL\}$ be the character attached to $X\in\{\calW_\hbar,\bM_\hbar,\bM_\hbar(\psiL)\}$. By Frobenius reciprocity
\begin{equation}\label{eq:frobenius}
(V\otimes X)^{\fraku,\chi}=\Hom_{\Uh(\frakg)}(X,\,V\otimes X),\qquad\xi\longmapsto f_\xi,\quad f_\xi(\one)=\xi .
\end{equation}
Write $R_s$ for right multiplication by $s\in\calR_X$ on $X$. Then
\begin{equation}\label{eq:pre-post}
f_{s\cdot\xi}=f_\xi\circ R_s,\qquad f_{\xi\cdot s}=(\id_V\otimes R_s)\circ f_\xi ,
\end{equation}
since both sides are $\Uh(\frakg)$-linear and agree on $\one$; for the first, $f_\xi(\one\cdot s)=f_\xi(s\cdot\one)=s\cdot\xi$. In particular both actions preserve $(V\otimes X)^{\fraku,\chi}$: \emph{the left action is precomposition, and the right action is postcomposition}. The maps $\Theta$ (\S\ref{sec:case1}) and $\Phi$ (\S\ref{sec:case2}) commute with the right action and with the left action of $\Zh(\frakg)$.

\begin{definition}[Composition]\label{def:composition}
For $\xi\in(V_1\otimes X)^{\fraku,\chi}$ and $\xi'\in(V_2\otimes X)^{\fraku,\chi}$ put
\[
\xi'\circ\xi:=(\id_{V_2}\otimes f_\xi)(\xi')\ \in\ (V_2\otimes V_1\otimes X)^{\fraku,\chi}.
\]
\end{definition}

\begin{lemma}[Composition of intertwiners]\label{lem:composition}
Let $\xi,\xi'$ be as in Definition~\ref{def:composition} and $a\in\calR_X$.
\begin{enumerate}[label=(\alph*)]
\item $f_{\xi'\circ\xi}=(\id_{V_2}\otimes f_\xi)\circ f_{\xi'}$. From this we see $\circ$ is associative, natural in $V_1,V_2$, with unit $\one\in(\mathrm{triv}\otimes X)^{\fraku,\chi}$.
\item We have the following:
\[
a\cdot(\xi'\circ\xi)=(a\cdot\xi')\circ\xi,\qquad(\xi'\cdot a)\circ\xi=\xi'\circ(a\cdot\xi),\qquad(\xi'\circ\xi)\cdot a=\xi'\circ(\xi\cdot a).
\]
\end{enumerate}
\end{lemma}

\begin{proof}
(a) Both sides are $\Uh(\frakg)$-linear, and they both send $\one$ to $(\id\otimes f_\xi)(\xi')$; associativity follows. 

(b) We use (a) \eqref{eq:pre-post}. The first equality: $a\cdot(\xi'\circ\xi)=f_{\xi'\circ\xi}(\one\cdot a)=(\id\otimes f_\xi)(a\cdot\xi')$. The next equality: $(\id\otimes f_\xi)(\xi'\cdot a)=(\id\otimes f_\xi\circ R_a)(\xi')=(\id\otimes f_{a\cdot\xi})(\xi')$. The last equality: $(\id\otimes f_\xi)(\xi')\cdot a=(\id\otimes(\id\otimes R_a)\circ f_\xi)(\xi')=(\id\otimes f_{\xi\cdot a})(\xi')$.
\end{proof}

By Lemma~\ref{lem:composition}, composition makes $\bigoplus_V V^*\otimes(V\otimes X)^{\fraku,\chi}$ into an algebra. It is $\Dh(G/U,\psi)$ for $X=\calW_\hbar$, $\Dh(G/U)$ for $X=\bM_\hbar$ and $\Dh(G/U,\psiL)$ for $X=\bM_\hbar(\psiL)$ (Remark~\ref{rem:U-inv}). In other words, we have understood the ring structure of these ring of (twisted) differential operators via the composition of intertwiners.

For $X=\bM_\hbar$, combining Lemma~\ref{lem:composition} with \eqref{eq:left-right} gives the formula used throughout \S\S\ref{sec:case2}--\ref{sec:discrepancy}. We denote $\sfA(V):=(V\otimes\bM_\hbar)^{\fraku}=\bigoplus_\nu\sfA_\nu(V)$.

\begin{corollary}[Composition and the Verma parameter]\label{cor:comp-verma}
Let $\xi\in\sfA_\nu(V_1)$, $\xi'\in\sfA_\mu(V_2)$, and let $a,b$ lie in $\Sh$ or in one of its localizations along walls (Definition~\ref{def:loci}).
\begin{enumerate}[label=(\alph*)]
\item $f_\xi(\one\cdot a)=a\cdot\xi=\xi\cdot a(y+\hbar\nu)$.
\item $\xi'\circ\xi\in\sfA_{\mu+\nu}(V_2\otimes V_1)$, and
\begin{equation}\label{eq:comp-shift}
(\xi'\cdot a)\circ(\xi\cdot b)=(\xi'\circ\xi)\cdot a(y+\hbar\nu)\,b(y).
\end{equation}
In other words: specialized at Verma parameter $y$, $f_\xi$ is an intertwiner $M(y+\hbar\nu)\to V_1\otimes M(y)$, and in the composite
\[
M\bigl(y+\hbar(\mu+\nu)\bigr)\xrightarrow{\ f_{\xi'}\ }V_2\otimes M(y+\hbar\nu)\xrightarrow{\ \id\otimes f_\xi\ }V_2\otimes V_1\otimes M(y)
\]
the coefficients of $\xi'$ are evaluated at the parameter $y+\hbar\nu$ of the middle term.
\item (The gauges). Let $c=(c_\nu)_\nu$ be a family of units of a localization of $\Sh$ along walls, acting on $\sfA$ by $\xi\mapsto\xi\cdot c_\nu$ on $\sfA_\nu$. If $c$ satisfies the cocycle condition
\begin{equation}\label{eq:gauge-cocycle}
c_{\mu+\nu}(y)=c_\mu(y+\hbar\nu)\,c_\nu(y)\qquad\text{for all weights }\mu,\nu,
\end{equation}
then $c$ is multiplicative: $c(\xi'\circ\xi)=c(\xi')\circ c(\xi)$. Examples are characters $c_\nu=\nu(t)$, $t\in T(\bbC)$, and coboundaries $c_\nu(y)=\varphi(y)/\varphi(y+\hbar\nu)$ for a (possibly transcendental) function $\varphi$ of $y$ for which these ratios are units of the localization. A coboundary is conjugation by right multiplication by $\varphi$: $f_{c(\xi)}=R_\varphi\circ f_\xi\circ R_\varphi^{-1}$.
\end{enumerate}
\end{corollary}

\begin{proof}
(a) is \eqref{eq:pre-post} together with \eqref{eq:left-right}. (b) Let $h\in\frakt$. By Lemma~\ref{lem:composition}(b) and (a),
\[
h\cdot(\xi'\circ\xi)=(h\cdot\xi')\circ\xi=(\xi'\cdot h)\circ\xi+\hbar\mu(h)\,\xi'\circ\xi=
\]
\[
\xi'\circ(h\cdot\xi)+\hbar\mu(h)\,\xi'\circ\xi=(\xi'\circ\xi)\cdot h+\hbar(\mu+\nu)(h)\,\xi'\circ\xi ,
\]
From this we see $\xi'\circ\xi$ has relative weight $\mu+\nu$. Similarly, $(\xi'\cdot a)\circ(\xi\cdot b)=\xi'\circ\bigl(a\cdot(\xi\cdot b)\bigr)=\xi'\circ\bigl(\xi\cdot a(y+\hbar\nu)b(y)\bigr)=(\xi'\circ\xi)\cdot a(y+\hbar\nu)b(y)$. 

(c) By \eqref{eq:comp-shift}, we compute: $c(\xi')\circ c(\xi)=(\xi'\circ\xi)\cdot c_\mu(y+\hbar\nu)c_\nu(y)$, and it equals $c(\xi'\circ\xi)$ by \eqref{eq:gauge-cocycle}. Both families satisfy \eqref{eq:gauge-cocycle} by direct computation. For the coboundary, by (a), we have:
\[
R_\varphi f_\xi R_\varphi^{-1}(\one)=R_\varphi\bigl(\varphi^{-1}\cdot\xi\bigr)=R_\varphi\bigl(\xi\cdot\varphi(y+\hbar\nu)^{-1}\bigr)=\xi\cdot\varphi(y)/\varphi(y+\hbar\nu)
\]\end{proof}

\subsection{The category \texorpdfstring{$\calN(\psiL)$}{N(psi\_L)}}\label{sec:N}
In \cite{Milicic1997}, Mili\v{c}i\'c--Soergel considered the category $\calN$ of $U(\frakg)$-modules that are finitely generated and locally finite over $U(\fraku)$ and over $Z(\frakg)$. It decomposes as $\bigoplus_\eta\calN(\eta)$ over characters $\eta$ of $\fraku$, $\calN(\eta)$ consisting of modules on which $\fraku^\eta$ acts locally nilpotently. $\calN(0)$ is category $\calO$ (without $\frakt$-semisimplicity), $\calN(\psi)$ consists of Whittaker modules, and $M(\psiL,\lambda)\in\calN(\psiL)$.

\subsubsection{Basic properties}\label{sec:N-properties}
The following lemmas are due to McDowell \cite{MCDOWELL1985161} and Mili\v{c}i\'c--Soergel \cite{Milicic1997}; see also \cite[\S2.2]{eberhardt2024proofsjantzenconjecturewhittaker}.

\begin{lemma}[{\cite[Proposition 2.1]{Milicic1997}}]\label{lem:MS}
$M(\psiL,\lambda)\cong M(\psiL,\mu)$ iff $W_L\cdot\lambda=W_L\cdot\mu$. $M(\psiL,\lambda)$ has a unique simple quotient $L(\psiL,\lambda)$, these are pairwise non-isomorphic for distinct $W_L$-dot-orbits, and they exhaust the simple objects of $\calN(\psiL)$.
\end{lemma}

We put an order on $\frakz^*$, defined by $\mu\leq\omega$ iff $\omega-\mu\in\sum_\alpha\bbZ_{\geq0}\,\alpha|_\frakz$.

\begin{lemma}[{\cite[Proposition 2.4]
{MCDOWELL1985161}}]\label{lem:McD-weights}
$\frakz$ acts semisimply on $M(\psiL,\lambda)$, and with $\omega=\lambda|_\frakz$,
\[
M(\psiL,\lambda)=\bigoplus_{\gamma\geq0}M(\psiL,\lambda)^{\omega-\gamma},\qquad M(\psiL,\lambda)^{\omega-\gamma}\cong U(\bar\fraku_P)^{-\gamma}\otimes\calW_L(\lambda)
\]
as $\frakl$-modules; in particular the top piece is $\calW_L(\lambda)$.
\end{lemma}

A \emph{$\psiL$-vector} is an element of $M^{\fraku,\psiL}$; a non-zero $v$ is called \emph{maximal} if $\fraku_Pv=0$.

\begin{lemma}[{\cite[Lemma 2.12]{MCDOWELL1985161}}]\label{lem:McD-max}
Let $M\neq0$ be $M(\psiL,\lambda)$ or a submodule of it, and $\omega=\lambda|_\frakz$.
\begin{enumerate}[label=(\arabic*)]
\item $M^{\omega'}$ contains a maximal vector iff it contains a non-zero $\psiL$-vector.
\item $M$ is reducible iff $M^{\omega'}$ contains a maximal vector for some $\omega'<\omega$.
\end{enumerate}
\end{lemma}

So $M(\psiL,\lambda)^{\fraku,\psiL}$ is $\frakz$-graded with top part $\bbC\cdot1$, and $M(\psiL,\lambda)$ is simple iff it has no $\psiL$-vectors of weight $<\omega$.

\subsubsection{Backelin's functor}\label{sec:backelin}
For a $\frakt$-semisimple module $M=\bigoplus_\mu M_\mu$, we denote its completion by $\widehat M:=\prod_\mu M_\mu$. For any $\frakg$-module $N$, let $\Gamma_{\psiL}(N)$ be the submodule of vectors killed by some power of every $x-\psiL(x)$, $x\in\fraku$ \cite[\S6]{Backelinfunctor}. Backelin's functor is
\begin{equation}\label{eq:backelin}
\widehat\Gamma_{\psiL}\colon\calO\longrightarrow\calN(\psiL),\qquad M\longmapsto\Gamma_{\psiL}(\widehat M).
\end{equation}
It commutes with $E\otimes(-)$ for $E$ finite-dimensional, and with parabolic induction: $\widehat\Gamma_{\psiL}(U(\frakg)\otimes_{U(\frakp_L)}N)\cong U(\frakg)\otimes_{U(\frakp_L)}\widehat\Gamma_{\psiL}(N)$ for $N$ in the category $\calO$ of $\frakl$ \cite[Remark 2.5.5]{eberhardt2024proofsjantzenconjecturewhittaker}.

\begin{lemma}[{\cite[Proposition 6.9]{Backelinfunctor}, \cite[Theorem 2.5.4]{eberhardt2024proofsjantzenconjecturewhittaker}}]\label{lem:backelin-verma}
If the Verma module $M_L(\lambda)$ of $\frakl$ is simple, then $\widehat\Gamma_{\psiL}(M(\lambda))\cong M(\psiL,\lambda)$.
\end{lemma}

For $\frakl$ this is $\widehat\Gamma_{\psiL}(M_L(\lambda))\cong\calW_L(\lambda)$ \cite[Lemmas 6.5, 6.6]{Backelinfunctor}; the general case follows by parabolic induction.

\begin{remark}\label{rem:tensor-filtration}
For $E$ finite-dimensional with weights $\nu_1,\dots,\nu_r$, $E\otimes M(\psiL,\lambda)\cong U(\frakg)\otimes_{U(\frakp_L)}(E|_{\frakp_L}\otimes\calW_L(\lambda))$ has a filtration with subquotients $M(\psiL,\lambda+\nu_i)$ (cf. \cite[\S3]{Kostant1978}, \cite[\S2]{Milicic1997}). Hence $\widehat\Gamma_{\psiL}(E\otimes M(\lambda))\cong E\otimes M(\psiL,\lambda)$ whenever $M_L(\lambda)$ is simple.
\end{remark}

\section{Equivariant localization}\label{sec:localization}

This whole section focus on the geometric side, i.e. we will discuss the $A$-equivariant localization on $\GrG$.
Following \cite{macerato2023leviequivariantrestrictionsphericalperverse}, we localize in two groups of roots: the $\alpha$-walls (Case~1, over $\aal$) and the $\beta$-walls (Case~2, over $\abe$).

\subsection{Localization on the geometric side}\label{sec:geom-loc}
Fix $\calF\in\Perv_{\cG(\calO)}(\GrG)$ and $V=\bbS_G(\calF)$. The fixed points of $A$ on $\GrG$ are $\GrT=\Lambda$, and $A$ preserves $\GrL$. Let $i_T^L\colon\GrT\hookrightarrow\GrL$, so $i_T=i_L\circ i_T^L$. Let $\Res^{\cG}_{\cL}$ be hyperbolic restriction, normalized as in \cite[\S\S1.2, 2.4]{macerato2023leviequivariantrestrictionsphericalperverse}, so that $\Res^{\cG}_{\cL}\IC_V\cong\bbS_L^{-1}(V|_L)$ \cite[Theorem 4.1.1]{macerato2023leviequivariantrestrictionsphericalperverse}, and let $\Xi\colon i_L^!\to\Res^{\cG}_{\cL}$ be the natural transformation of \cite[(1.3.1)]{macerato2023leviequivariantrestrictionsphericalperverse}.

\begin{lemma}[{\cite[Proposition 2.1.1]{macerato2023leviequivariantrestrictionsphericalperverse}}]\label{lem:formal}
$i_L^!\calF$ is $A$-equivariantly formal. Hence $\H_A(\GrL,i_L^!\calF)$ is free over $\Sh$ of rank $\dim V$, and
\[
\H_{\LO}(\GrL,i_L^!\calF)\otimes_{\Zh(\frakl)}\Sh\cong\H_A(\GrL,i_L^!\calF),\qquad \H_{\LO}(\GrL,i_L^!\calF)=\H_A(\GrL,i_L^!\calF)^{W_L}.
\]
\end{lemma}

\begin{proof}
Loc.\ cit.\ proves this for $\cT$ by a parity argument; since $\H_A(\pt)$ is also even, it applies to $A$. Using that $\H_{\LO}(\pt)=\Sh^{W_L}=\Zh(\frakl)$, the result follows. The rank is computed after localization, by Proposition~\ref{prop:geom-loc}(3).
\end{proof}

\begin{proposition}\label{prop:geom-loc}
Let $\calF\in\Perv_{\cG(\calO)}(\GrG)$.
\begin{enumerate}[label=(\arabic*)]
\item (Case 1.) $\H_A(\Xi)\colon\H_A(\GrL,i_L^!\calF)\to\H_A(\GrL,\Res^{\cG}_{\cL}\calF)$ is injective, and it becomes an isomorphism after removing all the $\alpha$-walls:
\begin{equation}\label{eq:case1-geom}
\H_A(\GrL,i_L^!\calF)_\alpha\xrightarrow{\ \sim\ }\H_A(\GrL,\Res^{\cG}_{\cL}\calF)_\alpha .
\end{equation}
\item (Case 2.) The counit $\H_A(\GrT,i_T^!\calF)\to\H_A(\GrL,i_L^!\calF)$ is injective, and it becomes an isomorphism after removing all the $\beta$-walls:
\begin{equation}\label{eq:case2-geom}
\H_A(\GrT,i_T^!\calF)_\beta\xrightarrow{\ \sim\ }\H_A(\GrL,i_L^!\calF)_\beta .
\end{equation}
\item (Regular semisimple locus.) After removing all walls,
\begin{equation}\label{eq:rs-geom}
\begin{gathered}
\H_A(\GrT,i_T^!\calF)_\rs\xrightarrow{\ \sim\ }\H_A(\GrL,i_L^!\calF)_\rs\xrightarrow{\ \sim\ }\H_A(\GrL,\Res^{\cG}_{\cL}\calF)_\rs,\\
\H_A(\GrT,i_T^!\calF)_\rs\xrightarrow{\ \sim\ }\H_A(\GrG,\calF)_\rs ,
\end{gathered}
\end{equation}
and all these modules are identified with $\bigoplus_{\lambda}V_\lambda\otimes\Crs$ through the Mirkovi\'c--Vilonen weight functors, cf \cite{mirkovic2018geometriclanglandsdualityrepresentations}.
\end{enumerate}
\end{proposition}

\begin{proof}
Without loop rotation this is \cite[Propositions 2.1.6, 2.6.7]{macerato2023leviequivariantrestrictionsphericalperverse}. The proofs carry over to $A$, with one minor change in \cite[Proposition 2.1.5]{macerato2023leviequivariantrestrictionsphericalperverse}: a point $x\in\GrG\setminus\GrL$ has  non-trivial coordinates in some $\cU_{\alpha}(t^{-1}\bbC[t^{-1}])$, and $A$ acts on the coefficient of $t^{-n}$ by $\check\alpha-n\hbar$; so some $\check\alpha-n\hbar$ is trivial on the stabilizer $A_x$. Only finitely many $n$ occur, and the localization theorem \cite[Theorem A.1.13]{macerato2023leviequivariantrestrictionsphericalperverse} applies. Injectivity follows from Lemma~\ref{lem:formal}. Part (3) is the case $L=T$ for $\cG$ and $\cL$, together with \cite[\S3.5]{bezrukavnikov2026equivariantsatakecategorykostantwhittaker}.
\end{proof}

\subsection{Gluing on the geometric side}\label{sec:geom-glue}
The gluing on the geometric side comes for free:
\begin{proposition}\label{prop:geom-glue}
The natural map
\[
\H_A(\GrL,i_L^!\calF)\longrightarrow\H_A(\GrL,i_L^!\calF)_\alpha\times_{\H_A(\GrL,i_L^!\calF)_\rs}\H_A(\GrL,i_L^!\calF)_\beta
\]
is an isomorphism; likewise for $\H_{\LO}$ over $\Zh(\frakl)$.
\end{proposition}

\begin{proof}
The module is free (Lemma~\ref{lem:formal}), and $\Cal\cap\Cbe=\Sh$ in $\Crs$ since $\Sh$ is a UFD and no $\alpha$-wall is a $\beta$-wall. Take $W_L$-invariants for the second claim.
\end{proof}

By Lemma~\ref{lem:formal} and Remark~\ref{rem:U-inv}, Theorem~\ref{thm:main} reduces to constructing, for each $V\in\Rep(G)$, $W_L$-equivariant isomorphisms
\begin{align}
\widetilde\Psi_{V,\alpha}&\colon\ (V\otimes\bM_\hbar(\psiL))^{\fraku,\psiL}\otimes_{\Zh(\frakl)}\Cal\xrightarrow{\ \sim\ }\H_A(\GrL,i_L^!\IC_V)_\alpha, \label{eq:PsiA}\\
\widetilde\Psi_{V,\beta}&\colon\ (V\otimes\bM_\hbar(\psiL))^{\fraku,\psiL}\otimes_{\Zh(\frakl)}\Cbe\xrightarrow{\ \sim\ }\H_A(\GrL,i_L^!\IC_V)_\beta, \label{eq:PsiB}
\end{align}
which agree over $\ars$, and to proving the algebraic analogue of Proposition~\ref{prop:geom-glue}:
\begin{equation}\label{eq:alg-glue}
\Dh(G/U,\psiL)\xrightarrow{\ \sim\ }\Dh(G/U,\psiL)_\alpha\times_{\Dh(G/U,\psiL)_\rs}\Dh(G/U,\psiL)_\beta .
\end{equation}
We construct $\widetilde\Psi_{V,\alpha}$ in \S\ref{sec:case1}, then a possible candidate for $\widetilde\Psi_{V,\beta}$ in \S\ref{sec:case2}, and prove \eqref{eq:alg-glue} in \S\ref{sec:hartogs}. Unforunately, the possible candidate does not agree with $\widetilde\Psi_{V,\alpha}$ over $\ars$; We correct it in \S\S\ref{sec:rs-compat}--\ref{sec:gluing} .

\section{Case 1: removing the \texorpdfstring{$\alpha$}{alpha}-walls}\label{sec:case1}
We work over $\Cal$ throughout this section. Geometrically, Proposition~\ref{prop:geom-loc}(1) reduces to equivariant Satake for $\cL$. Algebraically, we show that after inverting the $\alpha$-walls $\bM_\hbar(\psiL)$ behaves like $\calW_{L,\hbar}$, as far as only the Whittaker vectors in these modules are concerned.

\subsection{Simplicity of standard modules off the \texorpdfstring{$\alpha$}{alpha}-walls}\label{sec:simple}

We first work with non-asymtopic version.
\begin{lemma}\label{lem:submodule}
If $M(\psiL,\lambda)$ has a proper non-zero submodule, then $\langle\lambda+\rho,\check\alpha\rangle\in\bbZ$ for some $\alpha\in R^+\setminus R^+_L$.
\end{lemma}

\begin{proof}
By Lemma~\ref{lem:McD-weights}, $M(\psiL,\lambda)\cong U(\bar\fraku_P)\otimes\calW_L(\lambda)$ decomposes into $\frakz$-weight spaces of weights $\omega-\gamma$, where $\omega=\lambda|_\frakz$ and $\gamma$ runs over $\bbZ_{\geq0}$-combinations of the restrictions $\alpha|_\frakz$, $\alpha\in R^+\setminus R^+_L$. Suppose that $M(\psiL,\lambda)$ is not simple. By Lemma~\ref{lem:McD-max} there is a non-zero $\psiL$-vector $v\in M(\psiL,\lambda)$ of $\frakz$-weight $\omega-\gamma$ with $\gamma\neq0$. The weight space $M(\psiL,\lambda)^{\omega-\gamma}=U(\bar\fraku_P)^{-\gamma}\otimes\calW_L(\lambda)$ is isomorphic as an $\frakl$-module to $E\otimes\calW_L(\lambda)$, with $E=U(\bar\fraku_P)^{-\gamma}$ finite-dimensional (Lemma~\ref{lem:McD-weights}). By Kostant's description of $E\otimes\calW_L(\lambda)$ , see Remark~\ref{rem:tensor-filtration} for $\frakl$, the space of $\psiL$-vectors of $E\otimes\calW_L(\lambda)$ is finite-dimensional, $Z(\frakl)$-stable, and the generalized eigencharacters of $Z(\frakl)$ on it are of the form $\Omega_{\lambda-\gamma'}$ with $-\gamma'$ a weight of $E$.

Therefore, we are able to choose $v$ to be a $Z(\frakl)$-eigenvector, with eigencharacter $\Omega_\mu$, $\mu=\lambda-\gamma'$; in particular $\lambda-\mu=\gamma'$ is a non-zero $\bbZ_{\geq0}$-combination of roots in $R^+\setminus R^+_L$. Since $\fraku_P\cdot v=\psiL(\fraku_P)v=0$, Frobenius reciprocity gives a non-zero homomorphism $M(\psiL,\mu)\to M(\psiL,\lambda)$ sending the generator to $v$. Since $Z(\frakg)$ acts on $M(\psiL,\lambda)$ through the central character of $M(\lambda)$ \footnote{beacuse the parabolic Harish-Chandra homomorphism $\hc_P$ is compatible with $\hc$}, we conclude that $W\cdot\mu=W\cdot\lambda$. Therefore, $\mu+\rho=w(\lambda+\rho)$, for some $w\in W$. Moreover we must have $w\notin W_L$: otherwise $\mu|_\frakz=\lambda|_\frakz$, whereas $\mu|_\frakz=\omega-\gamma\neq\omega$.

Thus $w(\lambda+\rho)-(\lambda+\rho)=-\gamma'$ lies in the root lattice. By a theorem of Jantzen (see e.g.\ \cite[\S3.4]{humphreysrepresentations}), the following set:
\begin{equation*}
    W_{[\lambda+\rho]}:=\{w\in W\mid w(\lambda+\rho)-(\lambda+\rho)\in\bbZ R\}
\end{equation*}
is the Weyl group of the root subsystem:
\begin{equation*}
    R_{[\lambda+\rho]}:=\{\gamma\in R\mid\langle\lambda+\rho,\check\gamma\rangle\in\bbZ\}
\end{equation*}
i.e.\ it is generated by the reflections $s_\gamma$ with $\langle\lambda+\rho,\check\gamma\rangle\in\bbZ$. Since $w\in W_{[\lambda+\rho]}\setminus W_L$, the subsystem $R_{[\lambda+\rho]}$ is not contained in $R_L$, which is the assertion of the lemma.
\end{proof}

\begin{corollary}\label{cor:simple}
If $\lambda\in\frakt^*_\alpha$, i.e.\ $\langle\lambda,\check\alpha\rangle\notin\bbZ$ for all $\alpha\in R^+\setminus R^+_L$, then $M(\psiL,\lambda)$ is simple. In particular $M(\psiL,\lambda)^{\fraku,\psiL}=\calW_L(\lambda)^{\fraku_L,\psiL}=\bbC\cdot1$.
\end{corollary}

\begin{proof}
Combine Lemma~\ref{lem:submodule}, Lemma~\ref{lem:McD-max} and Kostant's theorem for $\frakl$.
\end{proof}

By Remark~\ref{rem:specialization}, the asymptotic counterpart of the hypothesis of Corollary~\ref{cor:simple} is $(y,a)\in\aal$, $a\neq0$: the specialization of $\bM_\hbar(\psiL)$ at such a point is the simple module $M(\psiL,y/a)$.

\subsection{The Verma--Whittaker module off the \texorpdfstring{$\alpha$}{alpha}-walls}\label{sec:VW-generic}
Let $\mathrm{pr}\colon\bM_\hbar(\psiL)\to\bar\fraku_P\backslash\bM_\hbar(\psiL)=\calW_{L,\hbar}$ be the projection onto the top $\frakz$-weight space (Lemma~\ref{lem:McD-weights}). It is a map of $(\Uh(\bar{\mathfrak p}),\Zh(\frakl))$-bimodules, where $\bar{\mathfrak p} =\frakl\oplus\bar\fraku_P$. Therefore $\mathrm{pr}$ is $U_L$-equivariant for the adjoint actions, but it is not $\fraku_P$-equivariant. By Frobenius reciprocity it corresponds to a map of $(U,\Zh(\frakl))$-bimodules
\begin{equation}\label{eq:Theta}
\Theta\colon\bM_\hbar(\psiL)\longrightarrow\Ind_{U_L}^U(\calW_{L,\hbar})=\Hom_{U_L}(U,\calW_{L,\hbar})\cong\calO(U_P)\otimes\calW_{L,\hbar},\qquad \Theta(v)(u)=\mathrm{pr}(u\cdot v).
\end{equation}
More explicitly, after identifying $\calO(U_P)$ with $\Uh(\fraku_P)^*$, we have:
\begin{equation*}
    \Theta(y\otimes w)(u)=\mathrm{pr}_\frakl(uy)\cdot w, \quad \forall y\in\Uh(\bar\fraku_P), u\in\Uh(\fraku_P)
\end{equation*}
 where $\mathrm{pr}_\frakl$ is the projection $ \Uh(\frakg) \twoheadrightarrow \Uh(\frakg) /(\bar\fraku_P\Uh(\frakg)+\Uh(\frakg)\fraku_P)$, which can also be viewed as a parabolic Shapovalov pairing.

\begin{theorem}\label{thm:Theta}
The homomorphism $\Theta$ becomes an isomorphism after removing the $\alpha$-walls:
\[
\Theta_\alpha\colon\bM_\hbar(\psiL)\otimes_{\Zh(\frakl)}\Cal\xrightarrow{\ \sim\ }\bigl(\Ind_{U_L}^U(\calW_{L,\hbar})\bigr)\otimes_{\Zh(\frakl)}\Cal
\]
as $(U,\Cal)$-bimodules.
\end{theorem}

We first treat the specializations at $(y,a)\in\aal$, $a\neq0$.

\begin{lemma}\label{lem:Theta-inj}
For $\lambda\in\frakt^*_\alpha$, the specialization $\Theta_\lambda\colon M(\psiL,\lambda)\to\Ind_{U_L}^U(\calW_L(\lambda))$ is injective.
\end{lemma}

\begin{proof}
The kernel of $\Theta_\lambda$ is a $U$-submodule of $M(\psiL,\lambda)$ on which $\fraku^{\psiL}$ acts locally nilpotently; if it is non-zero it contains a non-zero $\psiL$-vector $v$. For a $\psiL$-vector one has $\Theta_\lambda(v)(u)=\psiL(u)\,\mathrm{pr}(v)$, so $\Theta_\lambda(v)=0$ forces $\mathrm{pr}(v)=0$, i.e.\ all $\frakz$-weight components of $v$ have weight $<\omega:=\lambda|_\frakz$. Since $M(\psiL,\lambda)^{\fraku,\psiL}$ is $\frakz$-graded (\S\ref{sec:N-properties}), this produces a $\psiL$-vector of weight $<\omega$, contradicting the simplicity of $M(\psiL,\lambda)$ (Corollary~\ref{cor:simple} and Lemma~\ref{lem:McD-max}).
\end{proof}

Since $\mathrm{pr}$ only sees the top $\frakz$-weight, $\Theta_\lambda=\bigoplus_{\gamma}\Theta_\gamma$ with
\[
\Theta_\gamma\colon U(\bar\fraku_P)_{-\gamma}\otimes\calW_L(\lambda)\longrightarrow U(\fraku_P)_\gamma^{\,*}\otimes\calW_L(\lambda),\qquad
\Theta_\gamma(y\otimes w)(x)=\mathrm{pr}_\frakl(xy)\cdot w .
\]

\begin{lemma}\label{lem:Theta-l-linear}
Let $\frakl$ act on both sides of $\Theta_\gamma$ by the tensor product action, via the adjoint action $\ad(\frakl)$ on $U(\bar\fraku_P)_{-\gamma}$ and on $U(\fraku_P)^{\,*}_\gamma$. Then $\Theta_\gamma$ is $\frakl$-linear, and both sides are $F\otimes\calW_L(\lambda)$ with $\dim F=\dim U(\fraku_P)_{\gamma}$.
\end{lemma}

\begin{proof}
In $U(\frakg)\otimes_{U(\frakp_L)}\calW_L(\lambda)$ we have $l\cdot(y\otimes w)=[l,y]\otimes w+y\otimes lw$, because $l\in\frakp_L$. The projection $\mathrm{pr}$ is $\frakl$-linear. For $x\in U(\fraku_P)_\gamma$:
\[
\Theta(lv)(x)=\mathrm{pr}(xlv)=\mathrm{pr}(l\,xv)-\mathrm{pr}([l,x]v)=l\,\Theta(v)(x)-\Theta(v)([l,x])=(l\cdot\Theta(v))(x).
\]
Here $\ad(\frakl)$ preserves $U(\fraku_P)_\gamma$, since $\frakl$ normalizes $\fraku_P$ and commutes with $\frakz$.
\end{proof}

\begin{proposition}\label{prop:Theta-gamma-iso}
For $\lambda\in\frakt^*_\alpha$ every $\Theta_\gamma$, hence $\Theta_\lambda$, is an isomorphism.
\end{proposition}

\begin{proof}
$\Theta_\gamma$ is injective by Lemma~\ref{lem:Theta-inj}, and by (K2) of Lemma~\ref{lem:Kostant-for-s} the source and the target have the spaces of the Whittaker vectors of the same finite dimension. Since $\mathrm{Wh}$ is exact, $\mathrm{Wh}(\operatorname{coker}\Theta_\gamma)=0$, so $\operatorname{coker}\Theta_\gamma=0$ by (K1) of Lemma~\ref{lem:Kostant-for-s}.
\end{proof}

\begin{proof}[Proof of Theorem~\ref{thm:Theta}]
The same argument works in families. In the universal setting, $\Theta_\gamma$ is a map of $(\Uh(\frakl),\Zh(\frakl))$-bimodules $F_\gamma\otimes\calW_{L,\hbar}\to F'_\gamma\otimes\calW_{L,\hbar}$. Applying $\mathrm{Wh}$ to this gives a map:
\[
\mathrm{Wh}(\Theta_\gamma)\colon(F_\gamma\otimes\calW_{L,\hbar})^{\fraku_L,\psiL}\to(F'_\gamma\otimes\calW_{L,\hbar})^{\fraku_L,\psiL}.
\]
Both the domain and the codomain are free $\Zh(\frakl)$-modules of rank $\dim F_\gamma$, by Kostant's freeness theorem, or by Lemma~\ref{lem:GK} together with $\bar\fraku_L\backslash(F\otimes\calW_{L,\hbar})\cong F\otimes\Sh$. The formation of $\mathrm{Wh}$ commutes with specialization, by exactness of $\mathrm{Wh}$ and freeness of $\calW_{L,\hbar}$ over $\Zh(\frakl)$.

Let $d_\gamma\in\Zh(\frakl)\subset\Sh$ be the determinant of $\mathrm{Wh}(\Theta_\gamma)$. By Proposition~\ref{prop:Theta-gamma-iso}, $d_\gamma$ does not vanish at any point $(y,a)\in\aal$ with $a\neq0$. Hence every irreducible component of $\{d_\gamma=0\}$ is either an $\alpha$-wall or the hyperplane $\{\hbar=0\}$. The second possibility is excluded because $d_\gamma|_{\hbar=0}\neq0$, which is the classical statement at $\hbar=0$ \footnote{this follows from \cite[Proposition 3.4.12]{macerato2023leviequivariantrestrictionsphericalperverse}.}. Thus $d_\gamma$ is a unit in $\Cal$, so $\mathrm{Wh}(\Theta_\gamma)_\alpha$ is an isomorphism. The relative version of (K1) of Lemma~\ref{lem:Kostant-for-s} then shows that $(\Theta_\gamma)_\alpha$ is an isomorphism, by the same cokernel argument applied after localization.
\end{proof}

\begin{corollary}\label{cor:Theta-inv}
Let $V\in\Rep(G)$. After removing the $\alpha$-walls, the map $\id_V\otimes\mathrm{pr}$ induces an isomorphism of $(\Uh(\frakl),\Cal)$-bimodules
\begin{equation}\label{eq:step2}
\bigl(V\otimes\bM_\hbar(\psiL)\bigr)^{\fraku_P}\otimes_{\Zh(\frakl)}\Cal\ \xrightarrow{\ \sim\ }\ \bigl(V\otimes\calW_{L,\hbar}\bigr)\otimes_{\Zh(\frakl)}\Cal .
\end{equation}
\end{corollary}

\begin{proof}
By Theorem~\ref{thm:Theta} and since $V$ is finite-dimensional,
\[
\bigl(V\otimes\bM_\hbar(\psiL)\bigr)^{\fraku_P}_\alpha\cong\bigl(V\otimes\Ind_{U_L}^U\calW_{L,\hbar}\bigr)^{U_P}_\alpha\cong\bigl(\Ind_{U_L}^U(V|_{U_L}\otimes\calW_{L,\hbar})\bigr)^{U_P}_\alpha\cong (V\otimes\calW_{L,\hbar})_\alpha ,
\]
by the tensor identity and Frobenius reciprocity \footnote{roughly speaking, a $U_P$-invariant function in $\Hom_{U_L}(U,N)$ is determined by its value at $1$, because $U=U_LU_P$.}. The composition is induced by $\id_V\otimes\mathrm{pr}$.
\end{proof}

\subsection{Construction of \texorpdfstring{$\widetilde\Psi_{V,\alpha}$}{Psi\_alpha}}\label{sec:PsiA}
Denote $E:=(V\otimes\bM_\hbar(\psiL))^{\fraku_P}$ and $E':=V\otimes\calW_{L,\hbar}$, both are $(\Uh(\frakl),\Zh(\frakl))$-bimodules. Since $\psiL|_{\fraku_P}=0$, $(V\otimes\bM_\hbar(\psiL))^{\fraku,\psiL}=E^{\fraku_L,\psiL}$, and $\id_V\otimes\mathrm{pr}$ induces\begin{equation}\label{eq:prV}
\mathrm{pr}_{\mathcal{W}}\colon(V\otimes\bM_\hbar(\psiL))^{\fraku,\psiL}=E^{\fraku_L,\psiL}\longrightarrow E'^{\,\fraku_L,\psiL}=(V|_L\otimes\calW_{L,\hbar})^{\fraku_L,\psiL}.
\end{equation}

\begin{proposition}\label{prop:prV}
The homomorphism $\mathrm{pr}_{\mathcal{W}}$ becomes an isomorphism after removing the $\alpha$-walls:
\[
(\mathrm{pr}_{\mathcal{W}})_\alpha\colon(V\otimes\bM_\hbar(\psiL))^{\fraku,\psiL}_\alpha\xrightarrow{\ \sim\ }(V|_L\otimes\calW_{L,\hbar})^{\fraku_L,\psiL}_\alpha .
\]
\end{proposition}

\begin{proof}
$\fraku_L^{\psiL}$ acts locally nilpotently on $E$ and $E'$ (through the adjoint action, Remark~\ref{rem:U-inv}), so Lemma~\ref{lem:GK} for $(\fraks,\psi_{\mathrm{der}})$, tensored with $\Uh(\frakz)$, gives the commutative square below:
\begin{equation}\label{eq:GK-square}
\begin{tikzcd}[column sep=large]
\Sh\otimes_{\Zh(\frakl)}(V\otimes\bM_\hbar(\psiL))^{\fraku,\psiL}\arrow[r,"\sim"]\arrow[d,"\id\otimes\mathrm{pr}_{\mathcal{W}}"'] & \bar\fraku_L\backslash E\arrow[d,"\bar\fraku_L\backslash\mathrm{pr}_{\mathcal{W}}"]\\
\Sh\otimes_{\Zh(\frakl)}(V|_L\otimes\calW_{L,\hbar})^{\fraku_L,\psiL}\arrow[r,"\sim"] & \bar\fraku_L\backslash E' .
\end{tikzcd}
\end{equation}
The right vertical map is an isomorphism over $\Cal$ by Corollary~\ref{cor:Theta-inv}, and $\Sh$ is faithfully flat over $\Zh(\frakl)$.
\end{proof}

By \eqref{eq:BF} for $\cL$,
\begin{equation}\label{eq:levi-satake}
(V|_L\otimes\calW_{L,\hbar})^{\fraku_L,\psiL}\cong\H_{\LO}(\GrL,\bbS_L^{-1}(V|_L))=\H_{\LO}(\GrL,\Res^{\cG}_{\cL}\IC_V).
\end{equation}
Together with Proposition~\ref{prop:geom-loc}(1) we define the isomorphism of $\Zh(\frakl)_\alpha$-modules
\begin{equation}\label{eq:PsiA-def}
\Psi_{V,\alpha}:=\H_{\LO}(\Xi)_\alpha^{-1}\circ\eqref{eq:levi-satake}_\alpha\circ(\mathrm{pr}_{\mathcal{W}})_\alpha\colon(V\otimes\bM_\hbar(\psiL))^{\fraku,\psiL}_\alpha\xrightarrow{\ \sim\ }\H_{\LO}(\GrL,i_L^!\IC_V)_\alpha ,
\end{equation}
and $\widetilde\Psi_{V,\alpha}:=\Sh\otimes_{\Zh(\frakl)}\Psi_{V,\alpha}$.

\begin{proposition}\label{prop:PsiA}
$\widetilde\Psi_{V,\alpha}$ is a canonical $W_L$-equivariant isomorphism \eqref{eq:PsiA}.
\end{proposition}

In other words, off the $\alpha$-walls, an intertwiner $\bM_\hbar(\psiL)\to V\otimes\bM_\hbar(\psiL)$ is the same as an intertwiner $\calW_{L,\hbar}\to V|_L\otimes\calW_{L,\hbar}$, via projection to the top $\frakz$-weight. This is the quantized version of \cite[Corollary 3.3.10]{macerato2023leviequivariantrestrictionsphericalperverse}.

\section{Case 2: removing the \texorpdfstring{$\beta$}{beta}-walls}\label{sec:case2}

We work over $\Cbe$ throughout this section; in the asymptotic setting, $\hbar$ is moreover inverted (Remark~\ref{rem:KK-asymptotic}). Geometrically, Proposition~\ref{prop:geom-loc}(2) reduces to \eqref{eq:GR}. Algebraically, we use Backelin's functor to find $\calW_{L,\hbar}$ inside the completed $\bM_{L,\hbar}$. We first treat $L=G$ over $\bbC[\frakt^*_\rs]$, then apply it to $\frakl$ and induce (\S\ref{sec:case2-iso}).

\subsection{The Shapovalov form and Whittaker vectors in Verma modules}\label{sec:shapovalov}
Let $\frakg$ be semisimple. Grade $\bM=\bigoplus_{\eta\geq0}\bM_{-\eta}$, $\bM_{-\eta}=U(\bar\fraku)_{-\eta}\otimes \SymAlg(\frakt)$, by relative weight (\S\ref{sec:left-action}), and let $\widehat{\bM}:=\prod_\eta\bM_{-\eta}$, and similar notation after base change. A \emph{Whittaker vector} $\wv=\sum_\eta\wv_{-\eta}$ in a completion satisfies $x\wv=\psi(x)\wv$ for $x\in\fraku$; it is \emph{normalized} if $\wv_0=1$. It is well-known that over the fraction field $\bbC(\frakt^*)$, there is a unique normalized one, cf.\ \cite[\S8]{braverman2016instantonmodulispacesmathscr}. Indeed, it is not hard to use the Shapovalov form to control its denominators, refining the result in the loc.cit.

\begin{theorem}\label{thm:whittaker-vector}
There exists a unique normalized Whittaker vector $\wv_\psi$ in the completed universal Verma module $\widehat{\bM\otimes_{\SymAlg(\frakt)}\bbC[\frakt^*_\rs]}$ over $\bbC[\frakt^*_\rs]$. Every Whittaker vector in $\widehat{\bM\otimes_{\SymAlg(\frakt)}\bbC[\frakt^*_\rs]}$ is of the form $c\,\wv_\psi$ with $c\in\bbC[\frakt^*_\rs]$.
\end{theorem}

Let $\tau$ be the Chevalley anti-involution ($\tau(e_\gamma)=f_\gamma$, $\tau|_\frakt=\id$) and $\mathrm{pr}_\frakt$ the projection of \eqref{eq:HC-decomp}. The \emph{Shapovalov form}
\begin{equation}\label{eq:shapovalov}
F(x,y):=\mathrm{pr}_\frakt(\tau(x)\,y)\in\SymAlg(\frakt)
\end{equation}
is contravariant, $F(fx,y)=F(x,\tau(f)y)$, and descends to $\bM$ with distinct weight spaces orthogonal. Let $F_\eta$ be its restriction to $\bM_{-\eta}$, free of rank $\mathscr P(\eta)$, where $\mathscr P(\eta)$ is Kostant's partition function. We omit its definition here.

\begin{lemma}[{\cite[Theorem 1]{KAC197997}}]\label{lem:KK}
Fix a $W$-invariant non-degenerate symmetric bilinear form $(\,,)$ on $\frakt^*$ and use it to identify $\frakt^*\cong\frakt$, so that elements of $\frakt^*$ become linear functions on $\frakt^*$. Then, up to a non-zero constant depending on the choice of basis of $U(\bar\fraku)_{-\eta}$,
\begin{equation}\label{eq:KK}
\det F_\eta=\prod_{\gamma\in R^+}\ \prod_{n=1}^{\infty}\Bigl(\gamma+(\rho,\gamma)-\tfrac{n}{2}(\gamma,\gamma)\Bigr)^{\mathscr P(\eta-n\gamma)} .
\end{equation}
\end{lemma}

Since $(\lambda,\gamma)=\tfrac12(\gamma,\gamma)\langle\lambda,\check\gamma\rangle$ for $\lambda\in\frakt^*$, the factor indexed by $(\gamma,n)$ in \eqref{eq:KK} is, up to a non-zero constant, the linear function $\lambda\mapsto\langle\lambda+\rho,\check\gamma\rangle-n$, i.e.\ it is the linear form $\check\gamma+(\langle\rho,\check\gamma\rangle-n)$ on $\frakt^*$. Thus the zero locus of $\det F_\eta$ is contained in the union of the hyperplanes $\{\langle\lambda,\check\gamma\rangle\in\bbZ\}$, $\gamma\in R^+$, and $\det F_\eta$ becomes a unit in $\bbC[\frakt^*_\rs]$ (Definition~\ref{def:loci}).

\begin{proof}[Proof of Theorem~\ref{thm:whittaker-vector}]
Let $\{v_J\}$ be a PBW basis of $U(\bar\fraku)$, indexed by multi-indices $J$, so that $\{v_J\otimes1\}$ is an $\SymAlg(\frakt)$-basis of $\bM$ and each $v_J\otimes1$ has weight $-\eta_J$. Let $\{e_I\}$ be the dual PBW monomials in $U(\fraku)$, $e_I=\tau(v_I)$. A vector $\wv=\sum_Jc_J\,v_J\otimes1$ of $\widehat{\bM\otimes\bbC[\frakt^*_\rs]}$, $c_J\in\bbC[\frakt^*_\rs]$, is a Whittaker vector if and only if $e_I\wv=\psi(e_I)\wv$ for all $I$, where $\psi(e_I)=\psi(e_{i_1})\cdots\psi(e_{i_r})$ for $e_I=e_{i_1}\cdots e_{i_r}$. Fix a weight $\eta$ and an index $I$ with $\eta_I=\eta$. Pairing both sides of $e_I\wv=\psi(e_I)\wv$ with $1\in\bM_0$ (only the weight-$0$ components contribute) and using contravariance, $F(1,e_I\wv)=F(\tau(e_I)\cdot1,\wv)=F(v_I,\wv)=\sum_{J:\,\eta_J=\eta}F(v_I,v_J)\,c_J$, while $\psi(e_I)F(1,\wv)=\psi(e_I)c_0$. We obtain for each $\eta$ the finite linear system
\[
\sum_{J:\ \eta_J=\eta}F(v_I,v_J)\,c_J=\psi(e_I)\,c_0\qquad(\eta_I=\eta),
\]
whose matrix is $F_\eta$. Thus the weight-$(-\eta)$ component of $\wv$ is determined by $c_0=\wv_0$:
\begin{equation}\label{eq:cJ}
c_J=\sum_{I:\ \eta_I=\eta}(F_\eta^{-1})_{JI}\,\psi(e_I)\,c_0 .
\end{equation}
Since $F_\eta$ has entries in $\SymAlg(\frakt)$ and its determinant is a unit in $\bbC[\frakt^*_\rs]$ by Lemma~\ref{lem:KK}, the matrix $F_\eta^{-1}$ has entries in $\bbC[\frakt^*_\rs]$. Hence \eqref{eq:cJ} with $c_0=1$ defines a normalized Whittaker vector $\wv_\psi$ over $\bbC[\frakt^*_\rs]$, and every Whittaker vector is $c_0\,\wv_\psi$.
\end{proof}

\begin{remark}\label{rem:KK-asymptotic}
In the asymptotic setting $\det F_{\hbar,\eta}$ is a power of $\hbar$ times a product of walls, see Lemma~\ref{lem:integrality}. So Theorem~\ref{thm:whittaker-vector} holds over $\Crs[\hbar^{-1}]$: there is a unique normalized Whittaker vector $\wv_{\psi,\hbar}$ in $\widehat{\bM_\hbar\otimes_{\Sh}\Crs[\hbar^{-1}]}$, and every Whittaker vector there is a $\Crs[\hbar^{-1}]$-multiple of it. Its coefficients do have poles along $\hbar=0$ (Lemma~\ref{lem:sl2-whittaker}(a)). Accordingly, the asymptotic statements of \S\S\ref{sec:apply-backelin}--\ref{sec:PsiB} are over $\Cbe[\hbar^{-1}]$ (resp.\ $\Crs[\hbar^{-1}]$). The fiber $\hbar=0$ is treated in Appendix~\ref{app:hbar0}: the discrepancy of \S\ref{sec:discrepancy} is regular there, see Proposition~\ref{prop:Phi-integral}, and the final isomorphism extends across it (Corollary~\ref{cor:extension}).
\end{remark}

\subsection{Applying Backelin's functor}\label{sec:apply-backelin}
Let $\frakg$ be semisimple. By Theorem~\ref{thm:whittaker-vector} the Whittaker vector $\wv_\psi$ defines a homomorphism of $(U(\frakg),\bbC[\frakt^*_\rs])$-bimodules
\begin{equation}\label{eq:phi}
\phi\colon\calW\otimes_{Z(\frakg)}\bbC[\frakt^*_\rs]\longrightarrow\widehat\Gamma_\psi\bigl(\bM\otimes_{\SymAlg(\frakt)}\bbC[\frakt^*_\rs]\bigr),\qquad 1\otimes1\longmapsto \wv_\psi ,
\end{equation}
where $Z(\frakg)\to\bbC[\frakt^*_\rs]$ is the Harish-Chandra homomorphism followed by localization. The asymptotic version is
\begin{equation}\label{eq:phi-h}
\phi_\hbar\colon\calW_\hbar\otimes_{\Zh(\frakg)}\Crs[\hbar^{-1}]\longrightarrow\widehat\Gamma_\psi\bigl(\bM_\hbar\otimes_{\Sh}\Crs[\hbar^{-1}]\bigr),\qquad 1\otimes1\longmapsto \wv_{\psi,\hbar}.
\end{equation}

\begin{proposition}\label{prop:phi-iso}
The maps $\phi$ and $\phi_\hbar$ are isomorphisms.
\end{proposition}

\begin{proof}
We treat $\phi$; the asymptotic case is identical, over $\Crs[\hbar^{-1}]$ (Remark~\ref{rem:KK-asymptotic}). Put $R_\rs:=\bbC[\frakt^*_\rs]$, $E_1:=\calW\otimes_{Z(\frakg)}R_\rs$, $E_2:=\Gamma_\psi(\widehat{\bM}_{\rs})$. Lemma~\ref{lem:Kostant-for-s} for $(\frakg,\psi)$ does not need finite generation.

\emph{Step 1.} $\fraku^\psi$ acts locally nilpotently on $E_2$ by definition, and on $E_1$ because it acts by $\ad(\fraku)$ (Remark~\ref{rem:U-inv}).

\emph{Step 2: Whittaker vectors.} $R_\rs$ is flat over $Z(\frakg)$, so $E_1^{\fraku,\psi}=\calW^{\fraku,\psi}\otimes_{Z(\frakg)}R_\rs=R_\rs(1\otimes1)$ by Kostant. Whittaker vectors of $\widehat{\bM}_\rs$ are $\psi$-finite, so by Theorem~\ref{thm:whittaker-vector}, $E_2^{\fraku,\psi}=\wv_\psi R_\rs$, free since $(\wv_\psi r)_0=r$. So $\phi^{\fraku,\psi}$ is an isomorphism:
\[
\phi^{\fraku,\psi}\colon E_1^{\fraku,\psi}=R_\rs\cdot(1\otimes 1)\xrightarrow{\ \sim\ }\wv_\psi R_\rs=E_2^{\fraku,\psi}.
\]

\emph{Step 3: injectivity.} $K:=\ker\phi$ is a submodule of $E_1$, so $\fraku^\psi$ acts locally nilpotently on it, and $K^{\fraku,\psi}=\ker(\phi^{\fraku,\psi})=0$. By (K1), $K=0$.

\emph{Step 4: surjectivity.} Let $Q=\operatorname{coker}\phi$. By right exactness of $\bar\fraku\backslash(-)$ and Lemma~\ref{lem:GK}, $\bar\fraku\backslash E_1\to\bar\fraku\backslash E_2$ is $\SymAlg(\frakt)\otimes_{Z(\frakg)}\phi^{\fraku,\psi}$, an isomorphism; so $\bar\fraku\backslash Q=0$, hence $\SymAlg(\frakt)\otimes_{Z(\frakg)}Q^{\fraku,\psi}=0$, hence $Q^{\fraku,\psi}=0$ by faithful flatness, and $Q=0$ by (K1).
\end{proof}

\begin{remark}
Since $\calW$ is not finitely generated over $Z(\frakg)$, Proposition~\ref{prop:phi-iso} cannot be deduced from Lemma~\ref{lem:backelin-verma} by specializing at each point $\lambda\in\frakt^*_\rs$ and applying Nakayama's lemma. Conversely, specializing the proof of Proposition~\ref{prop:phi-iso} at $\lambda\in\frakt^*_\rs$ recovers Lemma~\ref{lem:backelin-verma} at $\lambda$.
\end{remark}

\subsection{Parabolic induction}\label{sec:case2-iso}
Applying \S\S\ref{sec:shapovalov}--\ref{sec:apply-backelin} to $(\fraks,\psi_{\mathrm{der}})$ and tensoring with $U(\frakz)$ gives, since the walls of $\fraks$ are the $\beta$-walls,
\begin{align}
\phi_L&\colon\calW_L\otimes_{Z(\frakl)}\bbC[\frakt^*_\beta]\xrightarrow{\ \sim\ }\widehat\Gamma_{\psiL}\bigl(\bM_L\otimes_{\SymAlg(\frakt)}\bbC[\frakt^*_\beta]\bigr),\label{eq:phiL}\\
\phi_{L,\hbar}&\colon\calW_{L,\hbar}\otimes_{\Zh(\frakl)}\Cbe[\hbar^{-1}]\xrightarrow{\ \sim\ }\widehat\Gamma_{\psiL}\bigl(\bM_{L,\hbar}\otimes_{\Sh}\Cbe[\hbar^{-1}]\bigr),\label{eq:phiLh}
\end{align}
Inducing from $\frakp_L$ to $\frakg$, and using that $\widehat\Gamma_{\psiL}$ commutes with parabolic induction (\S\ref{sec:backelin}), $\bM=U(\frakg)\otimes_{U(\frakp_L)}\bM_L$ and \eqref{eq:MpsiL-def}, we obtain
\begin{align}
\phi&\colon\bM(\psiL)\otimes_{Z(\frakl)}\bbC[\frakt^*_\beta]\xrightarrow{\ \sim\ }\widehat\Gamma_{\psiL}\bigl(\bM\otimes_{\SymAlg(\frakt)}\bbC[\frakt^*_\beta]\bigr),\label{eq:phi-G}\\
\phi_\hbar&\colon\bM_\hbar(\psiL)\otimes_{\Zh(\frakl)}\Cbe[\hbar^{-1}]\xrightarrow{\ \sim\ }\widehat\Gamma_{\psiL}\bigl(\bM_\hbar\otimes_{\Sh}\Cbe[\hbar^{-1}]\bigr).\label{eq:phi-Gh}
\end{align}
These are determined by the normalized Whittaker vector, and they commute with $V\otimes(-)$.

\subsection{The map \texorpdfstring{$\Phi$}{Phi}}\label{sec:Phi}
For $\xi\in(V\otimes\bM_\hbar)^{\fraku}$, let $f_\xi$ be the corresponding intertwiner, defined by $f_\xi(\one)=\xi$ (\S\ref{sec:intertwiners}). It extends to completions. Let $\wv_L=\sum_J y_J\one\,c_J$, $y_J\in\Uh(\bar\fraku_L)$, $c_J\in\Cbe[\hbar^{-1}]$, $c_\emptyset=1$, be the normalized $\psiL$-Whittaker vector of $\widehat{\bM_{L,\hbar}\otimes_{\Sh}\Cbe[\hbar^{-1}]}\subset\widehat{\bM_\hbar\otimes_{\Sh}\Cbe[\hbar^{-1}]}$ (Remark~\ref{rem:KK-asymptotic} for $\fraks$). Define
\begin{equation}\label{eq:Phi-def}
\Phi\colon(V\otimes\bM_\hbar)^{\fraku}\otimes_{\Sh}\Cbe[\hbar^{-1}]\longrightarrow(V\otimes\bM_\hbar(\psiL))^{\fraku,\psiL}\otimes_{\Zh(\frakl)}\Cbe[\hbar^{-1}],\qquad
\xi\longmapsto(V\otimes\phi_\hbar)^{-1}\bigl(f_\xi(\wv_L)\bigr),
\end{equation}
where $\phi_\hbar$ is \eqref{eq:phi-Gh}. Explicitly, $f_\xi(\wv_L)=\sum_Jy_J\,(c_J\cdot\xi)$, with the left action. Alternatively, the left action can be viewed as the right one via a shift of Verma parameter: $c_J\cdot\xi=\xi\cdot c_J(y+\hbar\nu)$, on the $\nu$-piece (Corollary~\ref{cor:comp-verma}(a)).

Point-wisely speaking, let the sub-category $\calO_\beta\subset\calO$, consisting of modules all of whose weights $\mu$ satisfy: $\langle\mu,\check\beta\rangle\notin\bbZ$, for $\beta\in R_L$. For instance, $E\otimes M(\lambda)$, for $\lambda\in\frakt^*_\beta$. With $\wv_L=\sum_J c_J\,y_J$ the non-asymptotic Whittaker vector, we define for $N\in\calO_\beta$:
\begin{equation}\label{eq:Phi}
\Phi_N\colon N^{\fraku}\longrightarrow\bigl(\widehat\Gamma_{\psiL}N\bigr)^{\fraku,\psiL},\qquad
v\longmapsto\sum_Jc_J(\mu)\,y_J\,v\qquad(v\in N^{\fraku}_\mu).
\end{equation}
This summation makes sense in $\widehat N$, because $c_J$ is regular at $\mu\in\frakt^*_\beta$. It is the image of $\wv_L$ under $M(\mu)\to N$, $1\mapsto v$, hence a $\psiL$-Whittaker vector; for $N=V\otimes M(\lambda)$, $\Phi_N$ is the specialization of \eqref{eq:Phi-def}.

\begin{proposition}\label{prop:Phi-iso}
For every $N\in\calO_\beta$ the map $\Phi_N$ is an isomorphism. In particular, for $\lambda\in\frakt^*_\beta$ and $E$ finite-dimensional,
\[
(E\otimes M(\lambda))^{\fraku}=\bigoplus_\mu\Hom_{\calO}\bigl(M(\mu),E\otimes M(\lambda)\bigr)\xrightarrow{\ \sim\ }\Hom_{U(\frakg)}\bigl(\bM(\psiL),\,E\otimes M(\psiL,\lambda)\bigr)=(E\otimes M(\psiL,\lambda))^{\fraku,\psiL}.
\]
\end{proposition}

\begin{proof}
\emph{Step 1: reduction to the Levi.} Denote by $N_L:=N^{\fraku_P}$, which is an $\frakl$-module. Since $\fraku_P$ acts by weight vectors, $\widehat N^{\fraku_P}=\widehat{N_L}$, and on it we have: $(u-\psiL(u))^k=(u_L-\psiL(u_L))^k$, for $u=u_L+u_P$. Therefore:
\[
(\widehat\Gamma_{\psiL}N)^{\fraku_P}=\widehat\Gamma^L_{\psiL}(N_L),\qquad
(\widehat\Gamma_{\psiL}N)^{\fraku,\psiL}=\bigl(\widehat\Gamma^L_{\psiL}N_L\bigr)^{\fraku_L,\psiL},\qquad
N^{\fraku}=N_L^{\fraku_L},
\]
where $\widehat\Gamma^L_{\psiL}$ is Backelin's functor for $(\frakl,\psiL|_{\fraku_L})$. Since $y_J\in U(\bar\fraku_L)$, the formula \eqref{eq:Phi} shows that $\Phi_N=\Phi^L_{N_L}$, which is the analogous map defined for $\frakl$.

\emph{Step 2: structure of $N_L$.} We decompose $N_L=\bigoplus_\omega N_L^\omega$, over $\frakz$-weights. As an $\frakl$-module, the $\frakz$-weight space $N^\omega$ has a finite filtration with subquotients: (finite dimensional modules)$\otimes M_L(\lambda_i)$. Thus, $N_L^\omega\in\calO_L$. Its weights $\mu$ satisfy $\langle\mu,\check\beta\rangle\notin\bbZ$. This implies each $M_L(\mu)$ is simple and projective in $\calO_L$ \cite[Proposition 3.8, Theorem 4.8]{humphreysrepresentations}. Therefore we have: $N_L^\omega\cong\bigoplus_i M_L(\mu_i)$, which is a finite sum.

\emph{Step 3: one Verma module.} Note that $M_L(\mu)^{\fraku_L}=\bbC v_\mu$, and $\widehat\Gamma^L_{\psiL}(M_L(\mu))\cong\calW_L(\mu)$ has Whittaker space $\bbC\cdot1$. We see that $\Phi^L$ sends $v_\mu \longmapsto \sum_Jc_J(\mu)y_Jv_\mu$, which is non-zero.

\emph{Step 4.} By Steps 2 and 3, we see that $\Phi^L_{N_L^\omega}$ is an isomorphism for each $\omega$. Since $\fraku_L$ commutes with $\frakz$, taking Whittaker vectors commutes with the $\frakz$-weight decomposition. As a priori, the domain is $\bigoplus_\omega$, and the codomain is $\prod_\omega$. But that doesn't matter, because:
\[
\dim(\widehat\Gamma^L_{\psiL}N_L^\omega)^{\fraku_L,\psiL}=\dim(N_L^\omega)^{\fraku_L}, \quad \mathrm{and} \quad \sum_\omega\dim(N_L^\omega)^{\fraku_L}=\dim N^{\fraku}<\infty
\]
Therefore, only finitely many $\omega$ contribute to either side.
\end{proof}

In order to pass from the pointwise isomorphism to the universal isomorphism, we need to use the following general algebraic fact:
\begin{lemma}\label{lem:height-one}
Let $R$ be a noetherian normal domain, with fraction field $K =\mathrm{Frac}(R)$. let $F$ be a free $R$-module of finite rank, let $T$ be a torsion-free $R$-module, and let $\Phi\colon F\to T$ be $R$-linear. Assume:
\begin{enumerate}[label=(\roman*)]
\item $\Phi\otimes_RK\colon F\otimes_RK\to T\otimes_RK$ is an isomorphism;
\item for every prime $\frakq\subset R$ of height one, $\Phi\otimes_Rk(\frakq)\colon F\otimes_Rk(\frakq)\to T\otimes_Rk(\frakq)$ is injective.
\end{enumerate}
Then $\Phi$ is an isomorphism.
\end{lemma}

\begin{proof}
$\Phi$ is injective because $F\subset F\otimes_RK$. For surjectivity, let $\zeta\in T$ and put $\xi:=(\Phi\otimes K)^{-1}(\zeta)\in F\otimes_RK$. Since $R$ is noetherian and normal, $F$ is free, it follows that $F=\bigcap_{\operatorname{ht}\frakq=1}F_\frakq$ inside $F\otimes_RK$. Therefore, it suffices to show $\xi\in F_\frakq$ for each height-one prime $\frakq$. Suppose the contrary. Let $g$ be a uniformizer of the discrete valuation ring $R_\frakq$, and write $\xi=g^{-m}\xi'$ with $m\ge1$ and $\xi'\in F_\frakq\smallsetminus gF_\frakq$. Then $\Phi(\xi')=g^m\zeta\in gT_\frakq$, so $\Phi\otimes k(\frakq)$ kills the image of $\xi'$ in $F_\frakq/gF_\frakq=F\otimes k(\frakq)$. However, that image is non-zero, contradicting (ii).
\end{proof}

\begin{theorem}\label{thm:full-faithful}
The map \eqref{eq:Phi-def} is an isomorphism:
\begin{equation}\label{eq:case2-alg}
\Phi\colon(V\otimes\bM_\hbar)^{\fraku}\otimes_{\Sh}\Cbe[\hbar^{-1}]\xrightarrow{\ \sim\ }(V\otimes\bM_\hbar(\psiL))^{\fraku,\psiL}\otimes_{\Zh(\frakl)}\Cbe[\hbar^{-1}] .
\end{equation}
In particular the target is free of rank $\dim V$ over $\Cbe[\hbar^{-1}]$.
\end{theorem}

\begin{proof}
We apply Lemma~\ref{lem:height-one} with $R:=\Cbe[\hbar^{-1}]$, $K:=\operatorname{Frac}R$, $F$ the domain and $T$ the codomain of $\Phi$. Put $E:=V\otimes\bM_\hbar(\psiL)\otimes_{\Zh(\frakl)}R$. Since $R$ is flat over $\Zh(\frakl)$, we have $T=E^{\fraku,\psiL}$.

\emph{Hypotheses on $F$ and $T$.} $F$ is free of rank $\dim V$ by \eqref{eq:GR} and Lemma~\ref{lem:formal}. $E$ is free over $R$, because $\calW_{L,\hbar}$ is free over $\Zh(\frakl)$, by Kostant, and $\bM_\hbar(\psiL)=\Uh(\bar\fraku_P)\otimes\calW_{L,\hbar}$. Hence $T\subset E$ is torsion-free.

\emph{Specialization.} Since the $c_J$ lie in $R$, $\Phi\otimes k(\frakq)$ is given by the same formula. Every prime of $R$ avoids $\hbar$, so $y/\hbar$ is defined over $k(\frakq)$. For $\operatorname{ht}\frakq\le1$, $F\otimes k(\frakq)\to V\otimes\bM_\hbar\otimes k(\frakq)$ is injective \footnote{$F_\frakq$ is a kernel of a homomorphism between free modules, over a  DVR if $\operatorname{ht}\frakq =1$, over a field if $\frakq =0$.  Therefore it must be a direct summand}. Therefore, it suffices to show that the specialized $\Phi_N$, $N=V\otimes M(y/\hbar)$, over $k(\frakq)$, is injective, and bijective for $\frakq=0$.

Then let $y/\hbar\in\frakt^*_\beta$ over $k(\frakq)$, and Proposition~\ref{prop:Phi-iso}, whose proof works over any field of characteristic zero, gives injectivity. It also gives bijectivity for $\frakq=0$. This is (ii); (i) follows since localization is exact.
\end{proof}

\subsection{The map \texorpdfstring{$\widetilde\Psi_{V,\beta}$}{Psi\_beta}}\label{sec:PsiB}
Using the isomorphisms we have proved so far, let's put:
\begin{equation}\label{eq:PsiB-def}
\widetilde\Psi_{V,\beta}:=\eqref{eq:case2-geom}\circ\eqref{eq:GR}_\beta\circ\Phi^{-1}\colon(V\otimes\bM_\hbar(\psiL))^{\fraku,\psiL}\otimes_{\Zh(\frakl)}\Cbe[\hbar^{-1}]\xrightarrow{\ \sim\ }\H_A(\GrL,i_L^!\IC_V)_\beta[\hbar^{-1}] .
\end{equation}
It is an isomorphism. It must still be corrected (\S\ref{sec:gluing}), and extended across $\hbar=0$ (Appendix~\ref{app:hbar0}). Moreover, at this time, $\widetilde\Psi_{V,\beta}$'s $W_L$-equivariance property is not visible through the intermediate terms. Recall the $W_L$-symmetry of \eqref{eq:GR} involves the dynamical Weyl group \cite[\S\S2.5--2.6]{ginzburg2014differentialoperatorsguaffine}, whose construction is rather complicated. We are not able to construct the dynamical $W_L$ action directly on the algebraic side; instead, for the corrected map $\widetilde \Psi_{V,\beta}''$, its $W_L$-equivariance follows from the compatibility with $\widetilde \Psi_{V,\alpha}$ on the regular semisimple locus and the $W_L$-equivariance of Case~1 (Corollary~\ref{cor:WL}).

\section{Hartogs property for $D_\hbar(G/U,\psiL)$}\label{sec:hartogs}

Following \cite[\S3.5]{macerato2023leviequivariantrestrictionsphericalperverse}, we prove \eqref{eq:alg-glue} using the geometry of $T^*(G/U,\psiL)\to\frakc_L$. What's new here is that we also need to use the idea of localization to sheafify $\Dh(G/U,\psiL)$ over the base affine space $G/U$.

\begin{theorem}\label{thm:alg-glue}
The natural inclusion
\[
\Dh(G/U,\psiL)\longrightarrow\Dh(G/U,\psiL)_\alpha\times_{\Dh(G/U,\psiL)_\rs}\Dh(G/U,\psiL)_\beta
\]
is an isomorphism.
\end{theorem}
\subsection{Sheaf of TDO}
We will follow \cite[\S\S6.3--6.4]{ginzburg2022differentialoperatorsgugelfandgraev} to sheafify $\Dh(G/U,\psiL)$. Put $X:=G/U$ and let $p\colon G\to X$ be the projection. Under $p$, $G$ becomes a right $U$-torsor over $X$. Let
$D_\hbar(G)=\mathcal O(G)\otimes U_\hbar(\mathfrak g)$ be the Rees algebra for the order filtration on $D(G)$.
Here $x\in\mathfrak g$ is sent to the left-invariant vector field $\xi_x$ generating right translations, and commutator is
$[\xi_x,f]=\hbar\,\xi_x(f)$. Denote
\[
\mathcal M:=D_\hbar(G)/D_\hbar(G)\,\mathfrak u^{\psi_L}\;=\;\mathcal O(G)\otimes \bM_\hbar(\psiL)
\]
The right translation $r_u$ induced by $u\in U$ is an algebra automorphism of $D_\hbar(G)$, with $r_u(\xi_x)=\xi_{\mathrm{Ad}(u)x}$.
It preserves the left ideal $D_\hbar(G)\mathfrak u^{\psi_L}$, because $\psi_L\circ\mathrm{Ad}(u)=\psi_L$ on $\mathfrak u$. Therefore, $U$ does act on $\mathcal M$. We write $\one\in\mathcal M$ for the class of $1$, so that the class of $a\in D_\hbar(G)$ is $a\one$; note $r_u(\one)=\one$. We regard $\mathcal M$ as a quasi-coherent sheaf on the affine variety $G$. Then we denote:
\[
\mathcal D:=\mathcal D^{\psi_L}_{G/U,\hbar}:=\bigl(p_*\mathcal M\bigr)^U .
\]
\footnote{This is the asymptotic version of the TDO: $\mathcal D^{\psi}_X\cong(p_*(\mathcal D_P/\mathcal D_P\mathfrak n_P^\psi))^N$ of
\cite[\S6.3]{ginzburg2022differentialoperatorsgugelfandgraev}, in the special case of the constant group scheme $N=U\times X$ and the constant section $\psi=\psi_L$.
The $\mathfrak k$-invariance hypothesis of loc.\ cit.\ holds because $\psi_L$ is $\mathrm{Ad}(U)$-invariant.} We list the following properties of $\mathcal D$:

\begin{itemize}
\item \emph{Algebra structure.} If $b\one\in\mathcal M$ is $U$-invariant, it is $\mathfrak u$-invariant. More concretely speaking, 
$[\xi_x,b]\in D_\hbar(G)\mathfrak u^{\psi_L}$, for $x\in\mathfrak u$.
Therefore, $D_\hbar(G)\mathfrak u^{\psi_L}\cdot b\subset D_\hbar(G)\mathfrak u^{\psi_L}$, and $(a\one)(b\one):=ab\one$ is a well-defined
product on $\mathcal D$. Since $\mathcal O(G)^U=\mathcal O(X)$, we get $\mathcal O_X\subset\mathcal D$.
\item \emph{Right $Z_\hbar(\mathfrak l)$-action.} Put $I:=U_\hbar(\mathfrak p_L)\mathfrak u_P$. Since $\mathfrak u_P$ is an ideal of $\mathfrak p_L$, we have $U_\hbar(\mathfrak p_L)\mathfrak u_P=\mathfrak u_PU_\hbar(\mathfrak p_L)$, so $I$ is a two-sided ideal of $U_\hbar(\mathfrak p_L)$. Moreover $I\subset U_\hbar(\mathfrak g)\mathfrak u^{\psi_L}$, because $x-\psi_L(x)=x$ for $x\in\mathfrak u_P$.

For $x\in\mathfrak u$, $z\in Z_\hbar(\frakl)$, we show that $[x,z]\in I$. Indeed, if $x\in\mathfrak u_L$, then $[x,z]=0$, because $z$ is central in $U_\hbar(\mathfrak l)$.
If $x\in\mathfrak u_P$, we write $z$ as a sum of products $l_1\cdots l_k$ with $l_i\in\mathfrak l$. Then $[x, l_1\cdots l_k]=\sum_il_1\cdots[x,l_i]\cdots l_k$ with $[x,l_i]\in\mathfrak u_P$, and each term lies in $I$ because $I$ is a two-sided ideal of $U_\hbar(\mathfrak p_L)$.
Therefore, $\mathfrak u^{\psi_L}z\subset z\,\mathfrak u^{\psi_L}+I\subset U_\hbar(\mathfrak g)\mathfrak u^{\psi_L}$, and $a\one\mapsto az\one$ is well defined on $\mathcal M$. From this, we deduce that the action of $Z_\hbar(\frakl)$ commutes with $U$. Thus, $r_u(az\one)=r_u(a)\,r_u(z)\one=r_u(a)z\one$. Hence $\mathcal D$ is a sheaf of right $Z_\hbar(\mathfrak l)$-modules. This
is right multiplication by the image of $Z_\hbar(\mathfrak l)\to\Gamma(X,\mathcal D)$, $z\mapsto z\one$.
\item \emph{Global sections.} $\Gamma(X,\mathcal D)=\Gamma(G,\mathcal M)^U=(\mathcal O(G)\otimes\bM_\hbar(\psiL))^U
=D_\hbar(G/U,\psi_L)$.
\end{itemize}

\subsection{The twisted cotangent bundle}
Let $T^*(G/U,\psi_L):=G\times_U(\psi_L+\mathfrak u^\perp)$, where $U$ acts on $\psi_L+\mathfrak u^\perp\subset\mathfrak g^*$ by
the coadjoint action. The space $\psi_L+\mathfrak u^\perp$ is $\mathrm{Ad}(U)$-stable, due to the fact that $\psi_L$ is $U$-invariant.
Let $\pi\colon T^*(G/U,\psi_L)\to X$ be the projection; it is an affine bundle. Then denote:
\[
\pi_L\colon T^*(G/U,\psi_L)\to\mathfrak c_L=\mathfrak l^*/\!/L,\qquad [g,\xi]\mapsto\chi_L(\xi|_{\mathfrak l}),
\]
where $\chi_L\colon\mathfrak l^*\to\mathfrak c_L$ is the adjoint quotient. This is well defined. For $u\in U_L$ use the
$L$-invariance of $\chi_L$. For $u\in U_P$ and $y\in\mathfrak l$, we have $\mathrm{Ad}(u)y-y\in\mathfrak u_P$ , and
$\xi|_{\mathfrak u_P}=\psi_L|_{\mathfrak u_P}=0$. Thus, we have $(\mathrm{Ad}^*(u)\xi)|_{\mathfrak l}=\xi|_{\mathfrak l}$.

\subsection{Local trivializations}
For each $g\in G$, we define a sub-variety $\Omega_g:=g\,\bar{U}TU/U\subset X$, where $\bar{U}$ denotes the opposite maximal unipotent subgroup. It is open and affine, with an identification: $\Omega_g\cong \bar{U}\times T$. Moreover, these open affine sub-varieties $\Omega_g$'s cover $X$.
The natural inclusion, $\sigma_g\colon\Omega_g\to G$, defined by $nt\,U\mapsto gnt$, gives a section of the projection $p: G \to X$. Therefore, we have:
\[
\Omega\times U\xrightarrow{\ \sim\ }p^{-1}(\Omega),\qquad (v,u)\mapsto\sigma_g(v)u
\]
for every open sub-variety $\Omega\subset\Omega_g$. Under this isomorphism right translation by $U$ becomes right multiplication on the
second factor, and $\xi_x$ ($x\in\mathfrak u$) becomes $1\otimes\xi^U_x$, the left-invariant vector field on $U$. We call such $\Omega$ \emph{trivialized}.

\begin{proposition}\label{prop:D-structure} We have following properties of $\mathcal D$:
\leavevmode
\begin{enumerate}
\item For trivialized affine $\Omega$, the map
\[
\theta_\Omega\colon D_\hbar(\Omega)\to\Gamma(\Omega,\mathcal D),\qquad P\mapsto(P\otimes1)\one,
\]
is an isomorphism of $\mathbb C[\hbar]$-algebras. Here $D_\hbar(\Omega)$ is the Rees algebra of the order filtration. In particular $\mathcal D$ is $\hbar$-torsion-free, and
$\mathcal D/(\hbar-a)\mathcal D$ is a TDO on $X$ for $a\neq0$ (with twist $\psi_L/a$).
\item There is a canonical isomorphism between sheaves of commutative $\mathcal O_X$-algebras
\[
\mathcal D/\hbar\mathcal D\ \xrightarrow{\ \sim\ }\ \pi_*\mathcal O_{T^*(G/U,\psi_L)} .
\]
\item Under (2), the right action of $z\in Z_\hbar(\mathfrak l)$ reduces to multiplication by $\pi_L^*(\bar z)$, where
$\bar z$ is the image of $z$ in $Z_\hbar(\mathfrak l)/\hbar=\SymAlg(\mathfrak l)^L=\mathcal O(\mathfrak c_L)$.
\end{enumerate}
\end{proposition}

\begin{proof}
(1) Let $\Omega$ be trivialized and affine. Then $D_\hbar(p^{-1}\Omega)=D_\hbar(\Omega)\otimes_{\mathbb C[\hbar]}D_\hbar(U)$, and
$D_\hbar(U)=\mathcal O(U)\otimes U_\hbar(\mathfrak u)$ via left-invariant vector fields. Since $\psi_L$ vanishes on
$[\mathfrak u,\mathfrak u]$, it extends to an algebra map $U_\hbar(\mathfrak u)\to\mathbb C[\hbar]$, and
$U_\hbar(\mathfrak u)=\mathbb C[\hbar]\oplus U_\hbar(\mathfrak u)\mathfrak u^{\psi_L}$. Hence:
\[
\mathcal M(p^{-1}\Omega)=D_\hbar(\Omega)\otimes_{\mathbb C[\hbar]}\bigl(D_\hbar(U)/D_\hbar(U)\mathfrak u^{\psi_L}\bigr)
=D_\hbar(\Omega)\otimes_{\mathbb C}\mathcal O(U),\qquad P\otimes f\xi^I\mapsto \psi_L(\xi^I)\,P\otimes f .
\tag{$\ast$}
\]
Since $r_u(f\one)=r_u(f)\one$, the group $U$ acts on the right-hand side of $(\ast)$ only through right
translation on $\mathcal O(U)$. So by taking $U$-invariants, we have:
\[
\mathcal M(p^{-1}\Omega)^U=D_\hbar(\Omega)\otimes\mathcal O(U)^U=D_\hbar(\Omega)\otimes1
\]
which is exactly 
the image of $\theta_\Omega$. Multiplicativity is clear, and the same computation applies to every open subset of $\Omega$.
In particular $\mathcal D|_\Omega\cong\mathcal D_{\Omega,\hbar}$, which is $\hbar$-torsion-free; for $a\ne0$ the specialization is
locally isomorphic to $\mathcal D_\Omega$. \footnote{changing the section $\sigma_g$ to $\sigma_g\cdot\phi$, with $\phi\colon \Omega\to U$, where $\phi$ is a $U$-valued function on $\Omega$, has the effect of changing $\theta_\Omega$ by conjugation with the exponential function: $e^{\psi_L(\phi)/\hbar}$. It twists $\calD^{\leq 1}$  by: $\hbar\partial\mapsto\hbar\partial+\partial(\psi_L\circ\phi)$.
Here $\psi_L\colon U\to\mathbb G_a$ is the character \eqref{eq:psi-def}. This follows from the usual description of a TDO by exact transition
$1$-forms.}

(2) Since $D_\hbar(G)/\hbar=\mathcal O(G\times\mathfrak g^*)$, under the left trivialization of $T^*G$, we have:
\[
\mathcal M/\hbar\mathcal M=\mathcal O(G)\otimes \SymAlg(\mathfrak g)/(x-\psi_L(x))_{x\in\mathfrak u}=\mathcal O\bigl(G\times(\psi_L+\mathfrak u^\perp)\bigr),
\]
and $U$ acts by $(g,\xi)u=(gu,\mathrm{Ad}^*(u^{-1})\xi) = (gu, u^{-1}\xi u)$. Therefore,
$(p_*(\mathcal M/\hbar))^U=\pi_*\mathcal O_{T^*(G/U,\psi_L)}$. This gives a canonical homomorphism between sheaves:
$\mathcal D/\hbar\mathcal D\to\pi_*\mathcal O_{T^*(G/U,\psi_L)}$. We need to check that it is an isomorphism over trivialized affine $\Omega$. Over such $\Omega$, we have the following short exact sequence:
\[
0\to\mathcal M\xrightarrow{\hbar}\mathcal M\to\mathcal M/\hbar\to0
\]
each term above, by $(\ast)$, is of
the form $N\otimes\mathcal O(U)$, with $U$ acting only on $\mathcal O(U)$. Therefore, we have the desired vanishing:
\[
H^1(U,N\otimes\mathcal O(U))=N\otimes H^1(U,\mathcal O(U))=0
\]
Therefore, taking $U$-invariants is exact, and it gives:
\[
\Gamma(\Omega,\mathcal D)/\hbar=(\mathcal M/\hbar)(p^{-1}\Omega)^U=\mathcal O(\pi^{-1}\Omega)
\]

(3) Modulo $\hbar$, the right multiplication by $z\in U_\hbar(\mathfrak l)$ on $\mathcal M$ becomes the multiplication by the
symbol $\bar z\in \SymAlg(\mathfrak l)^L\subset \SymAlg(\mathfrak g)=\mathcal O(\mathfrak g^*)$. In other words, multiply by the function: 
$(g,\xi)\mapsto\bar z(\xi|_{\mathfrak l})=\pi_L^*(\bar z)([g,\xi])$.
\end{proof}

\subsection{Proof of Theorem~\ref{thm:alg-glue}}\label{sec:hartogs-proof}
Denote $Z:=Z_\hbar(\mathfrak l)$, $Z_\alpha,Z_\beta,Z_{rs}$ for its localizations as in Definition~\ref{def:loci}, and for a sheaf, or a module, $E$ of right $Z$-modules, we write $E_\alpha:=E\otimes_ZZ_\alpha$, etc.

\begin{lemma}\label{lem:local-hartogs}
Let $\Omega$ be a trivialized affine open subset of $X$ and $E:=\Gamma(\Omega,\mathcal D)$. Then $E$ is torsion-free with
respect to the walls, and $E\to E_\alpha\times_{E_{rs}}E_\beta$ is an isomorphism.
\end{lemma}

\begin{proof}
Identify $E=D_\hbar(\Omega)$ via $\theta_\Omega$, and $E/\hbar E=\mathcal O(\pi^{-1}\Omega)=:\mathsf B$ via
Proposition~\ref{prop:D-structure}(2); write $\sigma\colon E\to \mathsf B$ for the reduction. Let $G_dE$ be the span of the
homogeneous components of order degree $\le d$ of the Rees algebra $D_\hbar(\Omega)$. \footnote{here $\hbar$ and vector fields have
degree $1$; this is half the grading of \S\ref{sec:notation}.} Then:
\begin{enumerate}
\item[(a)] $G_{-1}E=0$, $\bigcup_dG_dE=E$, and $G_dE\cap\hbar E=\hbar\,G_{d-1}E$.
\item[(b)] Using $\sigma_g$, identify $\pi^{-1}\Omega\cong \Omega\times(\psi_L+\mathfrak u^\perp)$, so that
$\mathsf B=\mathcal O(\Omega)[y_1,\cdots,y_n]$, where $y_i$'s are coordinates on the affine fiber. The trivialization of $T^*(p^{-1}\Omega)$ restricts to a fiberwise affine
isomorphism $T^*\Omega\cong\pi^{-1}\Omega$, so $\sigma(G_dE)=\mathsf B_{\le d}$, the functions of fiber degree $\le d$.
\item[(c)] If $z\in Z$ has $U_\hbar$-degree $\le e$ \footnote{with $\mathfrak l$ and $\hbar$ in degree $1$}, then
$(G_dE)\,z\subset G_{d+e}E$. Indeed, $(P\otimes1)z$ has order degree $\le d+e$ in $D_\hbar(p^{-1}\Omega)$, and the
projection $(\ast)$ does not increase degree. Moreover $\sigma(Pz)=\sigma(P)\,\pi_L^*(\bar z)$ by
Proposition~\ref{prop:D-structure}(3).
\item[(d)] Let $f\in Z$ be a finite product of walls, $\check\gamma+n\hbar$, of total degree $e$. Then its image:
$\bar f\in Z/\hbar Z=\SymAlg(\mathfrak t)^{W_L}=\mathcal O(\mathfrak c_L)$, is a product of coroots, homogeneous of degree $e$.

Over $\pi^{-1}\Omega\cong \Omega\times(\psi_L+\mathfrak u^\perp)$, the function $\pi_L^*\bar f$ is: $\xi\mapsto\bar f(\chi_L(\xi|_{\mathfrak l}))$, which is 
independent of the $\Omega$-coordinate. Its top-degree part is : $\xi_0\mapsto\bar f(\xi_0|_{\mathfrak t})$ on $\mathfrak u^\perp$
\footnote{invariant polynomials on $\mathfrak l^*$ restricted to $\mathfrak u_L^\perp$ only see the $\mathfrak t^*$-component}.
This is nonzero. Since $\mathcal O(\Omega)$ is a domain, we have: $\deg(b\,\pi_L^*\bar f)=\deg b+e$, for $b\in \mathsf B$.
\end{enumerate}
Using the above listed properties, it's easy to deduce $E$ is torsion free. Base on that we will prove surjectivity by an induction argument.
\\ \hspace{\fill} \\
\textbf{Torsion-freeness}. If $P\in G_dE$ and $Pf=0$, then $\sigma(P)\pi_L^*\bar f=0$, so $\sigma(P)=0$ , because $\mathsf B$ is a
domain, and $P=\hbar P_1$ with $P_1\in G_{d-1}E$ and $P_1f=0$. By induction on $d$, $P=0$.
\\ \hspace{\fill} \\
\textbf{Surjectivity}. Torsion-freeness reduces us to showing: if $P,Q\in E$, and $Pf_\beta=Qf_\alpha$, where
$f_\alpha$, (resp.\ $f_\beta$,) is a product of $\alpha$-walls, (resp.\ $\beta$-walls,) then
$P\in Ef_\alpha$. We argue by induction on $d$, with $P\in G_dE$. Then
by (c), $\sigma(P)\bar f_\beta=\sigma(Q)\bar f_\alpha$ in $\mathsf B$ , where we denote $\bar f$ for $\pi_L^*\bar f$.

The common zero locus of $\bar f_\alpha,\bar f_\beta$ in $\mathfrak c_L$ is the image of
$\bigcup_{\alpha,\beta}\{\check\alpha=0\}\cap\{\check\beta=0\}\subset\mathfrak t^*$; since $\check\alpha\neq\pm\check\beta$
it has codimension $2$. Because $\pi_L$ is flat (indeed, smooth, by \cite[Lemma 3.2.4]{macerato2023leviequivariantrestrictionsphericalperverse}), its preimage in the smooth variety
$\pi^{-1}\Omega$ has codimension $2$. So $\bar f_\alpha,\bar f_\beta$ have no common prime factor in any local ring, which is a UFD because of smoothness. Thus, $r:=\sigma(P)/\bar f_\alpha$ is regular on $\pi^{-1}\Omega$, i.e.\ $r\in \mathsf B$. By (d), $r\in \mathsf B_{\le d-e_\alpha}$.
By (b), we can choose $R\in G_{d-e_\alpha}E$, such that $\sigma(R)=r$.

By (c), $P-Rf_\alpha\in G_dE$ has symbol $0$, so by (a) we have: $P=Rf_\alpha+\hbar P_1$,  with some $P_1\in G_{d-1}E$.
Also, because $\sigma(Q-Rf_\beta)\,\bar f_\alpha=\sigma(Q)\bar f_\alpha-\sigma(P)\bar f_\beta=0$, we deduce that $Q=Rf_\beta+\hbar Q_1$.
Using the fact that $Z$ is commutative, we have:
\[
\hbar(P_1f_\beta-Q_1f_\alpha)=Pf_\beta-Qf_\alpha=0,
\]
Thus, $P_1f_\beta=Q_1f_\alpha$. By induction, $P_1\in Ef_\alpha$, and therefore $P\in Ef_\alpha$.
\end{proof}

\begin{proof}[Proof of Theorem~\ref{thm:alg-glue}]
The trivialized affine opens form a basis, which is closed under finite intersections. Cover $X$ by finitely many such opens $\Omega_i$. Then $\Omega_{ij}=\Omega_i\cap \Omega_j$ are again trivialized and
affine. The sequence
\[
0\to\Gamma(X,\mathcal D)\to\prod_i\Gamma(\Omega_i,\mathcal D)\to\prod_{i,j}\Gamma(\Omega_{ij},\mathcal D)
\]
is exact. Since $Z_\alpha,Z_\beta,Z_{rs}$ are flat over $Z$, it is still exact after applying $-\otimes_ZZ_\alpha$, etc.
In particular, $\Gamma(X,\mathcal D)_\alpha$ is the equalizer of the $\Gamma(\Omega_i,\mathcal D)_\alpha$, and similarly for
$\beta$ and $rs$. Since fiber products commute with equalizers, and by Lemma~\ref{lem:local-hartogs}, the map of
Theorem~\ref{thm:alg-glue} is an isomorphism on each $\Omega_i$ and $\Omega_{ij}$. Hence it is an isomorphism on $X$.
\end{proof}

\section{Comparison on the regular semisimple locus}\label{sec:rs-compat}
In this section, we will compare the isomorphisms $\widetilde \Psi_{\alpha}$ and $\widetilde \Psi_{\beta}$ we constructed in case 1 and 2, over the regular semisimple locus $\ars$. If they agree with each other, we can glue them together, by defining $\widetilde \Psi := \widetilde \Psi_{\alpha} \times_{\widetilde \Psi_{rs}} \widetilde \Psi_{\beta}$. Using the Hartogs property we proved in the previous section, it is not hard to see $\widetilde \Psi$ becomes the desired isomorphism in \ref{thm:main}. Unfortunately, this is not the case, and we have to make clear why $\widetilde \Psi_{\alpha}$ and $\widetilde \Psi_{\beta}$ differ over $\ars$, and what is that difference. As the first step towards finding the correction, we will explain why they differ in this section.

Although over $\ars$ both $\widetilde\Psi_{V,\alpha}$ and $\widetilde\Psi_{V,\beta}$ have the same image: $\H_A(\GrL,i_L^!\IC_V)_\rs\cong\bigoplus_\lambda V_\lambda\otimes\Crs$, when one goes around such a loop between algebraic side and geometric side, the composition is not identity, which means it won't come back to itself.
One underlying reason for this is, we used the isomorphism in \cite{bezrukavnikov2026equivariantsatakecategorykostantwhittaker} to construct $\widetilde \Psi_{\alpha}$, and the isomorphism in \cite{ginzburg2014differentialoperatorsguaffine} to construct $\widetilde \Psi_{\beta}$. However, the generic splittings of \cite{bezrukavnikov2026equivariantsatakecategorykostantwhittaker} and \cite{ginzburg2014differentialoperatorsguaffine} come from different filtrations. At the end of this section,  we will reduce the comparison to the group $L$ alone (Proposition~\ref{prop:reduction}).

\subsection{The maps involved}\label{sec:cmp-maps}
First, we give names to the homomorphisms that will be used frequently below. Fix $V \in \mathrm{Rep}(G)$, $\calF=\IC_V$ and $\calF_L:=\Res^{\cG}_{\cL}\calF$. All maps below are linear over $\Sh$. Since $\Phi$ is only defined after inverting $\hbar$ (Remark~\ref{rem:KK-asymptotic}), in this section the subscripts $\beta$ and $\rs$ mean localization over $\Cbe[\hbar^{-1}]$, resp.\ $\Crs[\hbar^{-1}]$, whenever $\Phi$ is involved; all modules are torsion-free, so nothing is lost.
\begin{itemize}
\item $\GR_G\colon(V\otimes\bM_\hbar)^{\fraku}\xrightarrow{\sim}\H_A(\GrT,i_T^!\calF)$, the isomorphism \eqref{eq:GR}.
\item $\GR_L\colon(V|_L\otimes\bM_{L,\hbar})^{\fraku_L}\xrightarrow{\sim}\H_A(\GrT,(i^L_T)^!\calF_L)$, the isomorphism \eqref{eq:GR} for $\cL$.
\item $\BF_L\colon\Sh\otimes_{\Zh(\frakl)}(V|_L\otimes\calW_{L,\hbar})^{\fraku_L,\psiL}\xrightarrow{\sim}\H_A(\GrL,\calF_L)$, the isomorphism \eqref{eq:BF} for $\cL$ together with Lemma~\ref{lem:formal}.
\item $\Phi^G:=\Phi\colon(V\otimes\bM_\hbar)^{\fraku}_\beta\xrightarrow{\sim}\Sh\otimes_{\Zh(\frakl)}(V\otimes\bM_\hbar(\psiL))^{\fraku,\psiL}_\beta$, the isomorphism \eqref{eq:case2-alg}, given by \eqref{eq:Phi-def}.
\item $\Phi^L\colon(V|_L\otimes\bM_{L,\hbar})^{\fraku_L}_\beta\xrightarrow{\sim}\Sh\otimes_{\Zh(\frakl)}(V|_L\otimes\calW_{L,\hbar})^{\fraku_L,\psiL}_\beta$, the same map for the pair $(\frakl,\frakl)$, built from the same Whittaker vector $\wv_L$.
\item $\mathrm{pr}_L\colon(V\otimes\bM_\hbar)^{\fraku}\to(V|_L\otimes\bM_{L,\hbar})^{\fraku_L}$, induced by $\id_V\otimes\mathrm{pr}$, where $\mathrm{pr}$ is the projection to the top $\frakz$-weight component; and its Verma--Whittaker analogue $\mathrm{pr}_{\mathcal{W}}$ of \eqref{eq:prV}, the map used to build $\widetilde\Psi_{V,\alpha}$.
\end{itemize}
By construction, after restriction to the regular semisimple locus,
\begin{equation}\label{eq:Psi-defs}
(\widetilde\Psi_{V,\beta})_\rs=\iota_L\circ\GR_G\circ\Phi^{-1},\qquad
(\widetilde\Psi_{V,\alpha})_\rs=\H_A(\Xi)^{-1}\circ\BF_L\circ\mathrm{pr}_{\mathcal{W}} .
\end{equation}
Here $\iota_L$ and $\H_A(\Xi)$ are the maps of Proposition~\ref{prop:geom-loc}(2) and (1).

\subsection{Three commuting squares}\label{sec:cmp-squares}
All the maps of \S\ref{sec:cmp-maps} fit into the following diagram of $\Crs$-modules, where everything is localized to $\ars$, we abbreviate $\H_A(-)=\H_A(-)_\rs$, and still denote $(i^L_T)^!\Xi$ by $\Xi$ by abusing the notation.
\begin{equation}\label{eq:big-diagram}
\begin{tikzcd}[column sep=3.4em,row sep=3.6em,cells={nodes={font=\footnotesize}}]
\H_A(\GrL,i_L^!\calF)\arrow[d,"\H_A(\Xi)"',"\wr"]
& \H_A(\GrT,i_T^!\calF)\arrow[l,"\iota_L"']\arrow[d,"\H_A(\Xi)"]\arrow[dl,phantom,"\mathrm{(I)}"]
& (V\otimes\bM_\hbar)^{\fraku}\arrow[l,"\GR_G"']\arrow[r,"\Phi^G"]\arrow[d,"\mathrm{pr}_L"]\arrow[dl,phantom,"\mathrm{(II)}"]
& (V\otimes\bM_\hbar(\psiL))^{\fraku,\psiL}\arrow[d,"\mathrm{pr}_{\mathcal{W}}"]\arrow[dl,phantom,"\mathrm{(III)}"]
\\
\H_A(\GrL,\calF_L)
& \H_A(\GrT,(i^L_T)^!\calF_L)\arrow[l,"\iota^L_T"]
& (V|_L\otimes\bM_{L,\hbar})^{\fraku_L}\arrow[l,"\GR_L"]\arrow[r,"\Phi^L"]
& (V|_L\otimes\calW_{L,\hbar})^{\fraku_L,\psiL}\arrow[lll,bend left=16,"\BF_L"]
\end{tikzcd}
\end{equation}
\begin{itemize}
\item The two isomorphisms to be compared are read off from \eqref{eq:Psi-defs}. $(\widetilde\Psi_{V,\beta})_\rs$ is the top row from right to left, $\iota_L\circ\GR_G\circ(\Phi^G)^{-1}$. $(\widetilde\Psi_{V,\alpha})_\rs$ goes down the right edge and back along the bottom arc, $\H_A(\Xi)^{-1}\circ\BF_L\circ\mathrm{pr}_{\mathcal{W}}$.
\item The squares (I), (II) and (III) commute. This is Lemmas~\ref{lem:geo-square}, \ref{lem:G2} and~\ref{lem:alg-square} below; none of them involves $\psiL$, Whittaker vectors or the isomorphism of \cite{bezrukavnikov2026equivariantsatakecategorykostantwhittaker}.
\item What remains is the bottom \emph{loop}: the bottom arc $\BF_L\circ\Phi^L$ against the bottom row $\iota^L_T\circ\GR_L$. This loop involves only the group $L$. It does \emph{not} commute, and correcting it is the subject of \S\S\ref{sec:cmp-reduction}--\ref{sec:discrepancy}.
\end{itemize}

\begin{lemma}[Algebraic square]\label{lem:alg-square}
On $(V\otimes\bM_\hbar)^{\fraku}_\rs$ we have $\mathrm{pr}_{\mathcal{W}}\circ\Phi^G=\Phi^L\circ\mathrm{pr}_L$.
\end{lemma}

\begin{proof}
By construction \eqref{eq:phi-G}, $\phi=U(\frakg)\otimes_{U(\frakp_L)}\phi_L$, that is, $\phi(y\otimes w)=y\cdot\phi_L(w)$ for $y\in U(\bar\fraku_P)$. Taking the top $\frakz$-weight component kills every $y$ of non-zero $\frakz$-weight. Hence $\mathrm{pr}\circ\phi=\phi_L\circ\mathrm{pr}$, where on the right $\mathrm{pr}$ is the projection $\widehat\bM\to\widehat\bM_L$ of the completed Verma module onto its top $\frakz$-weight piece. This projection is $\frakl$-linear, and so is $\id_V\otimes\mathrm{pr}$. Since $\wv_L=\sum_Jy_J\one\,c_J$ with $y_J\in U(\bar\fraku_L)\subset U(\frakl)$, and the coefficients act through $\frakt\subset\frakl$, we get $(\id_V\otimes\mathrm{pr})(f_\xi(\wv_L))=f_{\mathrm{pr}_L(\xi)}(\wv_L)$. The left hand side is $(V\otimes\phi_\hbar)\bigl(\mathrm{pr}_{\mathcal{W}}(\Phi^G(\xi))\bigr)$, and the right hand side is $(V\otimes\phi_{L,\hbar})\bigl(\Phi^L(\mathrm{pr}_L(\xi))\bigr)$.
\end{proof}

\begin{lemma}[Geometric square]\label{lem:geo-square}
On $\H_A(\GrT,i_T^!\calF)$ we have:
\begin{equation*}
    \H_A(\Xi)\circ \iota_L=\iota^L_T\circ\H_A\bigl((i^L_T)^!\Xi\bigr)
\end{equation*}
where $\H_A((i^L_T)^!\Xi)\colon\H_A(\GrT,i_T^!\calF)\to\H_A(\GrT,(i^L_T)^!\calF_L)$.
\end{lemma}

\begin{proof}
This is the naturality of the counit $\iota^L_T(i^L_T)^!\to\id$ with respect to $\Xi\colon i_L^!\calF\to\calF_L$.
\end{proof}

We write:
\begin{equation*}
    \Xi^G_T\colon i_T^!\to\Res^{\cG}_{\cT}, \quad \Xi^G_L=\Xi, \quad  \Xi^L_T\colon(i^L_T)^!\to\Res^{\cL}_{\cT}
\end{equation*}
for the natural transformations of \cite[(1.3.1)]{macerato2023leviequivariantrestrictionsphericalperverse} of the three pairs.

\begin{lemma}[Transitivity]\label{lem:Xi-trans}
Under the isomorphism $\Res^{\cL}_{\cT}\circ\Res^{\cG}_{\cL}\cong\Res^{\cG}_{\cT}$ one has
\begin{equation}\label{eq:Xi-trans}
\Xi^L_T\circ(i^L_T)^!\Xi^G_L=\Xi^G_T .
\end{equation}
\end{lemma}

\begin{proof}
This is a statement about sheaves only. It follows from the transitivity of hyperbolic localization and the definition of $\Xi$ as an adjunction map; compare \cite[\S1.3, \S2.4]{macerato2023leviequivariantrestrictionsphericalperverse}.
\end{proof}

\begin{lemma}[Ginzburg--Riche is compatible with hyperbolic restriction]\label{lem:G2}
\begin{equation}\label{eq:G2}
\H_A\bigl((i^L_T)^!\Xi\bigr)\circ\GR_G=\GR_L\circ\mathrm{pr}_L .
\end{equation}
\end{lemma}

\begin{proof}
We use the quantum form of \cite[Theorem 2.3.1]{ginzburg2014differentialoperatorsguaffine}, as restated in \cite[Theorem 1.4.4]{macerato2023leviequivariantrestrictionsphericalperverse}: for each $\mu\in\Lambda$, the isomorphism $\GR_G$ on the $\mu$-component is the \emph{unique} isomorphism compatible with
\begin{itemize}
\item the geometric map $\Xi^G_T\colon\H_A(\{t^\mu\},i_\mu^!\calF)\to\H_A(\{t^\mu\},\Res^{\cG}_{\cT}\calF)=V_\mu\otimes\Sh$.
\item the algebraic map $\Top\colon(V\otimes\bM_\hbar)^{\fraku}_\mu\to V_\mu\otimes\Sh$, the projection to the top Verma component  $V\otimes(\bM_\hbar)_0$, followed by the identification $(\bM_\hbar)_0=\Sh$.
\end{itemize}
In other words, we have $\Xi^G_T\circ\GR_G=\Top$. The map $\H_A(\Xi^L_T)\colon\H_A(\GrT,(i^L_T)^!\calF_L)\to V\otimes\Sh$ is injective, by Proposition~\ref{prop:geom-loc}(1) for the pair $(\cL,\cT)$, so it suffices to check \eqref{eq:G2} after composing with it. On the one hand, by Lemma~\ref{lem:Xi-trans},
\[
\Xi^L_T\circ(i^L_T)^!\Xi^G_L\circ\GR_G=\Xi^G_T\circ\GR_G=\Top .
\]
On the other hand, $\Xi^L_T\circ\GR_L\circ\mathrm{pr}_L=\Top_L\circ\mathrm{pr}_L$, and $\Top_L\circ\mathrm{pr}_L=\Top$. This is because $(\bM_\hbar)_0$ lies in the top $\frakz$-weight piece $\bM_{L,\hbar}$, and equals $(\bM_{L,\hbar})_0$.
\end{proof}

\subsection{Reduction to the Levi}\label{sec:cmp-reduction}
Let $(c_\nu)_\nu$, $c_\nu\in\Cbe^\times$, be invariant under translation along $\frakz^*$, and let $c$ act on $(V\otimes\bM_\hbar)^{\fraku}_\beta$ and $(V|_L\otimes\bM_{L,\hbar})^{\fraku_L}_\beta$ by $\xi\mapsto\xi\cdot c_\nu$ on the $\nu$-piece. Then $\mathrm{pr}_L\circ c=c\circ\mathrm{pr}_L$. Put
\[
\widetilde\Psi^{c}_{V,\beta}:=\eqref{eq:case2-geom}\circ\eqref{eq:GR}_\beta\circ c\circ\Phi^{-1}\colon(V\otimes\bM_\hbar(\psiL))^{\fraku,\psiL}\otimes_{\Zh(\frakl)}\Cbe\xrightarrow{\ \sim\ }\H_A(\GrL,i_L^!\IC_V)_\beta ,
\]
so that $\widetilde\Psi^{1}_{V,\beta}=\widetilde\Psi_{V,\beta}$.

In terms of \eqref{eq:big-diagram}, inserting $c$ replaces the bottom loop by the square
\begin{equation}\label{eq:loop}
\begin{tikzcd}[column sep=3.5em,cells={nodes={font=\footnotesize}}]
(V|_L\otimes\bM_{L,\hbar})^{\fraku_L}_\rs\arrow[r,"\Phi^L"]\arrow[d,"c"']
& (V|_L\otimes\calW_{L,\hbar})^{\fraku_L,\psiL}_\rs\arrow[d,"\BF_L"]
\\
(V|_L\otimes\bM_{L,\hbar})^{\fraku_L}_\rs\arrow[r,"\iota^L_T\circ\GR_L"']
& \H_A(\GrL,\calF_L)_\rs
\end{tikzcd}
\end{equation}
and the next proposition says that this square is all that has to be checked.

\begin{proposition}\label{prop:reduction}
$(\widetilde\Psi_{V,\alpha})_\rs=(\widetilde\Psi^c_{V,\beta})_\rs$ if and only if \eqref{eq:loop} commutes, i.e.
\begin{equation}\label{eq:G3}
\BF_L\circ\Phi^L=\iota^L_T\circ\GR_L\circ c\qquad\text{on }(V|_L\otimes\bM_{L,\hbar})^{\fraku_L}_\rs .
\end{equation}
\end{proposition}

\begin{proof}
Both sides of the claimed equality are isomorphisms over $\Crs$, so we may compose with $\Phi^G$ and with the isomorphism $\H_A(\Xi)_\rs$. By \eqref{eq:Psi-defs} the claim becomes $\H_A(\Xi)\circ \iota_L\circ\GR_G\circ c=\BF_L\circ\mathrm{pr}_{\mathcal{W}}\circ\Phi^G$ on $(V\otimes\bM_\hbar)^{\fraku}_\rs$. By Lemmas~\ref{lem:geo-square}, \ref{lem:G2} and~\ref{lem:alg-square}, we have:
\[
\H_A(\Xi)\circ\iota_L\circ\GR_G\circ c=\iota^L_T\circ\H_A((i^L_T)^!\Xi)\circ\GR_G\circ c=\iota^L_T\circ\GR_L\circ\mathrm{pr}_L\circ c=\iota^L_T\circ\GR_L\circ c\circ\mathrm{pr}_L,
\]
\[
\BF_L\circ\mathrm{pr}_{\mathcal{W}}\circ\Phi^G=\BF_L\circ\Phi^L\circ\mathrm{pr}_L .
\]
So \eqref{eq:G3} implies the claim. Conversely, $\mathrm{pr}_L$ becomes an isomorphism over $\Crs$ (Remark~\ref{rem:top-generic}), and all modules involved are torsion-free. Thus the claim implies \eqref{eq:G3}.
\end{proof}

The naive (before correction) choice $c=1$ doesn't work for $L\neq T$ (Theorem~\ref{thm:rank-one}); In Section~\ref{sec:discrepancy}, we will compute the correct $c$.

\section{The discrepancy for a reductive group}\label{sec:discrepancy}

\emph{In this section $G$ is any connected reductive group} . The result will be applied to our Levi subgroup $L$ under consideration. We will also denote Levi subgroups of $G$ in this section as "$L$". This is not to be confused with $L$ in previous sections.

\subsection{Set-up}\label{sec:disc-setup}

\paragraph{Top component and coinvariant class.}
Fix a basis $\{y_K\}$ of $\Uh(\ub)$ consisting of weight vectors, with $y_\emptyset=1$. Every $\xi\in V\otimes\bM_\hbar$ can be written uniquely as $\xi=\sum_Kv_K\otimes y_K\one$ with $v_K\in V\otimes\Sh$, the coefficients acting on the right. Put
\[
\Top(\xi):=v_\emptyset,\qquad\coinv(\xi):=\sum_KS(y_K)\,v_K\qquad(\in V\otimes\Sh),
\]
where $S$ is the antipode of $\Uh(\frakg)$:
\begin{equation*}
    S: U_\hbar(\mathfrak g) \to U_\hbar (\mathfrak g)^{\operatorname{op}}, \quad S(y_{(1)}y_{(2)}) = S(y_{(2)})S(y_{(1)})
\end{equation*}
Both maps are $\Sh$-linear. They extend to the completion $V\otimes\Mhat$, the product of the weight spaces; $\coinv$ stays a finite sum, because $S(y_K)$ kills $V$ when $y_K$ has large weight.

\paragraph{Weight pieces.}
$\sfA:=(V\otimes\bM_\hbar)^{\fraku}=\bigoplus_\nu\sfA_\nu$ is graded by relative weight (\S\ref{sec:left-action}); $\sfA_\nu$ is free over $\Sh$ of rank $m_\nu:=\dim V_\nu$, and $\Top(\sfA_\nu)\subset V_\nu\otimes\Sh$ \cite{ginzburg2014differentialoperatorsguaffine}, \cite[Lemma 7.3.3]{ginzburg2022differentialoperatorsgugelfandgraev}. On $\sfA_\nu$ put
\[
Y:=y+\hbar\rho,\qquad Y_\nu:=y+\hbar(\nu+\rho),
\]
\footnote{The shift by $\hbar\rho$ has a reason: The $\hbar$-dot action $w\cdot y=w(y+\hbar\rho)-\hbar\rho$ (\S\ref{sec:notation}) becomes linear in $Y$: $w\cdot y+\hbar\rho=w(Y)$. So $W$ permutes the functions $\langle Y,\check\gamma\rangle$, $\gamma\in R$, exactly as it permutes the coroots; this is used in \S\ref{sec:disc-weyl}. Moreover, $\Gamma$-functions below is better expressed when using $Y$ as its argument.} Indeed, by \eqref{eq:left-right}, $s\in\Sh$ acts on $\sfA_\nu$ on the right through $s(Y-\hbar\rho)=s(y)$ and on the left through $s(Y_\nu-\hbar\rho)=s(y+\hbar\nu)$. By Remark~\ref{rem:specialization} the set of walls is stable under translation by $\hbar\nu$, so localizing $\sfA_\nu$ along the right or along the left action gives the same module; we write $N_\rs$ for either.

\paragraph{The map $\Phi$.}
Let $\wv=\sum_Jy_J\one\,c_J\in(\Mhat)_\rs$ be the normalized Whittaker vector ($c_\emptyset=1$). For $\xi\in\sfA_\rs$ put
\begin{equation}\label{eq:Phi-G}
\Phi(\xi):=f_\xi(\wv)=\sum_Jy_J\,(c_J\cdot\xi)\ \in\ (V\otimes\Mhat)_\rs^{\fraku,\psi},
\end{equation}
with $c_J\cdot\xi$ the left action (Corollary~\ref{cor:comp-verma}(a)). This is \eqref{eq:Phi-def} for $L=G$, up to $V\otimes\phi_\hbar$. A priori $\Phi$, and hence the discrepancy $d$ and the map $\calE$ below, are defined over $\Crs[\hbar^{-1}]$; by Proposition~\ref{prop:Phi-integral} they are defined over $\Crs$ and invertible at the prime $(\hbar)$. We use this tacitly in this section.

\subsection{Top components and coinvariant classes}\label{sec:disc-class}

\begin{remark}[The reason for coinvariants]\label{rem:why-coinv}
The two sides see $\H_A(\GrG,\IC_V)$ through different coordinates. On the Ginzburg--Riche side the coordinate is $\Top(\xi)$, the costalk at $t^\nu$ ($\Xi^G_T\circ\GR=\Top$, proof of Lemma~\ref{lem:G2}). On the Bezrukavnikov--Finkelberg side, \eqref{eq:BF} is a statement over $\Zh(\frakg)=\H_{\GO}(\pt)$, and passing to $\Sh=\H_A(\pt)$ turns Whittaker invariants into $\ub$-coinvariants (Lemma~\ref{lem:GK}): Whittaker invariants are sections over $\frakt^*/W$, coinvariants their pull-back to $\frakt^*\times\bbA^1$. So the coordinate there is the coinvariant class. Since $\Phi$ preserves it (Lemma~\ref{lem:class}(b)), the discrepancy $d$ compares $\coinv$ with $\Top$; Lemma~\ref{lem:sl2-A} shows that they differ already for $\mathfrak{sl}_2$.
\end{remark}

\begin{lemma}\label{lem:class}
\begin{enumerate}[label=(\alph*)]
\item The map $\coinv$ vanishes on $\ub\,(V\otimes\Mhat)_\rs$, and it induces an isomorphism $\ub\backslash(V\otimes\Mhat)_\rs\xrightarrow{\ \sim\ }V\otimes\Crs$ with inverse $v\mapsto[v\otimes\one]$.
\item For $\xi\in(V\otimes\bM_\hbar)^{\fraku}_\rs$ we have $\Phi(\xi)\equiv\xi\pmod{\ub\,(V\otimes\Mhat)_\rs}$. In particular $\coinv(\Phi(\xi))=\coinv(\xi)$.
\end{enumerate}
\end{lemma}

\begin{proof}
(a) Since $\Uh(\frakg)$ is cocommutative \footnote{In a cocommutative Hopf algebra, $\Delta = \Delta^{\operatorname{op}}$, thus the composition: $(m\otimes 1) \circ (1\otimes \Delta) \circ (S \otimes 1) \circ \Delta^{\operatorname{op}}$ collapses to $(m\otimes 1) \circ (S \otimes 1\otimes 1) \circ (\Delta \otimes 1) \circ \Delta = (\eta \epsilon \otimes 1)\circ \Delta = \eta \otimes \id$; for $y \in U_\hbar(\frakg)$, it sends $y$ to $1\otimes y$. In expression \eqref{eq:untwist} below, we are using Sweedler notation, $\Delta (y) = \sum y_{(1)} \otimes y_{(2)}$. },
\begin{equation}\label{eq:untwist}
v\otimes y\one=\textstyle\sum y_{(1)}\cdot\bigl(S(y_{(2)})v\otimes\one\bigr)\qquad(y\in\Uh(\ub)).
\end{equation}
Hence $\zeta\equiv\coinv(\zeta)\otimes\one$ modulo $\ub\,(V\otimes\bM_\hbar)_\rs$ for $\zeta\in(V\otimes\bM_\hbar)_\rs$. For $\zeta$ in the completion, write $\zeta=\sum_\mu\zeta_\mu$ as a sum of weight components. Each $\zeta_\mu$ is a finite sum, and $\zeta_\mu-\coinv(\zeta_\mu)\otimes\one=\sum_\gamma f_\gamma Z_{\mu,\gamma}$, where $\{f_\gamma\}$ is a basis of $\ub$ and $Z_{\mu,\gamma}$ has weight $\mu+\gamma$. For fixed $\gamma$ the sum $\sum_\mu Z_{\mu,\gamma}$ converges in the completion, so $\zeta-\coinv(\zeta)\otimes\one=\sum_\gamma f_\gamma\sum_\mu Z_{\mu,\gamma}\in\ub\,(V\otimes\Mhat)_\rs$. Conversely, $\coinv(f\zeta)=0$ for $f\in\ub$, because $S(fy)=S(y)S(f)$ and $S(f)=-f$.

(b) By \eqref{eq:Phi-G}, $\Phi(\xi)-\xi=\sum_{J\neq\emptyset}y_J\,(c_J\cdot\xi)$, because $c_\emptyset=1$. For $J\neq\emptyset$, $y_J$ lies in the augmentation ideal $\ub\,\Uh(\ub)$, so each term lies in $\ub\,(V\otimes\Mhat)_\rs$. Regrouping the (convergent) sum according to the first factor of $y_J$, as in (a), shows that the whole sum lies there too. Now apply (a).
\end{proof}

\begin{remark}\label{rem:model}
$V\otimes\phi_\hbar$, $v\otimes\one_{\calW}\mapsto v\otimes\wv$, induces
\[
J\colon\ub\backslash\bigl(V\otimes\calW_\hbar\otimes_{\Zh(\frakg)}\Crs\bigr)\longrightarrow\ub\backslash(V\otimes\Mhat)_\rs\cong V\otimes\Crs .
\]
The source has rank $|W|\dim V$ and splits over $\ars$ into parts indexed by $w\in W$ \footnote{comparing the left action of $\Sh$, which descends to $\bar\fraku$-coinvariants, with the right action}; the target is the part $w=1$. Since $J[v\otimes\one_{\calW}]=v+(\text{lower weights})$, $J$ is an isomorphism on the part $w=1$ and kills the others. Thus $V\otimes\Crs$ with coordinate $\coinv$ is a model of \eqref{eq:BF-torus} after generic splitting, and $J\circ\Phi=\coinv$ by Lemma~\ref{lem:class}(b).
\end{remark}

\begin{lemma}\label{lem:J-integral}
The map $J$ of Remark~\ref{rem:model} is defined on the lattice $\sfW$, i.e.\ $J(\sfW)\subset V\otimes\Crs$, and $J\colon\sfW\otimes_{\Sh}\Crs\to V\otimes\Crs$ is an isomorphism.
\end{lemma}

\begin{proof}
In the parametrization $v\otimes s\mapsto[v\otimes s(h)\one_{\calW}]$ of \S\ref{sec:disc-lattice},
$J[v\otimes s(h)\one_{\calW}]=\calJ(v\otimes s)=\sum_JS(y_J)v\cdot s(y+\hbar\eta_J)\,c_J$. Since $\Uh(\frakg)$ acts on $V$ by $\hbar$ times the action of $\frakg$, $S(y_J)v\in\hbar^{|J|}V$, and $\hbar^{|J|}c_J\in\Crs$ by Lemma~\ref{lem:integrality}; so $\calJ(v\otimes s)\in V\otimes\Crs$. The term $J=\emptyset$ is $v\otimes s$ and the others lie in lower weights, so $J$ is unitriangular for the filtrations by weights, whose graded pieces are $V_\mu\otimes\Crs$ on both sides; hence it is an isomorphism over $\Crs$.
\end{proof}

\paragraph{The Verma flag.}
Fix a total order on the weights of $V$ refining the dominance order. Let $\mathrm{Fil}^{\ge\mu}\subset V\otimes\bM_\hbar$ be the $\Uh(\frakg)$-submodule generated by $V_{\ge\mu}\otimes\one\,\Sh$, and define $\mathrm{Fil}^{>\mu}$ in the same way. This is the Verma flag: $\mathrm{Fil}^{\ge\mu}/\mathrm{Fil}^{>\mu}$ is free over $\Uh(\ub)\otimes\Sh$ on $V_\mu\otimes\one$. Its image in $\varphi(V)$ is the canonical filtration of \cite[\S4.3, Lemma 6]{bezrukavnikov2026equivariantsatakecategorykostantwhittaker}, with $\operatorname{gr}_\mu\varphi(V)\cong V_\mu\otimes\Sh$ via $[v\otimes\one]\mapsto v$.

\begin{lemma}\label{lem:flag}
For $\xi\in(V\otimes\bM_\hbar)^{\fraku}_\nu$ we have $\xi\in\mathrm{Fil}^{\ge\nu}$ and $\xi\equiv\coinv(\xi)\otimes\one\pmod{\mathrm{Fil}^{>\nu}}$. In particular, the image of $\mathrm{cl}(\xi)$ in $\operatorname{gr}_\nu\varphi(V)$ is $\coinv(\xi)$. In general this differs from $\Top(\xi)$ (Lemma~\ref{lem:sl2-A}).
\end{lemma}

\begin{proof}
Write $\xi=\sum_Kv_K\otimes y_K\one$, and let $-|y_K|$ be the weight of $y_K$. Then $v_K$ has weight $\nu+|y_K|$. By \eqref{eq:untwist},
\[
v_K\otimes y_K\one=\textstyle\sum y_{K(1)}\cdot\bigl(S(y_{K(2)})v_K\otimes\one\bigr),
\]
and $S(y_{K(2)})v_K$ has weight $\nu+|y_{K(1)}|\ge\nu$. Equality holds only for the term with $y_{K(1)}=1$, which contributes $S(y_K)v_K\otimes\one$.
\end{proof}

\begin{remark}\label{rem:top-generic}
Although they differ, $\Top$ and $\coinv$ are both injective on $\sfA_\nu$ and become isomorphisms $\sfA_\nu\otimes_{\Sh}\Crs\xrightarrow{\sim}V_\nu\otimes\Crs$. For $\Top$ this is \cite[\S7.4, (7.4.2)]{ginzburg2014differentialoperatorsguaffine} (their $\kappa^{\mathrm{alg}}_{V,\lambda}$ is $\Top$); equivalently, $\Top=\Xi^G_T\circ\GR$ (proof of Lemma~\ref{lem:G2}) with $\Xi^G_T$ an isomorphism over $\Crs$ by Proposition~\ref{prop:geom-loc}(1) for $(\cG,\cT)$. For $\coinv$: over $\Crs[\hbar^{-1}]$ we have $\coinv=J\circ\Phi$ (Lemma~\ref{lem:class}(b), Remark~\ref{rem:model}), with $J$ and $\Phi$ isomorphisms (Lemma~\ref{lem:J-integral}, Theorem~\ref{thm:full-faithful} for $L=G$); and $\coinv$ is an isomorphism at the prime $(\hbar)$ (proof of Proposition~\ref{prop:Phi-integral}, Step~2). Since $\Crs$ is a UFD in which $\hbar$ is prime, its determinant is a unit in $\Crs$. Likewise $\mathrm{pr}_L\colon\sfA\to\sfA^L$ is injective and an isomorphism over $\Crs$, because $\Top=\Top_L\circ\mathrm{pr}_L$ (proof of Lemma~\ref{lem:G2}).
\end{remark}

\subsection{The discrepancy and the lattice criterion}\label{sec:disc-lattice}
Let $\sfH:=\H_A(\GrG,\IC_V)$. We denote the composition: $\iota_T\circ\GR\colon\sfA\to\sfH$, again as $\GR$, by abusing the notation. It is injective, and becomes an isomorphism over $\ars$, $\sfW:=\ub\backslash(V\otimes\calW_\hbar)$, and $\BF\colon\sfW\xrightarrow{\sim}\sfH$ the isomorphism \eqref{eq:BF-torus}. $\sfW,\sfH$ are free of rank $\dim V$. View $\Phi$ as a map $\sfA_\rs\xrightarrow{\sim}\sfW_\rs$ via $V\otimes\phi_\hbar$. Over $\ars$ all three modules split by the weights $\nu$ of $V$, and $\GR(\sfA_{\nu,\rs})=\sfH_\nu$.

\begin{convention}\label{conv:BF}
The isomorphism \eqref{eq:BF} is stated in \cite{bezrukavnikov2026equivariantsatakecategorykostantwhittaker} with the opposite conventions (\S\ref{sec:conventions}). We fix its transport to our conventions, for $G$ and for each of its Levi subgroups, so that the following hold.
\begin{itemize}
\item[(T)] $\BF\circ\Phi$ maps $\sfA_{\nu,\rs}$ into $\sfH_\nu$, for every weight $\nu$. Equivalently, the left and right actions of $\Zh(\frakg)$ on $(V\otimes\calW_\hbar)^{\fraku,\psi}$ (\S\ref{sec:intertwiners}) are matched with the two actions on $\sfH$, just as \eqref{eq:GR} matches them on $\sfA$.
\item[(M)] $\BF$ is compatible with the multiplication: composition of Whittaker intertwiners (Definition~\ref{def:composition}) on the source, the convolution product on the target \cite[Theorem 6(b)]{bezrukavnikov2026equivariantsatakecategorykostantwhittaker}.
\end{itemize}
Two transports satisfying (T) differ by a $\Zh(\frakg)\otimes\Zh(\frakg)$-linear automorphism of $(V\otimes\calW_\hbar)^{\fraku,\psi}$; by Theorem~\ref{thm:DF}, whose proof uses only (T) and (M), this affects only the constant $t_G$.
\end{convention}

We assume Convention~\ref{conv:BF} from now on and define the \emph{discrepancy}
\[
d:=\GR^{-1}\circ\BF\circ\Phi\ \in\operatorname{Aut}_{\Crs}(\sfA_\rs),\qquad d=\bigoplus_\nu d_\nu,\quad d_\nu\in GL(\sfA_{\nu,\rs}).
\]
For a Levi subgroup $L\subset G$ we write $\sfA^L,\sfW^L,\sfH^L,\GR_L,d^L,\dots$ for the same objects for $L$ and $V|_L$. Thus condition \eqref{eq:G3} of Proposition~\ref{prop:reduction} says exactly $d^L=c$; the square \eqref{eq:loop} for $L$ is the right hand square below.

All the maps of this section are summarized in the diagram (over $\Crs$)
\begin{equation}\label{eq:disc-diagram}
\begin{tikzcd}[column sep=4em,row sep=3em]
V\otimes\Crs
& \sfW_\rs\arrow[l,"J"',"\sim"]\arrow[r,"\BF","\sim"']
& \sfH_\rs
\\
\sfA_\rs\arrow[u,"\coinv","\wr"']\arrow[r,equal]
& \sfA_\rs\arrow[u,"\Phi"',"\wr"]\arrow[r,"d"]
& \sfA_\rs\arrow[u,"\GR"',"\wr"]
\end{tikzcd}
\end{equation}
\begin{itemize}
\item The right square commutes by the definition of $d$.
\item The left square commutes by Lemma~\ref{lem:class}(b) and Remark~\ref{rem:model}: $J\circ\Phi=\coinv$.
\end{itemize}
Each module carries an \emph{integral structure} over $\Sh$, compatible with the decomposition into weights $\nu$:
\begin{itemize}
\item $\sfA_\nu\subset\sfA_{\nu,\rs}$ on the bottom row;
\item $\sfW\cap\sfW_\nu$ and $\sfH\cap\sfH_\nu$ in the middle and on the right of the top row;
\item on the left of the top row, the lattices $\calC_\nu$ and $\calJ_\nu$ in $V_\nu\otimes\Crs$, which are the images $\coinv(\sfA_\nu)$ and $J(\sfW\cap\sfW_\nu)$ (defined below).
\end{itemize}
The maps $\GR$, $\BF$ and $J$ respect the integral structures (Proposition~\ref{prop:lattice}(a)), and so does $\coinv$ by definition. The map $\Phi$ does not, and $d$ measures the failure. Proposition~\ref{prop:lattice}(b) expresses it through the two lattices on the left: $d_\nu=Q_\nu^{-1}\circ\coinv$, where $Q_\nu\colon\sfA_\nu\xrightarrow{\sim}\calJ_\nu$ is the composite along the integral maps.

For $v\in V$ and $s\in\Sh$, we put:
\[
\calJ(v\otimes s):=\coinv\bigl(v\otimes s(h)\cdot\wv\bigr)=\sum_JS(y_J)v\cdot s(y+\hbar\eta_J)\,c_J(y)\ \in V\otimes\Crs,
\]
where $\wv=\sum_Jy_J\one c_J$ is the normalized Whittaker vector and $\eta_J\in\bbZ_{\le0}R^+$ is the weight of $y_J$. The map $v\otimes s\mapsto[v\otimes s(h)\cdot\one_{\calW}]$ is a bijection $V\otimes\Sh\to\sfW$ (PBW and Lemma~\ref{lem:GK}), and $\calJ$ is the map $J$ of Remark~\ref{rem:model} composed with it. Put
\[
\calJ_\nu:=\calJ(V\otimes\Sh)\cap(V_\nu\otimes\Crs),\qquad\calC_\nu:=\coinv(\sfA_\nu)\subset V_\nu\otimes\Sh .
\]
Both are free of rank $m_\nu$ and span $V_\nu\otimes\Crs$.

\begin{proposition}[Lattice criterion]\label{prop:lattice}
\begin{enumerate}[label=(\alph*)]
\item $\BF^{-1}\circ \GR(\sfA_\nu)=\sfW\cap\sfW_\nu$, and $J$ maps $\sfW\cap\sfW_\nu$ isomorphically onto $\calJ_\nu$.
\item Let $Q_\nu:=J\circ\BF^{-1}\circ\GR\colon\sfA_\nu\xrightarrow{\sim}\calJ_\nu$. Then $d_\nu=Q_\nu^{-1}\circ\coinv$ on $\sfA_{\nu,\rs}$.
\item If $m_\nu=1$, let $c_\nu$ and $q_\nu$ generate $\calC_\nu$ and $\calJ_\nu$. Then $d_\nu=\varepsilon_\nu\cdot c_\nu/q_\nu$ for some $\varepsilon_\nu\in\bbC^\times$. A different transport of \eqref{eq:BF} satisfying (T) only changes the constant $\varepsilon_\nu$.
\end{enumerate}
\end{proposition}

\begin{proof}
(a) By \cite[Theorem 1.4.1]{ginzburg2022differentialoperatorsgugelfandgraev}, we have: $\GR(\sfA_\nu)=\iota_{T,\nu}\H_A(i_\nu^!\IC_V)=\sfH\cap\sfH_\nu$ . $\BF$ is an isomorphism of lattices. By (T), $\BF$ maps $\sfW_\nu$ onto $\sfH_\nu$, thus, $\BF(\sfW\cap\sfW_\nu)=\sfH\cap\sfH_\nu$. The map $J$ is injective on $\sfW$, and it is unitriangular for the weight filtration: $J[v\otimes\one_{\calW}]=v+{}$lower weights. It is $(\Uh(\frakg), \Crs)$-bilinear, therefore it maps $\sfW_\nu$ into the $\nu$-piece of $\ub\backslash(V\otimes\Mhat)_\rs\cong V\otimes\Crs$, which is $V_\nu\otimes\Crs$. Note that on $V_\nu\otimes\Crs$ the left action of $s\in\Sh$ is the right action of $s(y+\hbar\nu)$ (\S\ref{sec:left-action}). Thus, $J(\sfW\cap\sfW_\nu)=J(\sfW)\cap(V_\nu\otimes\Crs)=\calJ_\nu$.

(b) For $\xi\in\sfA_\nu$, Lemma~\ref{lem:class}(b) and Remark~\ref{rem:model} give $J(\Phi(\xi))=\coinv(\xi)$. Moreover,  $J(\BF^{-1}\circ \GR(d_\nu\xi))=Q_\nu(d_\nu\xi)$. Since $\BF\circ\Phi=\GR\circ d$, we get $\coinv(\xi)=Q_\nu(d_\nu\xi)$.

(c) $Q_\nu(\xi_\nu)$ generates $\calJ_\nu$, so it is $\varepsilon_\nu^{-1}q_\nu$ with $\varepsilon_\nu\in\Sh^\times=\bbC^\times$. A change of transport is an automorphism of $\sfW$ preserving $\sfW\cap\sfW_\nu$; on a free module of rank one it must be a constant.
\end{proof}

\subsection{Rank one computation}\label{sec:disc-rank-one}
In this subsection $\frakg=\mathfrak{sl}_2$, with standard basis $e,f,h$ ($\fraku=\bbC e$, $\bar\fraku=\bbC f$, $h=\check\alpha$), and $\psi(e)=:\psi\neq0$. We write $y$ for $\langle y,\check\alpha\rangle$, so that $\Sh=\bbC[y,\hbar]$. In this case, $\langle\rho,\check\alpha\rangle=1$ and $\langle\alpha,\check\alpha\rangle=2$, so we have $Y=y+\hbar$ and $Y_\nu=y+\hbar(\nu+1)$. In $\Uh(\mathfrak{sl}_2)$ the relations read $[h,e]=2\hbar e$, $[h,f]=-2\hbar f$, $[e,f]=\hbar h$.

Let $V=V^n$ with basis $v_j:=f^jv_0$ ($0\le j\le n$) for the classical action. Then
\begin{equation}\label{eq:sl2-V}
f\cdot v_j=\hbar\,v_{j+1},\qquad h\cdot v_j=\hbar(n-2j)\,v_j,\qquad e\cdot v_j=\hbar\,j(n-j+1)\,v_{j-1}.
\end{equation}
The weights of $V$ are $\nu=n-2k$, $0\le k\le n$, with $V_\nu=\bbC v_k$; so all multiplicities are one.

For $k\in\bbZ$ put
\[
P_k(z):=\prod_{m=1}^{k}(z+m\hbar)\quad(k\ge0),\qquad P_k(z):=\prod_{m=0}^{-k-1}(z-m\hbar)^{-1}\quad(k<0),
\]
so that $P_k(z)=\hbar^k\,\Gamma(1+z/\hbar+k)/\Gamma(1+z/\hbar)$ for all $k$.

\begin{lemma}\label{lem:sl2-verma}
In $\bM_\hbar$, for $m\ge0$,
\[
h\cdot f^m\one=f^m\one\cdot(y-2m\hbar),\qquad e\cdot f^m\one=m\hbar\,f^{m-1}\one\cdot\bigl(y-(m-1)\hbar\bigr).
\]
\end{lemma}

\begin{proof}
From $hf=f(h-2\hbar)$ we get $hf^m=f^m(h-2m\hbar)$, which gives the first formula because $h\one=\one\cdot y$. For the second, $[e,f^m]=\sum_{i=0}^{m-1}f^i[e,f]f^{m-1-i}=\hbar\sum_{i=0}^{m-1}f^ihf^{m-1-i}$. Moving $h$ to the right gives
\[
\hbar\sum_{i=0}^{m-1}f^{m-1}\bigl(h-2(m-1-i)\hbar\bigr)=m\hbar\,f^{m-1}\bigl(h-(m-1)\hbar\bigr).
\]
Apply this to $\one$ and use $e\one=0$.
\end{proof}

\begin{lemma}[Ginzburg-Riche's computation]\label{lem:sl2-A}
Let $\nu=n-2k$. The module $\sfA_\nu=(V^n\otimes\bM_\hbar)^{\fraku}_\nu$ is free of rank one over $\Sh$, generated by:
\[
\xi_\nu:=\sum_{j=0}^kv_j\otimes f^{k-j}\one\cdot b_j,\qquad
b_j:=(-1)^{k-j}\binom kj\frac{(n-j)!}{(n-k)!}\prod_{i=k-j}^{k-1}(y-i\hbar).
\]
Moreover
\[
\Top(\xi_\nu)=\prod_{i=0}^{k-1}(y-i\hbar)\cdot v_k,\qquad
\coinv(\xi_\nu)=\prod_{i=0}^{k-1}\bigl(y+(\nu+2+i)\hbar\bigr)\cdot v_k .
\]
Consequently $\calC_\nu=\prod_{i=0}^{k-1}\bigl(y+(\nu+2+i)\hbar\bigr)\cdot\Sh\,v_k$.
\end{lemma}
This is the rank-one case of the computation of $(V\otimes\bM_\hbar)^{\fraku}$ in \cite{ginzburg2014differentialoperatorsguaffine}; we record the proof here in order to fix normalizations.

\begin{proof}
\emph{The equations.} The vectors $v_j\otimes f^m\one$ form an $\Sh$-basis of $V\otimes\bM_\hbar$ (right coefficients), and $v_j\otimes f^m\one$ has relative weight $n-2j-2m$. So an element of relative weight $\nu$ is $\xi=\sum_{j=0}^kv_j\otimes f^{k-j}\one\cdot a_j$ with $a_j\in\Sh$, and $\xi\in\sfA_\nu$ if and only if $e\cdot\xi=0$. By \eqref{eq:sl2-V} and Lemma~\ref{lem:sl2-verma},
\[
e\cdot\xi=\sum_{j=1}^k v_{j-1}\otimes f^{k-j}\one\cdot\hbar\Bigl(j(n-j+1)\,a_j+(k-j+1)\bigl(y-(k-j)\hbar\bigr)\,a_{j-1}\Bigr).
\]
Hence $\xi\in\sfA_\nu$ if and only if
\begin{equation}\label{eq:sl2-rec}
a_j=-\frac{(k-j+1)\bigl(y-(k-j)\hbar\bigr)}{j\,(n-j+1)}\,a_{j-1}\qquad(1\le j\le k).
\end{equation}

\emph{Freeness.} The denominators in \eqref{eq:sl2-rec} are non-zero integers. So for every $a_0\in\Sh$ the recursion has a unique solution, and it lies in $\Sh$. Thus $\xi\mapsto a_0$ is an isomorphism $\sfA_\nu\xrightarrow{\sim}\Sh$. The $b_j$ satisfy \eqref{eq:sl2-rec}: indeed
\[
\frac{b_j}{b_{j-1}}=-\frac{\binom kj}{\binom k{j-1}}\cdot\frac{(n-j)!}{(n-j+1)!}\cdot\bigl(y-(k-j)\hbar\bigr)=-\frac{k-j+1}{j}\cdot\frac1{n-j+1}\cdot\bigl(y-(k-j)\hbar\bigr).
\]
Moreover $b_0=(-1)^kn!/(n-k)!$ is a non-zero constant, so $\xi_\nu$ generates $\sfA_\nu$.

\emph{Top component.} The component along $\bM_\hbar$-degree $0$ is the term $j=k$, so $\Top(\xi_\nu)=b_k\,v_k$ with $b_k=\prod_{i=0}^{k-1}(y-i\hbar)$.

\emph{Coinvariant class.} Here $S(f^{k-j})=(-f)^{k-j}$ acts on $v_j$ by $(-\hbar)^{k-j}v_k$, by \eqref{eq:sl2-V}. Hence
\[
\coinv(\xi_\nu)=\Bigl(\sum_{j=0}^k(-\hbar)^{k-j}b_j\Bigr)v_k=\Bigl(\sum_{j=0}^k\hbar^{k-j}\binom kj\frac{(n-j)!}{(n-k)!}\prod_{i=k-j}^{k-1}(y-i\hbar)\Bigr)v_k .
\]
The $j$-th summand has degree $j$ in $y$, so the sum is a polynomial of degree $k$ in $y$ with leading coefficient $1$. To deal with the last expression, recall the following identity:
\begin{equation}\label{eq:Chu-Vander}
(x+a)_k = \sum_{j=0}^{k} \binom{k}{j}(x)_j (a)_{k-j}
\end{equation}
where $(b)_p = b(b+1)\cdots (b+p-1)$ is the Pochhammer symbol. This type of identity is sometimes called Chu-Vandermonde identity. In our situation, take $\hbar=1$, because the sum is homogeneous in $(y,\hbar)$, and take $x = y-(k-1)$, $a = n-k+1$: then $(x)_j=\prod_{i=k-j}^{k-1}(y-i)$, $(a)_{k-j}=(n-j)!/(n-k)!$ and $x+a = y+(n-2k+2) = y + \nu +2$. Restoring $\hbar$, the last expression becomes $\prod_{i=0}^{k-1}(y+(\nu + 2+i)\hbar)$. Therefore, $\coinv(\xi_\nu) = \prod_{i=0}^{k-1}(y+(\nu+2+i)\hbar)$.
\end{proof}

\begin{lemma}[Whittaker coefficients]\label{lem:sl2-whittaker}
\begin{enumerate}[label=(\alph*)]
\item The normalized Whittaker vector of $\Mhat$ is
\[
\wv=\sum_{r\ge0}f^r\one\cdot c_r,\qquad c_r=\frac{\psi^r}{\hbar^r\,r!\,\prod_{i=0}^{r-1}(y-i\hbar)} .
\]
The denominators vanish on the walls $H_{\alpha,-i}$, $i\ge0$, and on $\hbar=0$; so $c_r\in\Crs[\hbar^{-1}]$ and $\hbar^rc_r\in\Crs$, as Lemma~\ref{lem:integrality} predicts.
\item For $v_j\in V$ and $s\in\Crs$,
\[
\calJ(v_j\otimes s)=\sum_{r\ge0}\frac{(-\psi)^r\,s(y-2r\hbar)}{r!\,\prod_{i=0}^{r-1}(y-i\hbar)}\;v_{j+r}\qquad(v_{j+r}:=0\text{ for }j+r>n).
\]
\end{enumerate}
\end{lemma}

\begin{proof}
(a) Write $\wv=\sum_rf^r\one\cdot c_r$. By Lemma~\ref{lem:sl2-verma}, $e\cdot\wv=\sum_{r\ge1}f^{r-1}\one\cdot r\hbar\,(y-(r-1)\hbar)\,c_r$. So $e\cdot\wv=\psi\wv$ is equivalent to $r\hbar\,(y-(r-1)\hbar)\,c_r=\psi\,c_{r-1}$ for $r\ge1$, and with $c_0=1$ this gives the formula.

(b) By definition $\calJ(v\otimes s)=\coinv\bigl(v\otimes s(h)\cdot\wv\bigr)$. By Lemma~\ref{lem:sl2-verma}, $s(h)\cdot f^r\one=f^r\one\cdot s(y-2r\hbar)$. Hence $v_j\otimes s(h)\cdot\wv=\sum_rv_j\otimes f^r\one\cdot s(y-2r\hbar)c_r$, whose coinvariant class is $\sum_rS(f^r)v_j\cdot s(y-2r\hbar)c_r$. Finally $S(f^r)v_j=(-f)^rv_j=(-\hbar)^rv_{j+r}$, and the powers of $\hbar$ cancel against those in $c_r$ (cf.\ Lemma~\ref{lem:J-integral}).
\end{proof}

\begin{lemma}[Solving for $\calJ_\nu$]\label{lem:sl2-J}
Let $\nu=n-2k$.
\begin{enumerate}[label=(\alph*)]
\item Let $s\in\Crs$. There is a unique element of $\calJ(V\otimes\Crs)\cap(V_\nu\otimes\Crs)$ with $v_k$-component $s\,v_k$, namely $s\,v_k=\calJ\bigl(\sum_{r=0}^{n-k}v_{k+r}\otimes t_r\bigr)$ with
\[
t_r(y)=\frac{\psi^r\,s(y-2r\hbar)}{r!\,\prod_{i=r-1}^{2r-2}(y-i\hbar)}\qquad(t_0=s).
\]
\item $\calJ_\nu=\prod_{m=2}^{n-k+1}(y+m\hbar)\cdot\Sh\,v_k$. In particular $\calJ_\nu$ does not depend on $\psi$.
\end{enumerate}
\end{lemma}

\begin{proof}
(a) By Lemma~\ref{lem:sl2-whittaker}(b), $\calJ(v_j\otimes s)\in v_j\otimes s+\sum_{m>j}v_m\otimes\Crs$. So $\calJ$ is unitriangular for the order of the basis $(v_j)$, hence injective. Write an element of $V\otimes\Crs$ as $\sum_{j=k}^{n}v_j\otimes s_j$, and suppose its image under $\calJ$ lies in $V_\nu\otimes\Crs$. Thus, the vanishing of the weight $-n,-n+2, \cdots, n-2k-2 = \nu -2$-components of the image gives us a system of equations. Moreover, the weight $\nu$-component of the image is exactly $s_k$.For convenience, we relabel $t_r:=s_{k+r}$. By Lemma~\ref{lem:sl2-whittaker}(b), the $v_{k+M}$-component of the image ($1\le M\le n-k$) is
\begin{equation}\label{eq:sl2-EM}
E_M:=\sum_{r=0}^{M}\frac{(-\psi)^{M-r}\;t_r\bigl(y-2(M-r)\hbar\bigr)}{(M-r)!\,\prod_{i=0}^{M-r-1}(y-i\hbar)} .
\end{equation}
The coefficient of $t_M$ in $E_M$ is $1$. So the system $E_1=\dots=E_{n-k}=0$ determines $t_1,t_2,\dots$ recursively from $t_0=s$, and the solution is unique. After solving the first three terms:
\begin{equation*}
    t_1(y) = \frac{\psi s(y-2\hbar)}{y},\quad t_2(y)=\frac{\psi^2 s(y-4\hbar)}{2(y-\hbar)(y-2\hbar)}, \quad t_3(y) = \frac{\psi^3 s(y-6\hbar)}{6(y-2\hbar)(y-3\hbar)(y-4\hbar)}
\end{equation*}
one can find the pattern, and guess the general expression should be:
\begin{equation*}
    t_r(y) = \frac{\psi^r s(y-2r\hbar)}{r! (y-(2r-2)\hbar)(y-(2r-1)\hbar)\cdots (y-(r-1)\hbar)}
\end{equation*}
It remains to check that the $t_r$ above does solve the equations. Substituting, the $r$-th term of $E_M$ becomes
\[
\frac{(-1)^{M-r}\psi^M\,s(y-2M\hbar)}{r!\,(M-r)!}\cdot\frac{1}{\prod_{i=0}^{M-r-1}(y-i\hbar)\,\prod_{i=2M-r-1}^{2M-2}(y-i\hbar)},
\]
using $\prod_{i=r-1}^{2r-2}\bigl(y-2(M-r)\hbar-i\hbar\bigr)=\prod_{l=2M-r-1}^{2M-2}(y-l\hbar)$. Multiplying by $M!\prod_{l=0}^{2M-2}(y-l\hbar)\big/\bigl(\psi^Ms(y-2M\hbar)\bigr)$, which does not depend on $r$, gives
\[
E_M=\sum_{r=0}^M(-1)^{M-r}\binom Mr\,p_M(r),\qquad p_M(r):=\prod_{l=M-r}^{2M-r-2}(y-l\hbar)=\prod_{l=0}^{M-2}\bigl(y-(M-r+l)\hbar\bigr).
\]
For fixed $y,\hbar$, $p_M$ is a polynomial in $r$ of degree $M-1$. The sum is its $M$-th finite difference of $p_M$, hence equals $0$.

(b) Since $\calJ$ is injective, $s\,v_k\in\calJ(V\otimes\Sh)$ if and only if the preimage found in (a) lies in $V\otimes\Sh$, i.e.\ if and only if $t_r\in\Sh$ for $0\le r\le n-k$. As $\psi\neq0$ and the linear forms $y-i\hbar$ are pairwise coprime, $t_r\in\Sh$ if and only if $\prod_{i=r-1}^{2r-2}(y-i\hbar)$ divides $s(y-2r\hbar)$. Substituting $y\mapsto y+2r\hbar$, this says that $\prod_{i=r-1}^{2r-2}\bigl(y+(2r-i)\hbar\bigr)=\prod_{m=2}^{r+1}(y+m\hbar)$ divides $s$. These conditions are nested in $r$, so they are equivalent to $\prod_{m=2}^{n-k+1}(y+m\hbar)\mid s$.
\end{proof}

\begin{theorem}[Lattices in rank one]\label{thm:rank-one}
Let $\frakg=\mathfrak{sl}_2$, $V=V^n$, $\nu=n-2k$ with $0\le k\le n$, and $\psi(e)\neq0$. Then
\[
\calC_\nu=\prod_{i=0}^{k-1}\bigl(y+(\nu+2+i)\hbar\bigr)\cdot\Sh\,v_k,\qquad
\calJ_\nu=\prod_{m=2}^{n-k+1}(y+m\hbar)\cdot\Sh\,v_k,
\]
and $d_\nu=\delta_\nu$ up to a constant in $\bbC^\times$, where:
\[
\delta_\nu:=P_\nu(Y)^{-1}=\hbar^{-\nu}\,\frac{\Gamma(1+Y/\hbar)}{\Gamma(1+Y_\nu/\hbar)},\qquad\text{i.e.}\quad
\delta_\nu=\prod_{m=1}^{\nu}(Y+m\hbar)^{-1}\ (\nu\ge0),\qquad\delta_\nu=\prod_{m=0}^{-\nu-1}(Y-m\hbar)\ (\nu<0).
\]
In particular $d_\nu$ in general won't be a constant. So $c=1$ does not satisfy \eqref{eq:G3} even for $L$ of semisimple rank one.
\end{theorem}

\begin{proof}
The formulas for $\calC_\nu$ and $\calJ_\nu$ are Lemmas~\ref{lem:sl2-A} and~\ref{lem:sl2-J}(b). By Proposition~\ref{prop:lattice}(c), $d_\nu=\varepsilon_\nu\,c_\nu/q_\nu$ with
\[
\frac{c_\nu}{q_\nu}=\frac{\prod_{m=\nu+2}^{n-k+1}(y+m\hbar)}{\prod_{m=2}^{n-k+1}(y+m\hbar)}\qquad(\text{note that  } \nu+2+(k-1)=n-k+1).
\]
If $\nu\ge0$, the numerator is a sub-product of the denominator and $c_\nu/q_\nu=\prod_{m=2}^{\nu+1}(y+m\hbar)^{-1}=\prod_{m=1}^{\nu}(Y+m\hbar)^{-1}$. (Note that $\rho =1$ in rank one case, so $Y = y + \hbar$.) If $\nu<0$, the denominator is a sub-product of the numerator and $c_\nu/q_\nu=\prod_{m=\nu+2}^{1}(y+m\hbar)=\prod_{m=0}^{-\nu-1}(Y-m\hbar)$. In both cases this can be described by a factorial-like function $P_{k}(z)$, defined as following:
\begin{equation*}
    P_k(z)= \left\{
\begin{aligned}
    & \prod_{m=1}^{k}(z+m\hbar), \quad k\geq 0\\
    & \prod_{m=0}^{-k-1} (z-m\hbar)^{-1}, \quad k<0
\end{aligned}
\right.
\end{equation*}
such that  $ c_\nu /q_\nu =P_\nu(Y)^{-1}$. We can also express $P_k$ as a ration between Gamma-functions:
\begin{equation*}
    P_k(z)=\hbar^k\Gamma(1+z/\hbar+k)/\Gamma(1+z/\hbar)
\end{equation*}
Thus, if we put $Y_\nu=Y+\hbar\nu$, it follows that:
\begin{equation*}
    \delta_\nu = \hbar^{-\nu} \Gamma(1+Y/\hbar)/\Gamma(1+Y_\nu /\hbar)
\end{equation*}
\end{proof}

\subsection{The \texorpdfstring{$\Gamma$}{Gamma}-gauge}\label{sec:disc-gauge}

\begin{definition}[$\Gamma$-gauge]\label{def:gauge}
Let $L\supseteq T$ be a Levi subgroup of $G$ with positive roots $R_L^+$. We define the Gamma-function for $L$ by:
\begin{equation*}
    \Gamma_L(y):=\prod_{\beta\in R_L^+}\hbar^{\langle y,\check\beta\rangle/\hbar}\,\Gamma\bigl(1+\langle y,\check\beta\rangle/\hbar\bigr), \quad y\in S_\hbar
\end{equation*}
It is a product of the usual Gamma-functions for each positive roots of $L$, which is designed to cover the discrepancy at every wall of $L$. Following the expression of $\delta_\nu$ in rank one case, for each weight $\nu$, $Y=y+\hbar \rho$ and $Y_\nu = y+\hbar(\rho +\nu)$, we put
\[
\delta^L_\nu(y):=\prod_{\beta\in R_L^+}P_{\langle\nu,\check\beta\rangle}\bigl(\langle y+\hbar\rho,\check\beta\rangle\bigr)^{-1}
=\frac{\Gamma_L(Y)}{\Gamma_L(Y_\nu)}
\]
For $\sfA=(V\otimes\bM_\hbar)^{\fraku}$, or $\sfA = (V|_L \otimes \bM_{L,\hbar})^{\fraku_L}$, we define $\calG_L$ to be the automorphism of $\sfA_\beta$, such that:
\begin{equation*}
    \calG_L : \sfA_\nu \longrightarrow \sfA_\nu,\quad \xi \mapsto \xi \cdot\delta_\nu^L (y)
\end{equation*}
where $y\in S_\hbar$ stands for the Verma parameter. We call $\calG_L$ the $\Gamma$-gauge for the Levi subgroup $L$. Later we will see that the $\Gamma$-gauge for $L$ plays the role of the crucial part of the correction we are looking for.
\end{definition}

For $\beta\in R_L$ one has $\langle\rho,\check\beta\rangle=\langle\rho_L,\check\beta\rangle$, so $\delta^L_\nu$ depends only on $L$, not on $G$. For $L$ of semisimple rank one, $\delta^L_\nu$ is exactly the $\delta_\nu$ of Theorem~\ref{thm:rank-one}. The family $(\delta^L_\nu)_\nu$ is of the type $(c_\nu)_\nu$ considered in \S\ref{sec:cmp-reduction}.

\begin{lemma}\label{lem:gauge}
\begin{enumerate}[label=(\alph*)]
\item $\delta^L_\nu\in\Cbe^\times$. Therefore, $\calG_L$ is an automorphism of $\sfA_\beta$, which commutes with the left and right actions of $\Sh$.
\item $\mathrm{pr}_L\circ\calG_L=\calG_L\circ\mathrm{pr}_L$ . Here the  $\mathrm{pr}_L$ refers to the projection to the top $\mathfrak{z}$-component for universal Verma module, i.e. $\mathrm{pr}_L\colon(V\otimes\bM_\hbar)^{\fraku}\to(V|_L\otimes\bM_{L,\hbar})^{\fraku_L}$. We abuse the notation here and below.
\item $\calG_L$ is conjugation by $R_\Gamma$, the right multiplication by $\Gamma_L(y+\hbar\rho)$ on $\bM_\hbar$: for $\xi\in\sfA$, $f_{\calG_L(\xi)}=R_\Gamma\circ f_\xi\circ R_\Gamma^{-1}$. In particular $\calG_L$ is compatible with composition of intertwiners, so it is an automorphism of the algebra $\Dh(G/U)\otimes_{\Sh}\Cbe$.
\end{enumerate}
\end{lemma}

\begin{proof}
(a) The zeros and poles of $\delta^L_\nu$ lie on walls $H_{\beta,m}$ with $\beta\in R_L$. (b) $\mathrm{pr}_L$ preserves the relative weight and commutes with the right action. (c) $\delta^L_\nu=\varphi(y)/\varphi(y+\hbar\nu)$ for $\varphi(y)=\Gamma_L(y+\hbar\rho)$, so this is the coboundary case of Corollary~\ref{cor:comp-verma}(c). Each factor of $\Gamma_L(Y)/\Gamma_L(Y_\nu)$ is rational, because $Y_\nu-Y=\hbar\nu$ is integral.
\end{proof}

\subsection{The discrepancy formula}\label{sec:disc-formula}
Let $\calG_G$ be the $\Gamma$-gauge of $G$ (Definition~\ref{def:gauge} with $L=G$); it acts on $\sfA_\nu$ by $\delta_\nu:=\delta^G_\nu$. Put
\[
\calE:=d\circ\calG_G^{-1}=\bigoplus_\nu\calE_\nu,\qquad\calE_\nu\in GL(\sfA_{\nu,\rs}).
\]

\begin{theorem}[Discrepancy formula]\label{thm:DF}
There is a unique $t_G\in T(\bbC)$ such that $\calE_\nu=\nu(t_G)$ on $\sfA_\nu(V)$ for all $V\in\Rep(G)$ and all weights $\nu$. That is, $d=t_G\circ\calG_G$.
\end{theorem}

We prove integrality of $\calE$ at simple walls via minimal Levi subgroups (\S\ref{sec:disc-simple}) and at the other walls via Weyl reflections (\S\ref{sec:disc-weyl}), both resting on rank one; multiplicativity then forces $\calE\in T$ (\S\ref{sec:disc-rigid}). Consider the following condition:
\begin{center}
    $\calJ_\nu=\delta_\nu^{-1}\,\calC_\nu$, for all $V$ and $\nu$ \quad (DF$_{\rm lat}$)
\end{center}

\begin{lemma}\label{lem:int}
Let $\frakq\subset\Sh$ be a prime of height one. If $\frakq$ is not a wall, $\calE_\nu$ is an automorphism of $\sfA_\nu\otimes(\Sh)_\frakq$. If $\frakq$ is a wall, $\calE_\nu$ is an automorphism of $\sfA_\nu\otimes(\Sh)_\frakq$ if and only if $(\calJ_\nu)_\frakq=(\delta_\nu^{-1}\calC_\nu)_\frakq$. Hence (DF$_{\rm lat}$) holds if and only if $\calE_\nu\in GL(\sfA_\nu)$ for all $\nu$.
\end{lemma}

\begin{proof}
By Proposition~\ref{prop:lattice}(b), $\calE_\nu=Q_\nu^{-1}\circ(\delta_\nu^{-1}\coinv)$, where $Q_\nu\colon\sfA_\nu\to\calJ_\nu$ and $\delta_\nu^{-1}\coinv\colon\sfA_\nu\to\delta_\nu^{-1}\calC_\nu$ are isomorphisms. If $\frakq$ is not a wall and $\frakq\neq(\hbar)$, then $\Crs[\hbar^{-1}]\subset(\Sh)_\frakq$; for $\frakq=(\hbar)$ use Proposition~\ref{prop:Phi-integral}. Finally $\sfA_\nu$ is free, so $\sfA_\nu=\bigcap_\frakq\sfA_\nu\otimes(\Sh)_\frakq$.
\end{proof}

\begin{lemma}\label{lem:rank-one-all}
If $G$ has semisimple rank one, (DF$_{\rm lat}$) holds for every $V\in\Rep(G)$.
\end{lemma}

\begin{proof}
$\calC_\nu$ and $\calJ_\nu$ are additive in $V$, so we may take $V$ irreducible. Then $\nu$ has multiplicity one, the centre of $G$ only shifts the parameter, and Theorem~\ref{thm:rank-one} applies to the derived $\mathfrak{sl}_2$; its $\delta_\nu$ is $\delta^G_\nu$.
\end{proof}

\subsection{Localize to simple roots}\label{sec:disc-simple}
Fix a simple root $\alpha$. Let $L\supset T$ be the minimal Levi subgroup with roots $\pm\alpha$, $\frakp=\frakl\oplus\fraku_P$ the standard parabolic subalgebra, $\bar\fraku_P$ the opposite nilradical and $\frakz=\Lie Z(L)$ \footnote{Warning: this $L$ is better denoted as $L_\alpha$. But to keep the notation simple, we choose not to add the subscript $\alpha$. It is not the Levi subgroup in Theorem~\ref{thm:main}.} . The goal of this section is to show, after localizing to a wall of $\alpha$, the two lattices $\calC_\nu$, $\calJ_\nu$ become the corresponding lattices $\calC_\nu^L$, $\calJ_\nu^L$ for the rank one subgroup $L$. Let $\frakq$ be a wall of $\alpha$, i.e. $\frakq = H_{\alpha,n}$. 

\subsection*{Localize $\calC$ lattice}
To begin with, we have the following by PBW and Frobenius adjunction:
\begin{equation*}
    V\otimes \bM_\hbar = V \otimes \Ind_{\frakp}^{\frakg}(\bM_{L,\hbar}) \cong \Ind_{\frakp}^{\frakg}(V\otimes \bM_{L,\hbar}) = U_\hbar(\bar{\fraku}_P)\otimes V \otimes \bM_{L,\hbar}
\end{equation*}
Under the last expression, $\pr_C$ is just the projection to the $U_\hbar(\bar{\fraku}_P)$-degree-0 component, which is $\cong V\otimes\bM_{L,\hbar}$. Equivariantly, it is quotient by $\bar{\fraku}_P$ 
\begin{equation*}
    \pr_C(\xi) = \xi\, \,\mathrm{mod} \,\,\bar{\fraku}_P (V\otimes \bM_{\hbar}), \quad v\otimes y_P m \mapsto S(y_P)v\otimes m
\end{equation*}
Taking $\bar{\fraku}$-coinvariants in two stages, first quotient by $\bar{\fraku}_P$ and then by $\bar{\fraku}_L$ gives:
\begin{equation}\label{eq:prC}
    \coinv = \coinv_L \circ \pr_C
\end{equation}
Therefore, if $\pr_C: \sfA_\nu \to \sfA_\nu^L$ becomes an isomorphism after localizing at $\frakq$, then:
\begin{equation*}
   ( \calC_\nu)_\frakq = \coinv_L (\pr_C (\sfA_\nu)_\frakq) = \coinv_L((\sfA_\nu^L)_\frakq) = (\calC_\nu^L)_\frakq
\end{equation*}
This is proved in Lemma~\ref{lem:locC} below.

\begin{lemma}\label{lem:locC}
Let $\frakq$ be a wall of $\alpha$. Then $\mathrm{pr}_C\colon\sfA_\nu\otimes(\Sh)_\frakq\to\sfA^L_\nu\otimes(\Sh)_\frakq$ is an isomorphism. Hence $(\calC_\nu)_\frakq=(\calC^L_\nu)_\frakq$.
\end{lemma}

\begin{proof}
\emph{Step 1: reducing to injectivity}. Notice that both $\sfA_\nu$ and $\sfA_\nu^L$ are free of the same rank $m_\nu = \dim V_\nu$, and $(S_\hbar)_\frakq$ is a DVR, $\pr_C$ is an isomorphism iff its determinant is a unit. That happens iff $\pr_C$ is injective after reducing modulo $\frakq$, i.e. over the residue field $k(\frakq)$. We can also view $k(\frakq)$ as a generic point $y \in H_{\alpha,n}$, and $\sfA_{\nu}\otimes k(\frakq)$ can be viewed as the singular vectors at that point. (By "singular vectors", we mean the elements killed by $\fraku$ inside a module.) Therefore, this Lemma amounts to show: at the generic point of the wall $\frakq = H_{\alpha,n}$, no non-zero singular vector $\xi$ of weight $\nu$ can lie in $\bar{\fraku}_P (V\otimes \bM_{\hbar})$.    

\emph{Step 2: the filtration}. We fix an element $h_0 \in \frakz$, such that $\gamma (h_0) >0$ for every positive simple root $\gamma \neq \alpha$, $\alpha (h_0)=0$. For a weight $\mu$, we define the \emph{degree} of $\mu$ by $\mu(h_0)$. It is easy to see that $\bar{\fraku}_P$ lowers the degree, $\fraku_P$ raise it, and $\frakl$ preseves it. Denoting $V^{(k)} \subset V$ be the sum of the weight spaces of degree $k$, we filter $V$ by degree. Since each $V^{(\geq k)} \otimes \bM_{L,\hbar}$ is $\frakp$-stable, after applying $\Ind_{\frakp}^{\frakg}$, it gives a filtration of $V\otimes \bM_\hbar$ by $U_\hbar(\frakg)$-submodules. The $k$-th graded piece is:
\begin{equation*}
    \Ind_\frakp^{\frakg} (V^{(k)} \otimes \bM_{L,\hbar}) = U_\hbar(\bar{\fraku}_P)\otimes V^{(k)}\otimes \bM_{L,\hbar}
\end{equation*}
Note that this graded piece has a Verma flag, whose highest weights are of the form $y+\hbar \mu$ with $\mu$ of degree $k$ \footnote{Here and below in the proof of Lemma~\ref{lem:locC}, $y$ stands for a generic point on the $\alpha$-wall.}. Now we take a non-zero singular vector $\xi \in \sfA_\nu \otimes k(\frakq)$, the weight of $\xi$ equals $y+\hbar \nu$. Let $k$ be the deepest filtration step containing $\xi$ (i.e.take the maximal $k$). Denote $\bar\xi$ the image of $\xi$ in the $k$-th graded piece, then $\bar\xi$ is also non-zero and singular. 

\emph{Step 3: linkage}. Now $\bar\xi$ is a  singular vector in a module with a Verma flag, it must have the same central character with one of the flag's Verma. The module in consideration is the $k$-th graded piece, $U_\hbar(\bar{\fraku}_P)\otimes V^{(k)}\otimes \bM_{L,\hbar}$, thus we deduce $\chi_{y+\hbar\nu} = \chi_{y+\hbar \mu}$ for some weight $\mu \in V^{(k)}$ of degree $k$. This implies at the generic point of the wall $H_{\alpha,n}$, $y+\hbar\nu$ and $y+\hbar\mu$ are linked by some element in the Weyl group. Therefore, we must have $y+\hbar \mu = w \cdot(y+\hbar \nu)$ for some $w\in W$, and being generic implies it holds on the whole hyperplane $H_{\alpha,n}$. This forces $w$ to fix $\ker(s_\alpha)$ pointwise, and thus $w=s_\alpha$. Therefore, $\nu - \mu \in \bbZ \cdot \alpha$, and $\nu(h_0) - \mu(h_0) =\alpha(h_0)\cdot \bbZ =0$, which implies $\nu$ has the same degree with $\mu$, which equals $k$. In other words, for $\xi \in \sfA_\nu \otimes k(\frakq)$, $\bar\xi$ lies exactly in the $\nu(h_0)$-th graded piece.

\emph{Step 4: degree count}. From there it should be not hard to see that $\xi$ cannot lie in $\bar\fraku_P (V\otimes \bM_{\hbar})$ generically over the $\alpha$ wall. Indeed, $\xi$'s image in $k$-th graded piece $ U_\hbar(\bar{\fraku}_P)\otimes V^{(k)}\otimes \bM_{L,\hbar}$ has degree $k$, but any term with non-trivial $U_\hbar(\bar\fraku_P)$-terms on the left has degree $<k$. So the $V^{(k)}$ component of $\pr_C(\xi)$, $\bar\xi$, is non-zero, therefore $\pr_C(\xi) \neq 0$.
\end{proof}

\subsection*{Localize $\calJ$ Lattice}
\begin{lemma}\label{lem:locJ}
Let $\frakq$ be a wall of $\alpha$. Then $(\calJ_\nu)_\frakq=(\calJ^L_\nu)_\frakq$.
\end{lemma}

\begin{proof}
\emph{Step 1: the same lattice}. When we restrict $\calW_\hbar$ to $\bar\frakp = \frakl \oplus \bar\fraku_P$, PBW tells us that it is induced from $\calW_{L,\hbar}$: $\calW_\hbar|_{\bar\frakp} \cong \Ind_{\frakl}^{\bar\frakp} \calW_{L,\hbar}$. Thus, taking $\bar\fraku_P$-coinvariants first, $\bar\fraku_L$-coinvariants next, gives:
\begin{equation*}
    \sfW = \bar\fraku_L\backslash \bar\fraku_P \backslash (V\otimes \calW_{\hbar}) = \bar\fraku_L \backslash (V\otimes \calW_{L,\hbar})= \sfW^L
\end{equation*}
with the same parametrization:
\begin{equation}\label{eq:parametrize}
    V\otimes S_\hbar \xrightarrow{\sim} \sfW = \sfW^L, \quad v\otimes s(y) \mapsto [v\otimes s(h)\cdot \one_\calW] = [v\otimes s(h)\cdot \one_{\calW_L}]
\end{equation}

\emph{Step 2: compare $J$ and $J^L$}. To compare the two homomorphisms $J$ and $J^L$, we need first to understand the relation between the universal Whittaker vectors. Since $\fraku_L$ preserves the $\Uh(\bar\fraku_P)$-degree in $\Mhat=\Uh(\bar\fraku_P)\otimes\widehat\bM_{L,\hbar}$, the degree-$0$ \footnote{This degree refers to the degree defined by $h_0$, in Lemma~\ref{lem:locC} .} part of the Whittaker vector $\wv$ is a $\psiL$-Whittaker vector of $\widehat\bM_{L,\hbar}$ with leading coefficient $1$, i.e.\ it is $\wv_L$:
\begin{equation*}
    \wv = \wv_L + (\text{terms contain } \bar\fraku_P)
\end{equation*}
 where the remaining terms have strictly lower degrees. Therefore:
\begin{equation*}
    J = J^L +(\text{degree-lowering}),\quad J\circ (J^L)^{-1} = \id + (\text{degree-lowering})
\end{equation*}

\emph{Step 3: two splittings and one filtration}. 

The first splitting, which we call $L$-splitting, comes from Step 1: when we view $\sfW = \sfW^L$, there is a natural splitting of $\sfW$ into degree pieces:
$\sfW = \oplus_k \sfW^{(k)}$, where $\sfW^{(k)}$ is the image of $V^{(k)} \otimes S_\hbar$ under \eqref{eq:parametrize}. Since it is defined over $S_\hbar$, this splitting is integral. 

The second splitting, which we call $G$-splitting, comes from the generic splitting for $G$ in \cite[Section 5]{bezrukavnikov2026equivariantsatakecategorykostantwhittaker}: over the generic locus, i.e. base change to $\bbK = \mathrm{Frac} (S_\hbar)$, $\sfW$ splits, and we denote by $\sfW_\mu \subset \sfW_\bbK := \sfW \otimes_{S_\hbar} \bbK$ the $\mu$-weight piece \footnote{Be aware that $\sfW_\nu$ is only rational, not integral, even after localizing to the $\alpha$ wall, for $G$ and also for $L$}. After we group weights $\mu$'s by its degree, gives: $\sfW^{G,(k)} = \oplus_{\mu, \mu(h_0) =k} \sfW_\mu $. Different from the first splitting, as a prior, it is only defined over $\bbK$.

These two splittings actually define the same filtration. Indeed, the filtration in the first case is induced by $\mathrm{Fil}^{\leq k} := (J^L)^{-1} (V^{(\leq k)}\otimes \bbK ) $; the second one by: $\mathrm{Fil}^{\leq k} := (J)^{-1} (V^{(\leq k)}\otimes \bbK )$. Since by Step 2, $J\circ (J^L)^{-1}$ is unipotent, we have $(J)^{-1} (V^{(\leq k)}\otimes \bbK ) = (J^L)^{-1} (V^{(\leq k)}\otimes \bbK )$. Let $\pi_k : \mathrm{Fil}^{\leq k} \to \sfW_{\bbK}^{(k)} $ be the projection, then we have: $J = J^L \circ \pi_k$ on $\sfW^{G,(k)}$. In other words, on the $G$-piece of degree $k$, $J$ is first project to $L$-piece, then apply $J^L$.

\emph{Step 4: $G$-splitting is integral at the wall}. Although the $G$-splitting in Step 3 is not integral, the crucial point is that, it becomes integral after localizing at a $\alpha$ wall. More precisely, let $\frakq' = H_{\alpha,m}$ \footnote{$\frakq'$ is a shift of $\frakq$ by $\hbar \nu$ when dealing with weight-$\nu$ piece. The reason is that, on $\bar\fraku \backslash V\otimes \widehat\bM_\hbar \cong V\otimes \bbC [\ars]$, the left action and the right action of $s\in S_\hbar$ differ by shifting by $\hbar \nu$: $s\cdot [v\otimes \one] = [v\otimes \one] \cdot s(y+\hbar\nu)$, see \S\ref{sec:left-action}. For $\xi \in \sfW_\nu$, $J(s\cdot \xi) = J(\xi) \cdot s(y+\hbar\nu)$. Therefore, if $\frakq = \{ y:\, \langle y,\check{\alpha} \rangle +n\hbar =0$ \}, then $\frakq' = \{ y:\, \langle y,\check{\alpha} \rangle + (n-\langle\nu,\check{\alpha\rangle)\hbar =0\}}$; in other words, $\frakq' = \frakq + \hbar\nu$. This is again a wall of $\alpha$, because $\langle\nu, \check{\alpha}\rangle \in \bbZ$. We have $(\calJ_\nu)_\frakq = J(\sfW_{\frakq'}\cap \sfW_{\nu})$.} be a $\alpha$ wall, then:
\begin{equation}\label{eq:Wsplit}
    \sfW_{\frakq'} = \oplus_{k} (\sfW_{\frakq'} \cap  \sfW^{G,(k)} )
\end{equation}
as modules over $(S_\hbar)_{\frakq'}$. The idea is simple: we can find idemponents in $(S_\hbar)_{\frakq'}$, projecting $\sfW_{\frakq'}$ onto each degree $k$ pieces. 

The action of $Z_\hbar(\frakg)$ on $\sfW_\nu$ is by central character $\chi_{y+\hbar\nu}$. Here $y$ stands for the coordinate on $S_\hbar$, for every $z \in Z_\hbar(\frakg)$, $\hc (z) \in S_
\hbar$, and $\chi_{y+\hbar\mu}(z) :=\hc (z)|_{y\mapsto y+ \hbar\mu} \in S_\hbar
$ is polynomial in the coordinate $y$. For each degree $k$, we define a $S_\hbar$-valued polynomial:
\begin{equation*}
    p_k(T):=\prod_{\mu(h_0)=k}(T-\chi_{y+\hbar\mu}(z))\in\Sh[T]
\end{equation*}
where $z \in Z_\hbar(\frakg)$ is to be chosen later. 

Note that $(\Sh)_{\frakq'}$ is a DVR, whose units are elements not inside $\frakq'$, i.e. polynomials that do not vanish identically on the hyperplane $\frakq' = H_{\alpha,m}$. For $\mu$ and $\nu$
of different degrees, i.e. $\mu(h_0) \neq \nu(h_0)$, we have seen from Step 3 in Lemma~\ref{lem:locC} that:
\begin{equation*}
    \chi_{y+\hbar\mu}(z) - \chi_{y+\hbar \nu}(z) \notin \frakq', \quad \text{for generic} \, z\in Z_\hbar(\frakg)
\end{equation*}
That implies $T-\chi_{y+\hbar\mu} (z)$ and $T-\chi_{y+\hbar\nu} (z)$ generates the unit ideal in $(\Sh)_{\frakq'}[T]$. Since the products of pairwise comaximal elements are again comaximal, we see that $p_k(T)$ and $p_{k'}(T)$ are pairwisely comaximal for $k\neq k'$.(We choose $z$ to be generic in $Z_\hbar(\frakg)$.) Therefore, by Chinese remainder theorem:
\begin{equation*}
    (\Sh)_{\frakq'} [T]/\prod_k p_k(T) \cong \prod_k (\Sh)_{\frakq'} [T]/p_k(T)
\end{equation*}
and there are polynomials $e_k (T) \in (\Sh)_{\frakq'} [T]$ decomposing $1$ into sum of idemponents: $1 =\sum_k e_k$. Then it's clear that:
\begin{equation*}
    e_k(z) : \sfW_{\frakq'} \longrightarrow \sfW_{\frakq'} \cap \sfW^{G,(k)}
\end{equation*}
is a projector to degree $k$ piece, which is integral over $(\Sh)_{\frakq'}$.

\emph{Step 5: conclusion}. Now the situation is: $\sfW_{\frakq'}$ is a free module over $(S_\hbar)_{\frakq'}$, with a filtration in Step 3, $\mathrm{Fil}^{\leq k}$. It has two splittings: $\mathrm{Fil}^{\leq k} = \oplus_{k' \leq k} (\sfW_{\frakq'} \cap \sfW^{G,(k')}) = \oplus_{k'\leq k} \sfW_{\frakq'}^{(k')}$, the first one is $G$-splitting, the second one is $L$-splitting localized at $\frakq'$. Both of them are integral over $(\Sh)_{\frakq'}$, thus, the projection along $\mathrm{Fil}^{<k}$ gives an isomorphism: 
\begin{equation*}
    \pi_k: \sfW_{\frakq'} \cap \sfW^{G,(k)} \xrightarrow{\sim} \sfW_{\frakq'}^{(k)}
\end{equation*}
In particular, for $\nu$ of degree $k$, we have $\pi_k (\sfW_{\frakq'} \cap \sfW_\nu) = \sfW_{\frakq'} \cap \sfW_\nu^L$. Finally apply $J = J^L\circ \pi_k$, we get:
\begin{equation*}
    (\calJ_\nu)_\frakq = J(\sfW_{\frakq'} \cap \sfW_\nu) = J^L\circ \pi_k (\sfW_{\frakq'} \cap \sfW_\nu) = J^L (\sfW_{\frakq'} \cap \sfW_\nu^L) = (\calJ_\nu^L)_\frakq
\end{equation*}
This finishes the proof of the Lemma.
\end{proof}

\subsection*{Integrality at walls of simple roots}
\begin{proposition}\label{prop:simple-walls}
For every simple root $\alpha$, every wall $\frakq$ of $\alpha$ and every $\nu$, $\calE_\nu$ is an automorphism of $\sfA_\nu\otimes(\Sh)_\frakq$.
\end{proposition}

\begin{proof}
By Lemma~\ref{lem:rank-one-all} for $L$, $\calJ^L_\nu=(\delta^L_\nu)^{-1}\calC^L_\nu$. The factors of $\delta_\nu/\delta^L_\nu$ are walls of roots $\gamma\neq\pm\alpha$, units in $(\Sh)_\frakq$. By Lemmas~\ref{lem:locC} and~\ref{lem:locJ}, $(\calJ_\nu)_\frakq=(\calJ^L_\nu)_\frakq=(\delta^L_\nu)^{-1}(\calC^L_\nu)_\frakq=\delta_\nu^{-1}(\calC_\nu)_\frakq$. Apply Lemma~\ref{lem:int}.
\end{proof}

\subsection{Weyl reflections}\label{sec:disc-weyl}
With the knowledge of \emph{simple roots} in hand, we can handle \emph{all roots}, using the method of Weyl reflection. Such ideas have appeared in lots of places, such as \cite{ginzburg2014differentialoperatorsguaffine} and \cite{bezrukavnikov2026equivariantsatakecategorykostantwhittaker}. Let's briefly recall the underlying philosophy of Weyl reflection: let $\gamma \in R^+$ be a positive root, which is not simple. Then we can pick a simple root $\alpha_i$, with $\langle \gamma, \check{\alpha}_i \rangle >0$. In that case, $\gamma':= s_i \gamma$ is a positive root of smaller height, where $s_i :=s_{\alpha_i}$ is the corresponding simple reflection. Let $y$ denotes the coordinate on $S_\hbar$, then the dot action $y\mapsto s_i\cdot y$, sends walls of $\gamma$ to walls of $\gamma'$. Now, assume we want to deduce, a property $(P)$ holds for all roots, and we already know that $(P)$ holds for simple roots. If we can "lift" simple reflections $s_i$'s to the reflections on the objects (e.g. modules, automorphisms, etc.) we care about, then we can assert $(P)$ holds for all roots.

In our situation, the property $(P)$ is the following:
\begin{equation*}
    \calE_\nu \in GL(\sfA_\nu \otimes (\Sh)_\frakq), \quad \text{for}\, \,\frakq \, \,\text{a wall of root} \,\gamma, \quad \forall\,\,\text{weights}\,\, \nu.
\end{equation*}
and the known input from Section~\ref{sec:disc-simple} is:
\begin{equation*}
    \calE_\nu \in GL(\sfA_\nu \otimes (\Sh)_\frakq), \quad \text{for}\, \,\frakq \, \,\text{a wall of simple root} \,\alpha, \quad \forall\,\,\text{weights}\,\, \nu.
\end{equation*}
Therefore, to conduct Weyl reflection, we need to find a "lift" of $s_i \in W$ to an $\Sh$-semilinear automorphism of $\sfA$:
\begin{equation*}
    s_{i,\calG}: \sfA_\nu \xrightarrow{\sim} \sfA_{s_i\nu}, \quad s_{i,\calG} \,\, \text{conjuagtes}\, \, \calE_\nu \mapsto \calE_{s_i\nu}, \quad \forall\,\,\text{weights}\,\, \nu.
\end{equation*}
Then by semi-linearity, $s_{i.\calG}$ moves the parameter in $\Sh$ by $s_i$, it carries localization at a wall of $\gamma$ to localization at a wall of $\gamma'$. Therefore, if $\calE_\nu$ is invertible at walls of $\gamma'$ (for all weights $\nu$), then $\calE_{s_i\nu}$ is invertible at walls of $\gamma$ (for all $s_i\nu$, thus all the weights). By induction on height of $\gamma$, we can propagates property $(P)$ to all roots. So the core of the question reduces to the construction of the desired $s_{i,\calG}$. We will define $w_{\calG}$ for all elements $w \in W$ below.

\subsection*{Definition of the reflections}
Recall from \eqref{eq:disc-diagram} the main object of studying is the following diagram over $\ars$:
\begin{equation}\label{eq:commu-diagram}
\begin{tikzcd}[column sep=5em,row sep=3em]
\sfW_\nu\arrow[r,"\BF","\sim"']
& \sfH_\nu
\\
\sfA_\nu\arrow[u,"\wr"',"\Phi"]\arrow[r,"\calG_\nu"',"\sim"]
& \sfA_\nu\arrow[u,"\GR"',"\wr"]
\end{tikzcd}
\end{equation}
In general, this diagram still does not commute, and what measures its non-commutativity is exactly
\begin{equation*}
    \calE = \GR^{-1}\circ \BF \circ \Phi \circ \calG_G^{-1}
\end{equation*}

Notice that the Weyl group $W$ does act on $\sfH$ and $\sfW$. More precisely, $W$ acts on $\sfW=\Sh\otimes_{\Zh(\frakg)}(V\otimes\calW_\hbar)^{\fraku,\psi}$ through $\Sh$, and similarly on $\sfH=\Sh\otimes_{\Zh(\frakg)}\H_{\GO}(\GrG,\IC_V)$ through $\Sh$. Moreover, $\BF$ is $W$-equivariant: it is the base change of \eqref{eq:BF} along the left action of $\Zh(\frakg)$. Let's denote the $w \in W$ action on $\sfW$ by $w_\sfW$, and on $\sfH$ by $w_\sfH$. If we want to translate those $W$-action to $\sfA$, we can start at the piece $\sfA_\nu$ sitting in the bottom-right corner of \eqref{eq:commu-diagram}, traveling to $\sfH_\nu$ via $\GR$, and traveling to $\sfW_\nu$ via $\Phi \circ \calG_\nu^{-1}$, then pull back the $W$-actions there. Therefore, for $w\in W$, we define the following three automorphisms in $GL(\sfA_{rs})$:
\[
w_\Phi:=\Phi^{-1}w_{\sfW}\Phi,\qquad w_{\calG}:=\calG_G\,w_\Phi\,\calG_G^{-1},\qquad w_{\GR}:=\GR^{-1}w_{\sfH}\GR .
\]
Each maps $\sfA_{\nu}$ to $\sfA_{w\nu}$, and is semilinear \footnote{semilinear over $\Sh$ means, $w(a\cdot x) = (w\cdot a) \cdot w(x)$ for $w \in W$, $a \in \Sh$, where the dot action of $w$ on $a$, viewed as a function over $y\in \fraka^*$, is given by $(w\cdot a)(y) = a(w^{-1}\cdot y)$. For instance, $w_\sfW$ is $\Sh$ semilinear because it only acts on the first factor: $w_\sfW (a\otimes u) = (w\cdot a)\otimes u$. If we think those modules are fibering over $\Sh$,  then in fiber terms, $w$ sends fiber at $y$ to fiber at $w\cdot y$.} over $\Sh$: that is, let $y$ denote the coordinate of $\Sh$, then $y\mapsto w\cdot y$ under $w_\Phi$, $w_\calG$\footnote{semilinearity of $w_\calG$ is exactly what moves the walls. If $\frakq \subset \Sh$ is a prime, then $w_\calG (\frakq \sfA_\nu) = (w\cdot \frakq) \sfA_{w\nu}$. So conjugating by $w_\calG$ turns a property at a wall of $\gamma'$ to a wall of $\gamma=w\gamma'$.}, $w_\GR$. Below, let $y$ be \emph{generic} in $\ars$, with $\hbar\neq0$; identities between maps defined over $\Crs$ are checked at such points.

\subsection*{$\calE$ intertwins two reflections}

\begin{lemma}\label{lem:W-basic}
(a) $w_{\GR}(\sfA_\nu)=\sfA_{w\nu}$. (b) $\calE_{w\nu}\circ w_{\calG}=w_{\GR}\circ\calE_\nu$ on $\sfA_{\nu,\rs}$.
\end{lemma}

\begin{proof}
For part (a), just notice that $\GR(\sfA_\nu)=\sfH\cap\sfH_\nu$ \cite[Theorem 1.4.1]{ginzburg2022differentialoperatorsgugelfandgraev} \footnote{It is proved in \cite{ginzburg2022differentialoperatorsgugelfandgraev} that the image of $\sfA_\nu$ inside $\sfH$ is cut out by the ideal $J_\nu$, $\sfH^{J_\nu}$ can be equivalently described as $\sfH \cap \sfH_\nu$. }, $w_{\sfH}(\sfH)=\sfH$, and $w_{\sfH}(\sfH_\nu)=\sfH_{w\nu}$. \\
For part (b), by our definition, $\calE=\GR^{-1} \circ \BF\circ\Phi\circ\calG_G^{-1}$. Therefore, 
\begin{equation*}
    \calE \circ w_{\calG} = \GR^{-1}\circ \BF \circ \Phi \circ \calG_{G}^{-1} \circ w_{\calG} = 
\end{equation*}
\begin{equation*}
    \GR^{-1}\circ \BF \circ w_{\sfW} \circ \Phi \circ \calG_{\calG}^{-1} =   \GR^{-1} \circ w_{\sfH} \circ \BF \circ \Phi \circ \calG_{\calG}^{-1}
\end{equation*}
\begin{equation*}
    = w_{GR} \circ \GR^{-1}\circ \BF \circ \Phi \circ \calG_{\calG}^{-1} = w_{GR}\circ \calE
\end{equation*}
Since Weyl reflection by $w \in W$ maps $\nu$ weight space to $w\nu$, the conclusion follows.
\end{proof}
We make the following definition:
\begin{definition}[Integrality of reflections]\label{def:w-int}
    We call a $\Sh$-semilinear automorphism $w$ (which stands for $w_\GR$, $w_\calG$) on $\sfA_{rs}$ integral, iff it restricts to an isomorphism between lattices: $w: \sfA_\nu \xrightarrow{\sim} \sfA_{w\nu}$; Likewise, we call $w$ is integral at $\frakq$, iff it restricts to an isomorphism $w: \sfA_\nu \otimes (\Sh)_\frakq \xrightarrow{\sim} \sfA_{w\nu} \otimes (\Sh)_{w\cdot\frakq}$.
\end{definition}

Thus, part (a) of Lemma~\ref{lem:W-basic} says $w_\GR$ is integral. Part (b) of Lemma~\ref{lem:W-basic} tells us that, integrality of $\calE$ propagates from the walls of $\gamma$ to those of $w\gamma$, once we know $w_{\calG}$ is integral.

\subsection*{Integrality of $w_\calG$}
Our goal is reduced to showing $w_\calG$ is integral. The integrality of $w_\calG$ follows easily from the integrality of all the simple reflections $s_{i,\calG}$. To show each $s_{i,\calG}$ is integral, our strategy is to compare it with $s_{i,GR}$, which is known to be integral. Interestingly, here we need to apply the idea of reducing to rank one once again. Denote $L_i \subset G$ the minimal Levi subgroup corresponding to the simple reflection $s_i$, since we have the full knowledge about the rank one case, we will show: $s_{i,\GR}^{L_i} = c\cdot s_{i,\calG}^{L_i}$, for some constant $c \in \bbC^\times$, where $s_i^{L_i} \in  W_{L_i}$ is the simple reflection of the rank one subgroup $L_i$. Then, we transport this relation back to $G$, via $\pr_{L_i}$, and obtain: $s_{i,GR} = c\cdot s_{i,\calG}$ for the same constant $c$, proving the integrality.

\subsection*{Rank one $L_i$}
\begin{proposition}\label{prop:W-step}
Let $\alpha_i$ be a simple root, $L_i$ its minimal Levi subgroup. We assume that $\calE^{L_i}=t_i$ for some $t_i\in T$ (Theorem~\ref{thm:DF} for $L_i$). Then on $\sfA_\nu^{L_i}$, we have: $s_{i,\GR}^{L_i}=(s_i\nu-\nu)(t_i)\cdot s_{i,\calG}^{L_i}$. 
\end{proposition}

\begin{proof}
First recall the meaning of $\calE= t$ for some element $t$ in the torus: on the $\nu$-weight piece, $\calE_\nu$ acts by multiplying the scalar $\nu(t) \in \bbC^\times$. Here we view $\nu$ as a character of the torus, $\nu : T \to \bbC^\times$. Therefore, for $\xi \in \sfA_\nu^{L_i}$:
\begin{equation*}
    s_{i,GR}^{L_i}(\xi) = t_i s_{i,\calG}^{L_i} t_i^{-1} (\xi) = t_i s_{i,\calG}^{L_i} (\nu(t_i)^{-1}\cdot\xi)
\end{equation*}
\begin{equation*}
    =\nu(t_i)^{-1} \cdot t_i s_{i,\calG}^{L_i} (\xi) = \nu(t_i)^{-1} \cdot(s_i\nu)(t_i) \cdot  s_{i,\calG}^{L_i}(\xi)
\end{equation*}
moreover $\frac{(s_i\nu)(t_i)}{\nu (t_i)}$ can be viewed as the value of the character $s_i \nu - \nu$ at $t_i$, therefore $s_{i,\GR}^{L_i}=(s_i\nu-\nu)(t_i)\cdot s_{i,\calG}^{L_i}$. Here $(s_i\nu-\nu)(t_i)$ is what we denoted by the constant $c \in \bbC^\times$ above.
\end{proof}

\subsection*{Transport back to $G$}
The next step is to prove the following two relations:
\begin{equation*}
    s_{i,\calG} = \pr_{L_i}^{-1}\circ s_{i,\calG}^{L_i}\circ \pr_{L_i}
\end{equation*}
\begin{equation*}
     s_{i,\GR} = \pr_{L_i}^{-1}\circ s_{i,\GR}^{L_i}\circ \pr_{L_i}
\end{equation*}
Indeed, we will prove the following two more general statements.

\begin{lemma}\label{lem:WD}
Let $L$ be a standard Levi subgroup and $w\in W_L$. Then $w_{\calG}=\mathrm{pr}_L^{-1}\circ w^L_{\calG}\circ\mathrm{pr}_L$.
\end{lemma}

\begin{lemma}\label{lem:WP}
Let $L$ be a standard Levi subgroup and $w\in W_L$. Then $w_{\GR}=\mathrm{pr}_L^{-1}\circ w^L_{\GR}\circ\mathrm{pr}_L$.
\end{lemma}

Assuming these, it is then straight forward to see:
\begin{equation*}
    s_{i,\GR} = \pr_{L_i}^{-1}\circ s_{i,\GR}^{L_i}\circ \pr_{L_i} = \pr_{L_i}^{-1}\circ c\cdot s_{i,\calG}^{L_i}\circ \pr_{L_i} = c \cdot s_{i,\calG}
\end{equation*}
because $c$ is a scalar. Lemma~\ref{lem:WP} follows from \cite[Lemma 4.6.5]{ginzburg2014differentialoperatorsguaffine}.
To prove the relation in Lemma~\ref{lem:WD}, it suffices to prove it at a generic point. \footnote{ Put $R:=\Crs[\hbar^{-1}]$. All identities we want to prove are between linear or semilinear homomorphisms between free $R$-modules of finite rank, which are already defined over $R$. After choosing a $R$-basis, such a $R$-linear homomorphism can be written as a matrix with entries in $R$ (if it is semilinear, with a  twist by $y\mapsto w\cdot y$). Because $R$ is the ring of regular functions on an irreducible variety, two such matrices coincide iff their specializations coincide at all generic $y$.}

\paragraph{Notation at a generic point.} Call a point \emph{generic} if it is a point of $\ars(\bbC)$ with $\hbar\neq0$.  Below, we will work with one such point at a time. A point of $\ars(\bbC)$ with $\hbar$-coordinate $a\neq0$ is denoted simply by $y$. At such a point we write:
\begin{itemize}
\item $M(y):=\bM_\hbar\otimes_{\Sh}\bbC_y$ and $M_L(y):=\bM_{L,\hbar}\otimes_{\Sh}\bbC_y$, and $\widehat M(y)$, $\widehat M_L(y)$ for their completions; as $U(\frakg)$- resp.\ $U(\frakl)$-modules these are the Verma modules of highest weight $y/a$;
\item $\chi_y$ for the character by which $\Zh(\frakg)$ acts on $M(y)$, and $\chi^L_y$ for the character by which $\Zh(\frakl)$ acts on $M_L(y)$;
\item $\calW_\hbar(\chi):=\calW_\hbar\otimes_{\Zh(\frakg)}\bbC_\chi$ and $M(\psiL,\chi^L):=\bM_\hbar(\psiL)\otimes_{\Zh(\frakl)}\bbC_\chi$ for characters $\chi$ of $\Zh(\frakg)$, and characters $\chi^L$ of $\Zh(\frakl)$;
\item $\wv(y)\in\widehat M(y)$ and $\wv_L(y)\in\widehat M_L(y)\subset\widehat M(y)$ for the specializations of the normalized $\psi$- and $\psiL$-Whittaker vectors $\wv$, $\wv_L$.
\end{itemize}
The characters $\chi_y,\chi^L_y$ can be thought as the image of $y$ modulo the dot action of $W$, $W_L$ on $\ars$. By Harish-Chandra, $\chi_{w\cdot y}=\chi_y$ for $w\in W$ and $\chi^L_{w\cdot y}=\chi^L_y$ for $w\in W_L$; so the modules $\calW_\hbar(\chi_y)$,  $M(\psiL, \chi^L_y)$ only see $W$, \ $W_L$-invariant things. The Verma modules $M(y)$, $M_L(y)$ and the vectors $\wv(y)$, $\wv_L(y)$ depend on $y$ itself: for $w\in W_L$, $\wv_L(y)$ and $\wv_L(w\cdot y)$ live in different Verma modules with the same central character. The $W_L$-actions on the following modules:
\[
\sfW_\rs,\qquad \sfW^{G,L}_\rs:=(V\otimes\bM_\hbar(\psiL))^{\fraku,\psiL}\otimes_{\Zh(\frakl)}\Crs,\qquad \sfW^L_\rs=(V|_L\otimes\calW_{L,\hbar})^{\fraku_L,\psiL}\otimes_{\Zh(\frakl)}\Crs
\]
are all through the coefficients $\Crs$. Their fibers at $y$ are \footnote{at generic point $y$, fibers of invariants are invariants of fibers. For Verma modules, we have $\sfA\otimes_{\Sh}\bbC_y\xrightarrow{\sim}(V\otimes M(y))^{\fraku}$, and similarly for $\sfW$, $\sfW^{G,L}$, $\sfW^L$. For $\sfA$, the map is injective because its composite with $\Top$ at $y$ is the specialization of $\Top$, an isomorphism over $\Crs$ , see Remark~\ref{rem:top-generic}, and both sides have dimension $\dim V$ because $V\otimes M(y/a)$ is a direct sum of simple Verma modules. For the Whittaker modules, this follows from exactness of $\mathrm{Wh}$ and Kostant's freeness (Lemma~\ref{lem:Kostant-for-s}), together with Proposition~\ref{prop:prV} for $\sfW^{G,L}$. This is one of the merit of working with generic point.}:
\[
\sfW(y)=(V\otimes\calW_\hbar(\chi_y))^{\fraku,\psi},\qquad \sfW^{G,L}(y)=(V\otimes M(\psiL,\chi^L_y))^{\fraku,\psiL},\qquad \sfW^L(y)=(V|_L\otimes\calW_{L,\hbar}(\chi^L_y))^{\fraku_L,\psiL},
\]
which depend on $y$ only through $\chi_y$ resp.\ $\chi^L_y$. So for $w\in W_L$ the fibers at $y$ and at $w\cdot y$ are the same vector space, and the action of $w$ is this identification. In particular, a map $F$ between two of these modules, defined fiberwise by maps $F_y$, commutes with $W_L$ if and only if $F_y=F_{w\cdot y}$ for all $y$ and all $w\in W_L$.
\\ \hspace{\fill} \\
The proof of Lemma~\ref{lem:WD} has two parts: the first part is to prove the same identity for $w_\Phi$, that is, $w_\Phi = \pr_L^{-1}\circ w_\Phi^{L} \circ \pr_L$. This is the most complicated part; The second part is gauging by $\calG$. Recall that $w_\calG = \calG \circ w_\Phi \circ \calG^{-1}$, and $w_\calG^L = \calG \circ w_\Phi^L \circ \calG^{-1}$. The effect of interchanging $\calG_G$, $\calG_L$ with $\pr_L$ follows from the expression of $\delta_\nu$, $\delta_\nu^L$. The difficulty in the first part is that, $w_\Phi$ is defined using the $\psi$-Whittaker vector $\wv$ of $G$ inside $\widehat\bM_\hbar$, while $w_\Phi^L$ is defined using $\psiL$-Whittaker vector $\wv_L$ of $L$ inside $\widehat\bM_{L,\hbar}$. The idea to compare them is to factor $\Phi$ through $\wv_L$, and such that the remaining factor commutes with $W_L$. More concretely speaking, the $\psi$-Whittaker vector of the Verma module can be built in two stages: first the $\psiL$-Whittaker vector $\wv_L(y)$, which generates a copy of the Verma--Whittaker module $M(\psiL,\chi^L_y)$ inside $\widehat M(y)$; then the $\psi$-Whittaker vector of that Verma--Whittaker module.

\begin{lemma}[Transitivity of Whittaker vectors]\label{lem:wv-trans}
Let $y$ be generic and $\chi:=\chi^L_y$. Write $\widehat{M(\psiL,\chi)}$ for the completion.Then the following map:
\begin{equation*}
    \iota_y\colon M(\psiL, \chi)\to\widehat M(y),\quad 1\mapsto\wv_L(y)
\end{equation*}
is a well-defined injective homomorphism between $\Uh(\frakg)$-modules \footnote{The map $\iota_y$ is only defined pointwise. As a consequence, the map $\vartheta_y$ below is only pointwise defined, too. This is the reason why we work pointwise at generic points.}. Moreover, $\widehat{M(\psiL,\chi)}$ contains a unique $\psi$-Whittaker vector $\wv(\chi)_\psi$ whose top $\frakz$-weight component is $1$. Finally $\iota_y$ satisfies: $\iota_y(\wv(\chi)_\psi)=\wv(y)$. Equivalently, the triangle
\begin{equation}\label{eq:wv-trans}
\begin{tikzcd}[column sep=4em]
\calW_\hbar(\chi_y)\arrow[r,"1\,\mapsto\,\wv(\chi)_\psi"]\arrow[rr,bend right=18,"1\,\mapsto\,\wv(y)"'] & \widehat{M(\psiL,\chi)}\arrow[r,"\iota_y"] & \widehat M(y)
\end{tikzcd}
\end{equation}
commutes. The bottom arrow is the specialization of $\phi_\hbar$ at $y$; the left arrow depends only on $\chi$.
\end{lemma}

\begin{proof}
By PBW, as $\frakl$-modules, we have:
\begin{equation*}
    \widehat M(y)=\prod_\gamma\Uh(\bar\fraku_P)_{-\gamma}\otimes\widehat M_L(y),\qquad
M(\psiL, \chi)=\bigoplus_\gamma\Uh(\bar\fraku_P)_{-\gamma}\otimes\calW_{L,\hbar}(\chi)
\end{equation*}
\emph{Step 1: the map $\iota_y$}. $\wv_L(y)$ lies in the top $\frakz$-weight piece $1\otimes\widehat M_L(y)$, so it is killed by $\fraku_P$, which raises the $\frakz$-weight; it is a $\psiL$-Whittaker vector for $\fraku_L$; and $\Zh(\frakl)$ acts on it by $\chi$. Since $M(\psiL, \chi)=\Uh(\frakg)\otimes_{\Uh(\frakp_L)}\calW_{L,\hbar}(\chi)$, Frobenius reciprocity gives $\iota_y$; on the top piece it is the isomorphism $\calW_{L,\hbar}(\chi)\xrightarrow{\sim}\Uh(\frakl)\wv_L(y)$, and on the $\gamma$-th piece it is
\[
\Uh(\bar\fraku_P)_{-\gamma}\otimes\calW_{L,\hbar}(\chi)\xrightarrow{\ \sim\ }\Uh(\bar\fraku_P)_{-\gamma}\otimes\Uh(\frakl)\wv_L(y)\ \subset\ \Uh(\bar\fraku_P)_{-\gamma}\otimes\widehat M_L(y).
\]
$\iota_y$ is non-zero and $M(\psiL, \chi)$ is simple, therefore $\iota_y$ is injective.

\emph{Step 2: decompose $\wv(y)$}. Write $\wv(y)=\sum_\gamma\wv_\gamma$ with $\wv_\gamma\in\Uh(\bar\fraku_P)_{-\gamma}\otimes\widehat M_L(y)$. Since $\fraku_L$ commutes with $\frakz$ and $\psi|_{\fraku_L}=\psiL$, each $\wv_\gamma$ is a $\psiL$-Whittaker vector for $\fraku_L$, i.e.\ $\wv_\gamma\in Y_\gamma:=\bigl(\Uh(\bar\fraku_P)_{-\gamma}\otimes\widehat M_L(y)\bigr)^{\fraku_L,\psiL}$.

\emph{Step 3: dimension count}. Put $X_\gamma:=\bigl(\Uh(\bar\fraku_P)_{-\gamma}\otimes\Uh(\frakl)\wv_L(y)\bigr)^{\fraku_L,\psiL}\subseteq Y_\gamma$. By Step~1, $\Uh(\bar\fraku_P)_{-\gamma}\otimes\Uh(\frakl)\wv_L(y)\cong\Uh(\bar\fraku_P)_{-\gamma}\otimes\calW_{L,\hbar}(\chi)$, so $\dim X_\gamma=\dim\Uh(\bar\fraku_P)_{-\gamma}$ by (K2) of Lemma~\ref{lem:Kostant-for-s}. On the other hand, Whittaker vectors of the completion are $\psiL$-finite, therefore:
\[
Y_\gamma\subset\widehat\Gamma_{\psiL}\bigl(\Uh(\bar\fraku_P)_{-\gamma}\otimes M_L(y)\bigr)=\Uh(\bar\fraku_P)_{-\gamma}\otimes\widehat\Gamma_{\psiL}(M_L(y))\cong\Uh(\bar\fraku_P)_{-\gamma}\otimes\calW_{L,\hbar}(\chi)
\]
hence $\dim Y_\gamma=\dim\Uh(\bar\fraku_P)_{-\gamma}$ as well. An inclusion of the same finite-dimensional vector spaces must be an equality: $X_\gamma=Y_\gamma$. Therefore
\[
\wv_\gamma\in\Uh(\bar\fraku_P)_{-\gamma}\otimes\Uh(\frakl)\wv_L(y)=\iota_y\bigl(\Uh(\bar\fraku_P)_{-\gamma}\otimes\calW_{L,\hbar}(\chi)\bigr).
\]

\emph{Step 4: conclusion.} By Step~3 and the injectivity of $\iota_y$ there are unique $\zeta_\gamma$ with $\iota_y(\zeta_\gamma)=\wv_\gamma$; put $\zeta:=\sum_\gamma\zeta_\gamma\in\widehat{M(\psiL, \chi)}$, so $\iota_y(\zeta)=\wv(y)$. For $x\in\fraku$, $\iota_y\bigl((x-\psi(x))\zeta\bigr)=(x-\psi(x))\wv(y)=0$, so $\zeta$ is a $\psi$-Whittaker vector; its top component is $1$, because $\wv_0=\wv_L(y)=\iota_y(1)$. This proves existence. For uniqueness, let $\zeta'$ be another such vector. Then $\iota_y(\zeta')$ is a $\psi$-Whittaker vector of $\widehat M(y)$ whose component of weight $y$ is $1$ (the top component of $\zeta'$ is $1$, and $\iota_y(1)=\wv_L(y)$ has leading term $\one$). Since $M(y)$ is simple, the $\psi$-Whittaker vectors inside $\widehat{M(y)}$ is one dimensional, thus $\iota_y(\zeta')=\wv(y)=\iota_y(\zeta)$. $\iota_y$ being injective, so $\zeta'=\zeta$.
\end{proof}

\begin{proof}[Proof of the Lemma~\ref{lem:WD}]
Recall:
\[
w_{\calG}=\calG_G\,w_\Phi\,\calG_G^{-1},\quad w_\Phi=\Phi^{-1}w_{\sfW}\Phi,\qquad\qquad
w^L_{\calG}=\calG_L\,w^L_\Phi\,\calG_L^{-1},\quad w^L_\Phi=(\Phi^L)^{-1}w_{\sfW}^L\Phi^L .
\]
We factor $\Phi$ through $\wv_L$ (Step~1), check that the remaining factor commutes with $W_L$ (Step~2), compare (Step~3), and finally deal with the gauges (Step~4). Let
\[
\Phi^{G,L}\colon\sfA_\rs\longrightarrow\sfW^{G,L}_\rs,\qquad \xi\longmapsto(V\otimes\phi_\hbar)^{-1}\bigl(f_\xi(\wv_L)\bigr),
\]
be the map \eqref{eq:Phi-def} of Case~2 (over $\Crs[\hbar^{-1}]$), and $\mathrm{pr}_{\mathcal W}\colon\sfW^{G,L}_\rs\to\sfW^L_\rs$ the map \eqref{eq:prV}, an isomorphism over $\Crs$ by Proposition~\ref{prop:prV}.

\emph{Step 1: $\Phi$ factors through $\Phi^{G,L}$.} Fix a generic $y$ and $\xi\in\sfA_\nu$, and put $\chi:=\chi^L_y$, $\chi':=\chi^L_{y+\hbar\nu}$. The element $\zeta:=\Phi^{G,L}(\xi)(y)\in\sfW^{G,L}(y)$ has the property that the left action of $\Zh(\frakl)$ on it is through $\chi'$ (\S\ref{sec:intertwiners}), so it corresponds to an intertwiner $g_\zeta\colon M(\psiL,\chi')\to V\otimes M(\psiL, \chi)$ with $g_\zeta(1)=\zeta$. The square
\[
\begin{tikzcd}[column sep=5em,row sep=3em]
M(\psiL,\chi')\arrow[r,"\iota_{y+\hbar\nu}"]\arrow[d,"g_\zeta"'] & \widehat M(y+\hbar\nu)\arrow[d,"f_\xi"]\\
V\otimes\widehat{M(\psiL,\chi)}\arrow[r,"V\otimes\iota_y"'] & V\otimes\widehat M(y)
\end{tikzcd}
\]
commutes: both composites are $\Uh(\frakg)$-linear, and both send $1$ to $f_\xi(\wv_L(y+\hbar\nu))=(V\otimes\iota_y)(\zeta)$, by the definition of $\Phi^{G,L}$ (the specialization of $\phi_\hbar$ at $y$ is $\iota_y$). Define
\[
\vartheta_y\colon\sfW^{G,L}(y)\to\sfW(y),\qquad \zeta\longmapsto(V\otimes\kappa_\chi)^{-1}\bigl(g_\zeta(\wv(\chi')_\psi)\bigr),
\]
where $\kappa_\chi\colon\calW_\hbar(\chi_y)\to\widehat{M(\psiL, \chi)}$, $1\mapsto\wv(\chi)_\psi$, is the left arrow of \eqref{eq:wv-trans}. On $\psi$-Whittaker vectors, $V\otimes\kappa_\chi$ is an isomorphism: by \eqref{eq:wv-trans}, $(V\otimes\iota_y)\circ(V\otimes\kappa_\chi)$ is the specialization of $V\otimes\phi_\hbar$, which is an isomorphism on $\psi$-Whittaker vectors, and $V\otimes\iota_y$ is injective. Now use Lemma~\ref{lem:wv-trans}, gives:
\[
(V\otimes\iota_y)\bigl(g_\zeta(\wv(\chi')_\psi)\bigr)=f_\xi\bigl(\iota_{y+\hbar\nu}(\wv(\chi')_\psi)\bigr)=f_\xi\bigl(\wv(y+\hbar\nu)\bigr)=(V\otimes\iota_y)(V\otimes\kappa_\chi)\bigl(\Phi(\xi)(y)\bigr).
\]
Cancelling the injective map $V\otimes\iota_y$ gives $\Phi(\xi)(y)=\vartheta_y\bigl(\Phi^{G,L}(\xi)(y)\bigr)$. Hence, fiberwise at every generic point,
\begin{equation}\label{eq:Phi-factor}
\begin{tikzcd}[column sep=4em]
\sfA_\rs\arrow[r,"\Phi^{G,L}"]\arrow[rr,bend right=18,"\Phi"'] & \sfW^{G,L}_\rs\arrow[r,"\vartheta"] & \sfW_\rs .
\end{tikzcd}
\end{equation}
In words: evaluating an intertwiner on $\wv$ is the same as first evaluating it on $\wv_L$, and then evaluating the resulting Verma--Whittaker intertwiner on $\wv(\chi')_\psi$.

\emph{Step 2: homomorphisms commute with $W_L$.} By the notation paragraph, a fiber-wise map commutes with $W_L$ iff its fibers at $y$ and $w\cdot y$ agree.
\begin{itemize}
\item $\vartheta$ commutes with $W_L$. Indeed $\vartheta_y$ is built from $M(\psiL,\chi')$, $M(\psiL, \chi)$, $\calW_\hbar(\chi_y)$ and the vectors $\wv(\chi')_\psi$, $\wv(\chi)_\psi$, which by Lemma~\ref{lem:wv-trans}(b) depend only on the characters $\chi=\chi^L_y$, $\chi'=\chi^L_{y+\hbar\nu}$ and $\chi_y$. For $w\in W_L$ these are unchanged under $y\mapsto w\cdot y$, $\nu\mapsto w\nu$, because $w\cdot(y+\hbar\nu)=w\cdot y+\hbar w\nu$. So $\vartheta_{w\cdot y}=\vartheta_y$:
\[
\begin{tikzcd}[column sep=4em,row sep=2.5em]
\sfW^{G,L}(y)\arrow[r,"\vartheta_y"]\arrow[d,equal,"w"'] & \sfW(y)\arrow[d,equal,"w"]\\
\sfW^{G,L}(w\cdot y)\arrow[r,"\vartheta_{w\cdot y}"] & \sfW(w\cdot y)
\end{tikzcd}
\]
\item $\mathrm{pr}_{\mathcal W}$ commutes with $W_L$: it is the base change along $\Zh(\frakl)\to\Crs$ of a $\Zh(\frakl)$-linear map, so it does not touch the coefficients.
\item $\Phi^{G,L}$ and $\Phi^L$ do \emph{not} commute with $W_L$: they use $\wv_L(y)$, which depends on $y$ itself. Instead, they are related, by Lemma~\ref{lem:alg-square}: $\mathrm{pr}_{\mathcal W}\circ\Phi^{G,L}=\Phi^L\circ\mathrm{pr}_L$.
\end{itemize}
Combine, we have the following diagram over $\Crs$, in which the square commutes, and the upper triangle is exactly \eqref{eq:Phi-factor}:
\begin{equation}\label{eq:WD-diagram}
\begin{tikzcd}[column sep=4.5em,row sep=3em]
\sfA_\rs\arrow[r,"\Phi^{G,L}"]\arrow[d,"\mathrm{pr}_L"']\arrow[rr,bend left=22,"\Phi"] & \sfW^{G,L}_\rs\arrow[r,"\vartheta"]\arrow[d,"\mathrm{pr}_{\mathcal W}","\wr"'] & \sfW_\rs\\
\sfA^L_\rs\arrow[r,"\Phi^L"'] & \sfW^L_\rs &
\end{tikzcd}
\end{equation}
Here $\vartheta$ and $\mathrm{pr}_{\mathcal W}$ commute with $W_L$, and $\mathrm{pr}_L$, $\Phi^{G,L}$, $\Phi^L$ do not.

\emph{Step 3: the comparison of $w_\Phi$ and $w^L_\Phi$.} Reading \eqref{eq:WD-diagram},
\[
\begin{aligned}
w_\Phi&=\Phi^{-1}w_{\sfW}\Phi=(\Phi^{G,L})^{-1}\bigl(\vartheta^{-1}w_{\sfW}\vartheta\bigr)\Phi^{G,L} &&\text{by \eqref{eq:Phi-factor}}\\
&=(\Phi^{G,L})^{-1}\,w_{\sfW^{G,L}}\,\Phi^{G,L} &&\text{$\vartheta$ commutes with $W_L$}\\
&=(\Phi^{G,L})^{-1}\bigl(\mathrm{pr}_{\mathcal W}^{-1}w_{\sfW}^L\mathrm{pr}_{\mathcal W}\bigr)\Phi^{G,L} &&\text{$\mathrm{pr}_{\mathcal W}$ commutes with $W_L$}\\
&=\mathrm{pr}_L^{-1}\,(\Phi^L)^{-1}w_{\sfW}^L\Phi^L\,\mathrm{pr}_L &&\text{Lemma~\ref{lem:alg-square}}\\
&=\mathrm{pr}_L^{-1}\,w^L_\Phi\,\mathrm{pr}_L .
\end{aligned}
\]
Each line is an identity at every generic point, hence an identity over $\Crs$; only the invertibility of $\mathrm{pr}_L$ over $\Crs$ , cf. Remark~\ref{rem:top-generic}, is used here.

\emph{Step 4: the gauges.} Write $\calG_G=\calG_L\,\calG_{G/L}$, where $\calG_{G/L}$ multiplies $\sfA_\nu$ by
\[
\frac{\delta_\nu}{\delta^L_\nu}=\prod_{\gamma\in R^+\smallsetminus R^+_L}P_{\langle\nu,\check\gamma\rangle}\bigl(\langle y+\hbar\rho,\check\gamma\rangle\bigr)^{-1}.
\]
For $w\in W_L$ we have $\langle w\cdot y+\hbar\rho,\check\gamma\rangle=\langle y+\hbar\rho,w^{-1}\check\gamma\rangle$ and $\langle w\nu,\check\gamma\rangle=\langle\nu,w^{-1}\check\gamma\rangle$, and $\gamma\mapsto w^{-1}\gamma$ permutes $R^+\smallsetminus R^+_L$. Hence $(\delta_{w\nu}/\delta^L_{w\nu})(w\cdot y)=(\delta_\nu/\delta^L_\nu)(y)$, which says exactly that the semilinear map $w_\Phi\colon\sfA_\nu\to\sfA_{w\nu}$ commutes with $\calG_{G/L}$:
\[
\begin{tikzcd}[column sep=4em,row sep=2.5em]
\sfA_{\nu,\rs}\arrow[r,"w_\Phi"]\arrow[d,"\cdot\,\delta_\nu/\delta^L_\nu"'] & \sfA_{w\nu,\rs}\arrow[d,"\cdot\,\delta_{w\nu}/\delta^L_{w\nu}"]\\
\sfA_{\nu,\rs}\arrow[r,"w_\Phi"'] & \sfA_{w\nu,\rs}
\end{tikzcd}
\]
Moreover $\calG_L$ commutes with $\mathrm{pr}_L$, by Lemma~\ref{lem:gauge}(b). Therefore, using Step 3, we get:
\[
w_\calG=\calG_L\calG_{G/L}\,w_\Phi\,\calG_{G/L}^{-1}\calG_L^{-1}=\calG_L\,w_\Phi\,\calG_L^{-1}
=\calG_L\,\mathrm{pr}_L^{-1}w^L_\Phi\,\mathrm{pr}_L\,\calG_L^{-1}=\mathrm{pr}_L^{-1}\bigl(\calG_L\,w^L_\Phi\,\calG_L^{-1}\bigr)\mathrm{pr}_L=\mathrm{pr}_L^{-1}\,w^L_\calG\,\mathrm{pr}_L.\qedhere
\]
\end{proof}

\subsection{Multiplicativity and rigidity}\label{sec:disc-rigid}
For $V_1,V_2\in\Rep(G)$, $\xi\in\sfA_\nu(V_1)$ and $\xi'\in\sfA_\mu(V_2)$ let $\xi'\circ\xi\in\sfA_{\nu+\mu}(V_2\otimes V_1)$ be their composition (Definition~\ref{def:composition}, Corollary~\ref{cor:comp-verma}(b)); likewise on the $\sfW$'s, and convolution on the $\sfH$'s. $\GR$ is multiplicative by \eqref{eq:GR-reg} and the lax monoidal structure of $i_T^!$ \cite[\S2.2]{macerato2023leviequivariantrestrictionsphericalperverse}, and $\BF$ by Convention~\ref{conv:BF}(M).

\begin{lemma}\label{lem:E-mult}
$\calE$ is $\Sh$-linear, natural in $V\in\Rep(G)$, and multiplicative: $\calE(\xi'\circ\xi)=\calE(\xi')\circ\calE(\xi)$.
\end{lemma}

\begin{proof}
By Lemma~\ref{lem:composition}(a), $f_{\xi'\circ\xi}=(\id\otimes f_\xi)\circ f_{\xi'}$ for $\xi\in\sfA(V_1)$, $\xi'\in\sfA(V_2)$, and likewise $g_{\eta'\circ\eta}=(\id\otimes g_\eta)\circ g_{\eta'}$ for intertwiners on the Whittaker side \footnote{for $\eta \in (V\otimes\calW_\hbar)^{\fraku,\psi}$, denote by $g_\eta : \calW_\hbar \to V\otimes \calW_\hbar$ the corresponding intertwiner}, where $g_\eta(\one)=\eta$. Identify $\calW_\hbar$ with $\Uh(\frakg)\wv\subset\Mhat$ via $\phi_\hbar$, so that $\Phi(\xi)=f_\xi(\wv)$. Since $f_\xi(\wv)$ is a $\psi$-Whittaker vector, it lies in $V_1\otimes\Uh(\frakg)\wv$ (Proposition~\ref{prop:phi-iso}), so $f_\xi$ maps $\Uh(\frakg)\wv$ into $V_1\otimes\Uh(\frakg)\wv$, and there it equals $g_{\Phi(\xi)}$: both are $\Uh(\frakg)$-linear and send $\wv$ to $\Phi(\xi)$. Hence
\[
\Phi(\xi'\circ\xi)=f_{\xi'\circ\xi}(\wv)=(\id\otimes f_\xi)\bigl(f_{\xi'}(\wv)\bigr)=(\id\otimes g_{\Phi(\xi)})\bigl(\Phi(\xi')\bigr)=\Phi(\xi')\circ\Phi(\xi).
\]
$\calG_G$ is multiplicative by Lemma~\ref{lem:gauge}(c), and $\GR$, $\BF$ are multiplicative as explained above.
\end{proof}

\begin{theorem}[Rigidity]\label{thm:rigid}
Suppose that $\calE_\nu(V)\in GL(\sfA_\nu(V))$ for all $V$ and $\nu$. Then $\calE=t$ for a unique $t\in T(\bbC)$ \footnote{the concrete meaning of $\calE = t$ is as before, $\calE_\nu$ acts on $\sfA_\nu$ by multiplying scalar $\nu(t)$}.
\end{theorem}

We first explain the idea of the proof. We transport $\calE$ to the vector spaces $V_\nu\otimes\bbK$ by $\Top$, obtaining automorphisms $e(V)$. If $\Top$ were multiplicative, multiplicativity of $\calE$ would say $e(V_2\otimes V_1)=e(V_2)\otimes e(V_1)$, and Tannakian duality would finish the proof. Unfortunately, it is not: by \eqref{eq:top-product} below, $\Top(\xi'\circ\xi)$ differs from $\Top(\xi')\otimes\Top(\xi)$ in two ways: (1) by \emph{correction terms} of lower weight; (2) by a \emph{shift} $y\mapsto y+\hbar\nu$ of the parameter in the first factor. Step~1 uses a filtration to ignore the corrections: it shows that $e(V_2\otimes V_1)$ is block triangular, with diagonal blocks given by shifted tensor products. Step~2 compares the two orderings: $V_1\otimes V_2$, and its flip $V_2\otimes V_1$. They give opposite triangularities, so $e$ is block diagonal. Moreover, they put the shift on different factors, and that gives a functional equation for $e$. In Step~3 this equation, together with integrality, shows that $e$ actually does not depend on $y$. Then, the shifts disappear, and $e$ becomes a desired tensor automorphism. Step~4 applies Tannakian duality to conclude.

\begin{proof}[Proof of Theorem~\ref{thm:rigid}] 
Let $\bbK:=\operatorname{Frac}\Sh=\bbC(\hbar)(y)$. By Remark~\ref{rem:top-generic}, $\Top\colon\sfA_\nu(V)\otimes\bbK\xrightarrow{\sim}V_\nu\otimes\bbK$; put
\[
e_\nu(V):=\Top\circ\,\calE_\nu(V)\circ\Top^{-1}\in GL(V_\nu\otimes\bbK),\qquad e(V):=\textstyle\bigoplus_\nu e_\nu(V)\in GL(V\otimes\bbK).
\]
For a matrix $m$ over $\bbK$ and $\delta\in\Lambda$, $m(y+\hbar\delta)$ denotes the matrix obtained by substituting $y+\hbar\delta$ for $y$. Since $\Top$ and $\calE$ are natural in $V$, by Lemma~\ref{lem:E-mult}, so is $e$; and $e(\mathrm{triv})=1$, because $\calE(1)=\calE(1\circ1)=\calE(1)^2$ and $\calE(1)$ is invertible.

\emph{Step 1: products and filtration}. Let $\xi\in\sfA_\nu(V_1)$ and write $\xi'=\sum_K\xi'_K\otimes y_K\one\in\sfA(V_2)$ with $\xi'_K\in V_2\otimes\Sh$, where a term $v\otimes s$ of $\xi'_K$ stands for $v\otimes y_K\one\cdot s$. Then
\[
\xi'\circ\xi=(\id\otimes f_\xi)(\xi')=\sum_K\sum v\otimes y_K\cdot f_\xi(\one\cdot s).
\]
Now $f_\xi(\one\cdot s)=\xi\cdot s(y+\hbar\nu)$ by Corollary~\ref{cor:comp-verma}(a): as explained in Corollary~\ref{cor:comp-verma}(b), the coefficients of $\xi'$ are evaluated at the parameter $y+\hbar\nu$ of the middle term. Taking the $\bM_\hbar$-degree-$0$ component, only the part of the diagonal action of $y_K$ acting on $V_1$ survives, so
\begin{equation}\label{eq:top-product}
\Top(\xi'\circ\xi)=\underbrace{\xi'_\emptyset(y+\hbar\nu)\otimes\Top(\xi)}_{\text{leading term}}+\underbrace{\sum_{K\neq\emptyset}\xi'_K(y+\hbar\nu)\otimes y_K\Top(\xi)}_{\text{corrections}},
\end{equation}
with $y_K$ acting on $V_1$. Note $\xi'_\emptyset=\Top(\xi')$. In a correction term, $y_K\in\Uh(\ub)$ strictly lowers the $V_1$-weight. This suggests us to filter by the weight of the second factor: i.e. on $(V_2\otimes V_1)_\lambda\otimes\bbK$, denote:
\[
\mathrm{Fil}_{\le\gamma}:=\bigoplus_{\gamma'\le\gamma}(V_2)_{\lambda-\gamma'}\otimes(V_1)_{\gamma'}\otimes\bbK,\qquad\mathrm{gr}_\gamma=(V_2)_{\lambda-\gamma}\otimes(V_1)_\gamma\otimes\bbK .
\]
By \eqref{eq:top-product}, $\Top(\xi'\circ\xi)\in\mathrm{Fil}_{\le\nu}$, and its class in $\mathrm{gr}_\nu$ is the main term $\Top(\xi')(y+\hbar\nu)\otimes\Top(\xi)$. Since $\Top$ is onto over $\bbK$, and the shift $y\mapsto y+\hbar\nu$ is an automorphism of $\bbK$, these classes span $\mathrm{gr}_\nu$. Induction on $\nu$, we see the elements $\Top(\xi'\circ\xi)$ with $\xi\in\sfA_\nu(V_1)$, $\nu\le\gamma$, span $\mathrm{Fil}_{\le\gamma}$. By multiplicativity in Lemma~\ref{lem:E-mult}, $e(V_2\otimes V_1)$ sends $\Top(\xi'\circ\xi)$ to $\Top(\calE\xi'\circ\calE\xi)$, which lies in $\mathrm{Fil}_{\le\nu}$, because $\calE\xi\in\sfA_\nu(V_1)$. It's image in $\mathrm{gr}_\nu$ is:
\[
\Top(\calE\xi')(y+\hbar\nu)\otimes\Top(\calE\xi)=e(V_2)(y+\hbar\nu)\,\Top(\xi')(y+\hbar\nu)\otimes e(V_1)\,\Top(\xi)
\]
Therefore, we see $e(V_2\otimes V_1)$ is block triangular, due to the fact that it preserves each $\mathrm{Fil}_{\le\gamma}$. On $\mathrm{gr}_\gamma$ it acts by $e_{\lambda-\gamma}(V_2)(y+\hbar\gamma)\otimes e_\gamma(V_1)(y)$. Notice that, the tensor product is twisted: the outer factor is shifted by the weight of the inner factor.

\emph{Step 2: flip}. Both the correction of lower weight, and the twist, will disappear, when we compare the two orderings of the tensor factors. The flip $\mathrm{sw}\colon V_2\otimes V_1\to V_1\otimes V_2$ is $G$-linear, by naturality $e(V_1\otimes V_2)=\mathrm{sw}\circ e(V_2\otimes V_1)\circ\mathrm{sw}^{-1}$. Fix a total weight $\lambda$ and consider the pieces $(V_1)_\gamma\otimes(V_2)_\delta$ of $(V_1\otimes V_2)_\lambda$, $\gamma+\delta=\lambda$. Then Step 1 gives two descriptions of $e(V_1\otimes V_2)$:

(1) For $V_1\otimes V_2$: the filtration is increasing in $\delta$. Therefore, the matrix of $e(V_1\otimes V_2)$ is triangular for the weight $\delta$ of $V_2$, thus block upper triangular. Its diagonal block equals $e_\gamma(V_1)(y+\hbar\delta)\otimes e_\delta(V_2)(y)$.

(2) For $V_2\otimes V_1$: the filtration is increasing in $\gamma$. Since $\lambda$ is fixed, it is decreasing in $\delta =\lambda -\gamma$. Therefore, the matrix of $\mathrm{sw}\circ e(V_2\otimes V_1)\circ\mathrm{sw}^{-1}$  is triangular for the  weight $\gamma$ of $V_1$, thus block lower triangular. Its diagnal block equals $e_\gamma(V_1)(y)\otimes e_\delta(V_2)(y+\hbar\gamma)$.

Combine, we see the matrix of $e(V_1 \otimes V_2) = \mathrm{sw} \circ e(V_2 \otimes V_1) \circ \mathrm{sw}^{-1}$ is both upper block triangular and lower block triangular, thus must be block diagonal. In other words, $e(V_1\otimes V_2)$ preserves each piece $(V_1)_\gamma\otimes(V_2)_\delta\otimes\bbK$, and the two diagonal blocks must be the same:
\[
e_\gamma(V_1)(y+\hbar\delta)\otimes e_\delta(V_2)(y)\;=\;e_\gamma(V_1)(y)\otimes e_\delta(V_2)(y+\hbar\gamma).
\]
All four factors are invertible. Now if $A\otimes B=A'\otimes B'$ with all four invertible, then we can deduce $A=cA'$ for a non-zero function  $c \in \bbK^\times$. Therefore:
\begin{equation}\label{eq:F}
e_\gamma(V_1)(y+\hbar\delta)=c\cdot e_\gamma(V_1)(y),\qquad\text{for some }c\in\bbK^\times,
\end{equation}
for every weight $\delta$ of every $V_2\in\Rep(G)$.

\emph{Step 3: integrality implies constant}. Fix $V_1$ and $\gamma$, and denote $m:=m_\gamma(V_1) = \dim (V_1)_\gamma$. Obviously, \eqref{eq:F} holds for every $\delta\in\Lambda$. This implies, every ratio, $e_{i,j}/e_{k,l}$, of two matrix entries of $e_\gamma(V_1)$, should be a rational function of $y$, which is invariant under translation by $\hbar\Lambda$. Such a function must be independent of $y$: specialize $\hbar$ to a generic $a\in\bbC^\times$; It's well-known that a rational function on $\frakt^*$, periodic in a full-rank lattice $a\Lambda$, must be a constant. Therefore:
\[
e_\gamma(V_1)=\phi\,N,\qquad \phi\in\bbK^\times,\quad N\in GL_m(\bbC(\hbar)).
\]
We will use the integral property here, that is, the assumption $\calE_\nu (V)\in GL(\sfA_\nu (V))$, in an essential way. Since determinant is invariant under conjugation by $\Top$, we have:
\[ \phi^m \cdot \det(N) = \det e_\gamma(V_1)=\det\calE_\gamma(V_1)\in\Sh^\times=\bbC^\times
\]
Therefore, $\phi^m\in\bbC(\hbar)$, and because $\bbC(\hbar)$ is algebraically closed in $\bbK=\bbC(\hbar)(y)$, we deduce $\phi\in\bbC(\hbar)$. So $e_\gamma(V_1)\in GL\bigl((V_1)_\gamma\otimes\bbC(\hbar)\bigr)$, it does not depend on $y$. \footnote{Without integrality, $\phi$ could be any function with $\phi(y+\hbar\delta)/\phi(y)\in\bbK$. For instance, the ratio of $\Gamma$-functions defines a gauge as in \S\ref{sec:gauge-why}. The integral property is exactly what rules out this.} Now the shifts in Step 2 completely disappear, so we have:
\[
e(V_2\otimes V_1)=e(V_2)\otimes e(V_1).
\]

\emph{Step 4: Tannakian duality.} By previous steps, $e$ is an automorphism of the tensor functor $\Rep(G)\to\Vect_{\bbC(\hbar)}$, $V\mapsto V\otimes\bbC(\hbar)$, which moreover preserves weight spaces. By \cite[Proposition 2.8]{DeligneMilne}, $e=t$ for some $t\in G(\bbC(\hbar))$. Since $t$ preserves all weight spaces of all representations, it centralizes $T$, so $t\in T(\bbC(\hbar))$. Finally, $\nu(t)^{m_\nu}=\det e_\nu(V)\in\bbC^\times$, again $\bbC$ is algebraically closed in $\bbC(\hbar)$, so $\nu(t)\in\bbC^\times$ for every weight $\nu$. We conclude that $t\in T(\bbC)$. Transport this back, we see $\calE_\nu(V)=\nu(t)$. Obviously, such $t$ is unique.
\end{proof}

\subsection*{Proof of Theorem~\ref{thm:DF}}\label{sec:disc-proof}
If $G$ has semisimple rank one, Lemmas~\ref{lem:rank-one-all} and~\ref{lem:int} give $\calE_\nu\in GL(\sfA_\nu)$, and Theorem~\ref{thm:rigid} gives $\calE=t_G$. For a general $G$, we first apply the rank-one case to each minimal Levi subgroup $L_i \subset G$, corresponding to each simple root $\alpha_i$ of $G$. We get $\calE^{L_i} = t_i$ for some $t_i \in T$. The assumption of Proposition~\ref{prop:W-step} is then satisfied with the input from the rank one $L_i$'s. The integrality of  $w_\calG$'s follows from this. Proposition~\ref{prop:simple-walls} gives: $\calE_\nu$ is integral at every wall of every simple root. Applying Weyl reflection, we see that $\calE_\nu (V)$ is integral at every wall of every root. Since $\calE$ is already defined over $\Crs$, being integral at all walls of all roots implies $\calE_\nu\in GL(\sfA_\nu)$ for all $V$ and $\nu$. Finally we apply Theorem~\ref{thm:rigid}, obtaining $\calE = t$. \qed

\section{The proof of the main theorem}\label{sec:gluing}
We come back to our original notation, $L\subset G$ is the Levi subgroup in Theorem~\ref{thm:main}, and $T\subset L \subset G$ is the maximal torus.

\subsection{The corrected Case 2}\label{sec:corrected}
Let $t_L\in T(\bbC)$ be the element of Theorem~\ref{thm:DF} for the group $L$, and let $\calG_L$ be the $\Gamma$-gauge of $L$, see Definition~\ref{def:gauge}. $\calG_L$ acts on $(V\otimes\bM_\hbar)^{\fraku}_\beta$ and also on $(V|_L\otimes\bM_{L,\hbar})^{\fraku_L}_\beta$. We correct $\widetilde\Psi_{V,\beta}$ in Section~\ref{sec:case2}, by choosing the correction $c=t_L\calG_L$ in Proposition~\ref{prop:reduction}, as follows:
\begin{equation}\label{eq:PsiB-corrected}
\widetilde\Psi''_{V,\beta}:=\eqref{eq:case2-geom}\circ\eqref{eq:GR}_\beta\circ t_L\circ\calG_L\circ\Phi^{-1}\colon(V\otimes\bM_\hbar(\psiL))^{\fraku,\psiL}\otimes_{\Zh(\frakl)}\Cbe[\hbar^{-1}]\xrightarrow{\ \sim\ }\H_A(\GrL,i_L^!\IC_V)_\beta[\hbar^{-1}] ,
\end{equation}
where $t_L$ acts on the $\nu$-piece by $\nu(t_L)$. Just as $\widetilde\Psi_{V,\beta}$, it is defined after inverting $\hbar$. Equivalently, $\Phi$ is replaced by $\calG_L^{-1}t_L^{-1}\circ\Phi$. By Lemma~\ref{lem:gauge}(c), this correction is can be interperated as: evaluating an intertwiner on the $\Gamma$-renormalized Whittaker vector, and then renormalizes back:
\[ \xi\mapsto R_\Gamma^{-1}\bigl(f_\xi(\wv_L\cdot\Gamma_L(Y))\bigr), \quad \text{up to a torus element}\, \,t_L \in T(\bbC)
\]
This is a formal expression, and is rational on each weight piece.

\begin{corollary}\label{cor:rs-compat}
$(\widetilde\Psi_{V,\alpha})_\rs=(\widetilde\Psi''_{V,\beta})_\rs$ over $\Crs[\hbar^{-1}]$.
\end{corollary}

\begin{proof}
The family $c_\nu=\nu(t_L)\delta^L_\nu$ consists of units of $\Cbe$, by Lemma~\ref{lem:gauge}(a). Note that they are also invariant under translation along $\frakz^*$. Applying Proposition~\ref{prop:reduction}: we need to show $\BF_L\circ\Phi^L=\GR_L\circ t_L\calG_L$, in other words, $d^L=t_L\calG_L$. This follows from Theorem~\ref{thm:DF}, applying to $L$.
\end{proof}

By Corollary~\ref{cor:extension}, $\widetilde\Psi''_{V,\beta}$ extends uniquely to an isomorphism over $\Cbe$, and still agrees with $\widetilde\Psi_{V,\alpha}$ over $\Crs$. Therefore, from now on we let $\widetilde\Psi''_{V,\beta}$ denotes this extension.

\begin{corollary}\label{cor:WL}
$\widetilde\Psi''_{V,\beta}$ is $W_L$-equivariant.
\end{corollary}

\begin{proof}
The target of $\widetilde\Psi''_{V,\beta}$ is free over $\Cbe$ , by Lemma~\ref{lem:formal}, and also the source by Corollary~\ref{cor:extension}. Both of them embed $W_L$-equivariantly into their localizations over $\Crs$. Over $\Crs$, we have $(\widetilde\Psi''_{V,\beta})_\rs=(\widetilde\Psi_{V,\alpha})_\rs$, by Corollary~\ref{cor:rs-compat}, and the right hand side is $W_L$-equivariant by Proposition~\ref{prop:PsiA}, so does the left hand side.
\end{proof}

\begin{remark}[Multiplicativity]\label{rem:mult}
$\widetilde\Psi''_{V,\beta}$ is a composite of ring homomorphisms: \eqref{eq:case2-geom} (localization, by the lax monoidal structure of $i^!$), \eqref{eq:GR} (by \eqref{eq:GR-reg}), $t_L$, $\calG_L$ (Lemma~\ref{lem:gauge}(c)) and $\Phi^{-1}$ (composition of intertwiners). So $\Psi_\rs$, and hence $\Psi$, is multiplicative, since $\Dh(G/U,\psiL)$ embeds in its localization. This proves that $\Psi$ is an isomorphism of algebras.
\end{remark}

\subsection{Proof of Theorem~\ref{thm:main}}\label{sec:proof-main}
By Proposition~\ref{prop:PsiA} and Corollaries~\ref{cor:extension}, \ref{cor:WL}, $\widetilde\Psi_{V,\alpha}$ and $\widetilde\Psi''_{V,\beta}$ are $W_L$-equivariant and agree over $\ars$. Taking $W_L$-invariants (Lemma~\ref{lem:formal}) and summing $(V^\lambda)^*\otimes(-)$ over $\lambda\in\Lambda^+$ gives compatible isomorphisms
\[
\Psi_\alpha\colon\Dh(G/U,\psiL)_\alpha\xrightarrow{\ \sim\ }\H_{\LO}(\GrL,i_L^!\IC_\reg)_\alpha,\qquad
\Psi''_\beta\colon\Dh(G/U,\psiL)_\beta\xrightarrow{\ \sim\ }\H_{\LO}(\GrL,i_L^!\IC_\reg)_\beta .
\]
By Theorem~\ref{thm:alg-glue} and Proposition~\ref{prop:geom-glue}, $\Psi:=\Psi_\alpha\times_{\Psi_\rs}\Psi''_\beta$ is the isomorphism \eqref{eq:main}. It respects $V$-isotypic components, giving \eqref{eq:main-V}, and it is multiplicative by Remark~\ref{rem:mult}. \qed

\appendix
\section{The fiber \texorpdfstring{$\hbar=0$}{hbar=0}}\label{app:hbar0}
Case~2 (\S\ref{sec:case2}) and the comparison of \S\S\ref{sec:rs-compat}--\ref{sec:discrepancy} are carried out after inverting $\hbar$ (Remark~\ref{rem:KK-asymptotic}). Here we supply what happens at $\hbar=0$. In \S\ref{app:integrality} we show that the discrepancy of \S\ref{sec:discrepancy} is regular and invertible at $\hbar=0$; this is used for $t_G\in T(\bbC)$ in Theorem~\ref{thm:DF}. In \S\ref{app:extension} we extend the corrected Case~2 isomorphism across $\hbar=0$; this is used in the proof of Theorem~\ref{thm:main}.

\subsection{The vanishing order of the Whittaker vector}\label{app:integrality}
We use the notation of \S\ref{sec:shapovalov} and \S\ref{sec:discrepancy}. Here and below, $G$ stands for a reductive group, $|J|$ denotes the PBW degree of $v_J$. Denote $N_\eta:=\sum_{\eta_J=\eta}|J|$, which can be thought as the total PBW degree of the $-\eta$-weight piece, $\bM_{\hbar,-\eta}$, of the asymtopic universal Verma module $\bM_\hbar$ for $G$. Remember that in the proof of Lemma~\ref{lem:J-integral}, we need to control the order of poles at $\hbar$ in Whittaker coefficients $c_J$. This is done by the following Lemma:

\begin{lemma}[Order of poles at $\hbar$]\label{lem:integrality} We have 
$\hbar^{|J|}c_J\in\Crs$ for every $J$. In other words, $\wv_{\psi,\hbar}=\sum_J(\hbar^{-|J|}v_J)\one\cdot c'_J$, with coefficients  $c'_J:=\hbar^{|J|}c_J\in\Crs$.
\end{lemma}

\begin{proof}
Under the isomorphism $\Uh(\frakg)[\hbar^{-1}]\xrightarrow{\sim}U(\frakg)[\hbar^{\pm1}]$, $x\mapsto\hbar x$ ($x\in\frakg$), the module $\bM_\hbar[\hbar^{-1}]$ becomes the Verma module with parameter $\lambda=y/\hbar$. The PBW monomial $v_J$ becomes $\hbar^{|J|}v_J$,  $\tau$, $\mathrm{pr}_\frakt$ are the same as before. Therefore, the asymtopic and the usual Shapovalov forms are related by:
\[
F_\hbar(v_I\one,v_J\one)=\hbar^{|I|+|J|}\,F(v_I,v_J)(y/\hbar),
\]
where $F(v_I,v_J)=\mathrm{pr}_\frakt(e_Iv_J)\in\SymAlg(\frakt)$ is the usual Shapovalov form. The degree of the polynomial $F(v_I,v_J)$ is bounded $\min(|I|,|J|)$, because each element in $\frakt$ must come from a commutator: $[e_\gamma,f_\gamma]=\check\gamma$, which costs an element in $e_I$, and also an element in $v_J$. Also, commuting an element in $\frakt$ with a root vector only gives constant. Therefore
\[
G_{IJ}:=\hbar^{|I|}F(v_I,v_J)(y/\hbar)\in\Sh,\qquad F_\hbar(v_I\one,v_J\one)=G_{IJ}\,\hbar^{|J|}.
\]
The equation solving $\wv_{\psi,\hbar}$, in the proof of Theorem~\ref{thm:whittaker-vector}, is:  $\sum_{\eta_J=\eta}F_\hbar(v_I\one,v_J\one)c_J=\psi(e_I)$ with $\psi(e_I)\in\bbC$, which is the same as:  $\sum_{\eta_J=\eta}G_{IJ}\,c'_J=\psi(e_I)$. Let's compare the determinant of the matrix $(G_{IJ})$ with $F_{IJ}$: 
\begin{equation*}
    \det(G_{IJ})_{\eta_I=\eta_J=\eta}= \hbar^{\sum_{I: \eta_I = \eta} |I|}\det F_\eta(y/\hbar)= 
    \hbar^{N_\eta}\det F_\eta(y/\hbar)
\end{equation*}
By Lemma~\ref{lem:KK}, $\det F_\eta$ is, up to a constant, a product of $\sum_{\gamma,n}\mathscr P(\eta-n\gamma)$ linear factors $\check{\gamma} + (\langle\rho,\check\gamma\rangle-n)$. We claim that:
\begin{equation*}
    \sum_{\gamma,n}\mathscr P(\eta-n\gamma) = N_\eta
\end{equation*}
Indeed, both sides count the pairs $(J,m)$ with $\eta_J=\eta$ and $1\le m\le a_\gamma(J)$, where $v_J=\prod_\gamma f_\gamma^{a_\gamma(J)}$. Therefore the powers of $\hbar$ cancels, and we have:
\[
\det G=\mathrm{const}\cdot\prod_{\gamma\in R^+}\prod_{n\ge1}\bigl(\check\gamma+(\langle\rho,\check\gamma\rangle-n)\hbar\bigr)^{\mathscr P(\eta-n\gamma)}
\]
is a product of walls, which is a unit in $\Crs$; so $G$ is invertible over $\Crs$ and $c'=G^{-1}\bigl(\psi(e_I)\bigr)_I$ has entries in $\Crs$. 
\end{proof}

\begin{proposition}[$\Phi$ is regular and invertible at $\hbar=0$]\label{prop:Phi-integral}
Let $\Phi\colon\sfA_\rs\to\sfW_\rs$ be the map \eqref{eq:Phi-G} composed with $V\otimes\phi_\hbar$ and Lemma~\ref{lem:GK}, a priori defined over $\Crs[\hbar^{-1}]$. Then $\Phi(\sfA\otimes\Crs)\subset\sfW\otimes\Crs$, and $\Phi\otimes(\Sh)_{(\hbar)}$ is an isomorphism. Thus, $d$ and $\calE$ are automorphisms of $\sfA_\rs$ defined over $\Crs$, and they are automorphisms of $\sfA\otimes(\Sh)_{(\hbar)}$.
\end{proposition}

\begin{proof}
Recall $\Crs\subset(\Sh)_{(\hbar)}$, since no wall lies in $\hbar\Sh$.

\emph{Step 1: $\Top$ is an isomorphism at $(\hbar)$}. On the $\nu$-component $\Top=\Xi^G_T\circ\GR$, see the proof of Lemma~\ref{lem:G2}. $\GR$ is an isomorphism, and $\Xi^G_T$ becomes isomorphism over $\Crs\subset(\Sh)_{(\hbar)}$, by Proposition~\ref{prop:geom-loc}.

\emph{Step 2: $\coinv$ is an isomorphism at $(\hbar)$}. Write $\xi=\sum_Kv_K\otimes y_K\one$,  then:
\[
\coinv(\xi)-\Top(\xi)=\sum_{K\neq\emptyset}S(y_K)v_K\in\hbar\,(V\otimes\Sh)
\]
This is because an element $x \in \Uh(\frakg)$ acts on $V$ by $\hbar x$. So $\Top^{-1}\circ\coinv=\id+\hbar(\cdots)$ on $\sfA_\nu\otimes(\Sh)_{(\hbar)}$, which is invertible in the local ring $(S_\hbar)_{(\hbar)}$.

\emph{Step 3: $\Phi$ is an isomorphism at $(\hbar)$}. By Lemma~\ref{lem:class}(b) and Remark~\ref{rem:model}, $J\circ\Phi=\coinv$ over $\Crs[\hbar^{-1}]$. By Lemma~\ref{lem:J-integral}, $J$ is an isomorphism over $\Crs$. Thus, 
\[
\Phi=J^{-1}\circ\coinv: \sfA\otimes\Crs \xrightarrow{\sim} \sfW\otimes\Crs
\]
becomes an isomorphism at $(\hbar)$ by Step 2. Since $\GR$, $\BF$ and $\calG_G$ \footnote{for $\calG$, recall that $\delta_\nu$ is a product of walls, so it is invertible at $\hbar$.} are defined over $\Crs$ and invertible at $(\hbar)$, the same holds for $d=\GR^{-1}\circ\BF\circ\Phi$ and $\calE=d\,\calG_G^{-1}$.
\end{proof}

\subsection{Extension across \texorpdfstring{$\hbar=0$}{hbar=0}}\label{app:extension}
The map $\widetilde\Psi''_{V,\beta}$ of \eqref{eq:PsiB-corrected} is defined after inverting $\hbar$. The points of $\abe$ with $\hbar=0$ that are not already in $\aal$ lie on $\{\hbar=0\}\cap\bigcup_\alpha\{\check\alpha=0\}$. This is again a set of codimension two, so we are going to use the philosophy of extension once more. The following lemma does this: it extends $\widetilde\Psi''_{V,\beta}$ across $\{\hbar=0\}\cap\bigcup_\alpha\{\check\alpha=0\}$ using Case 1. The classical input enters only through Case 1, that is, through Theorem~\ref{thm:Theta}, which uses \cite[Proposition 3.4.12]{macerato2023leviequivariantrestrictionsphericalperverse} at $\hbar=0$.

\begin{lemma}[Extension across $\hbar=0$]\label{lem:extension}
Let $E$ be a free $\Sh$-module and $M\subset E$ a submodule, such that $M\cap\hbar E=\hbar M$. Then, inside $M_\rs[\hbar^{-1}]$,
\[
M_\beta\;=\;M_\beta[\hbar^{-1}]\ \cap\ M_\rs .
\]
\end{lemma}

\begin{proof}
All maps $M_\beta\to M_\beta[\hbar^{-1}]$, $M_\beta\to M_\rs$, $M_\rs\to M_\rs[\hbar^{-1}]$ are injective, since $M\subset E$ is torsion-free. Because $M\cap\hbar E=\hbar M$, we see that $M/\hbar M\to E/\hbar E$ is injective. Because $E/\hbar E$ is free over $\SymAlg(\frakt)$, multiplication by any non-zero element of $\SymAlg(\frakt)$ is injective on $M/\hbar M$. The same holds after localization: for $f$ a product of $\alpha$-walls, its image $\bar f\in\SymAlg(\frakt)$ is a product of coroots, multiplication by $\bar f$ is injective on $M_\beta/\hbar M_\beta=(M/\hbar M)_\beta$.

Let $x\in M_\beta[\hbar^{-1}]\cap M_\rs$ and write $x=m/\hbar^n=m'/f$, where $m,m'\in M_\beta$, $n\ge0$, and $f$ is a product of $\alpha$-walls. Then $fm=\hbar^nm'$ holds in $M_\beta$. Modulo $\hbar$, we get: $\bar f\bar m=0$ in $M_\beta/\hbar M_\beta$. This implies $\bar m=0$.
Thus, $\hbar$ divides $m$, and we can write $m =\hbar m_1$ then use induction.
\end{proof}

\begin{corollary}\label{cor:extension}
Let $V\in\Rep(G)$. The isomorphism $\widetilde\Psi''_{V,\beta}$ of \eqref{eq:PsiB-corrected}, defined over $\Cbe[\hbar^{-1}]$, extends uniquely to an isomorphism of $\Cbe$-modules
\[
\widetilde\Psi''_{V,\beta}\colon(V\otimes\bM_\hbar(\psiL))^{\fraku,\psiL}\otimes_{\Zh(\frakl)}\Cbe\xrightarrow{\ \sim\ }\H_A(\GrL,i_L^!\IC_V)_\beta ,
\]
and the extension still satisfies $(\widetilde\Psi''_{V,\beta})_\rs=(\widetilde\Psi_{V,\alpha})_\rs$. In particular the domain is free of rank $\dim V$ over $\Cbe$.
\end{corollary}

\begin{proof}
Put $E:=\Sh\otimes_{\Zh(\frakl)}\bigl(V\otimes\bM_\hbar(\psiL)\bigr)$ and $M:=\Sh\otimes_{\Zh(\frakl)}\bigl(V\otimes\bM_\hbar(\psiL)\bigr)^{\fraku,\psiL}$, so that the domain of $\widetilde\Psi''_{V,\beta}$ is $M_\beta[\hbar^{-1}]$. It's clear that $E$ is a free $\Sh$-module, and $M\cap\hbar E=\hbar M$. Let $N:=\H_A(\GrL,i_L^!\IC_V)$, which is free over $\Sh$, by Lemma~\ref{lem:formal}. In this case, Lemma~\ref{lem:extension} applies to $M$ and to $N$.

Let $m\in M_\beta$. Since $M_\beta\subset M_\rs$ and $\widetilde\Psi_{V,\alpha}$ is defined over $\Cal$, they are of course defined over $\Crs$, so we can talk about the element $\widetilde\Psi_{V,\alpha}(m)\in N_\rs$. By Corollary~\ref{cor:rs-compat}, it equals $\widetilde\Psi''_{V,\beta}(m)\in N_\beta[\hbar^{-1}]$ inside $N_\rs[\hbar^{-1}]$. Apply Lemma~\ref{lem:extension} to $N$, it lies in $N_\beta$. So $\widetilde\Psi_{V,\alpha}$ restricts to a map $M_\beta\to N_\beta$, which extends $\widetilde\Psi''_{V,\beta}$ and agrees with $\widetilde\Psi_{V,\alpha}$ over $\Crs$. Then we apply the same argument to their inverses: $(\widetilde\Psi_{V,\alpha})^{-1}$ and $(\widetilde\Psi''_{V,\beta})^{-1}$, and Lemma~\ref{lem:extension} to $M$. We see that it maps $M_\beta$ onto $N_\beta$.
\end{proof}

\section{How the \texorpdfstring{$\Gamma$}{Gamma}-gauge is determined}\label{sec:gauge-why}
Theorem~\ref{thm:DF} says that the discrepancy $d$ equals $\calG_G$ up to an element of $T(\bbC)$. The rank one computation in Theorem~\ref{thm:rank-one} only determines what happens along simple roots. We explain here why the additional structures force it to be the form of Definition~\ref{def:gauge}. The proof of Theorem~\ref{thm:DF} also confirms the $\Gamma$-gauge is the correct one.

Consider a \emph{gauge}, i.e.\ a family $c=(c_\nu)_\nu$, $c_\nu\in\Crs^\times$, acting on $\sfA_\rs$ by $\xi\mapsto\xi\cdot c_\nu$ on $\sfA_\nu$. We want $c$ to satisfy the following constraints.

\emph{(i) Multiplicativity}. Because $d$ is multiplicative, by Lemma~\ref{lem:E-mult}, and we also want  $d\circ c^{-1}$ to be multiplicative after gauging, $c$ must also be multiplicative. By \eqref{eq:comp-shift}, $(\xi'\cdot c_\mu)\circ(\xi\cdot c_\nu)=(\xi'\circ\xi)\cdot c_\mu(y+\hbar\nu)\,c_\nu(y)$ for $\xi\in\sfA_\nu(V_1)$, $\xi'\in\sfA_\mu(V_2)$, while multiplicativity asks for $(\xi'\circ\xi)\cdot c_{\mu+\nu}$. So a multiplicative $c$ satisfies the cocycle condition \eqref{eq:gauge-cocycle}, and conversely (Corollary~\ref{cor:comp-verma}(c)). There are two kinds of obvious solutions: (1) characters, $c_\nu=\nu(t)$ with $t\in T(\bbC)$; (2) coboundaries
\[
c_\nu(y)=\frac{\varphi(y)}{\varphi(y+\hbar\nu)}
\]
for a single function $\varphi$ of $y$, possibly a transcendental one, but with those ratios rational. A coboundary gauge is conjugation by right multiplication by $\varphi$ on $\bM_\hbar$ (Corollary~\ref{cor:comp-verma}(c)), exactly as in Lemma~\ref{lem:gauge}(c). The point of (i) is that we are not going to choose $\delta_\nu$ weight by weight. Instead we choose one function $\varphi$, and up to $T(\bbC)$, this is all the freedom we will use.

\emph{(ii) Rank one.} For $\mathfrak{sl}_2$ the gauge is computed: by Theorem~\ref{thm:rank-one}, $d_\nu=\delta_\nu=P_\nu(Y)^{-1}$ up to constants, where:
\[
P_\nu(Y)^{-1}=\frac{\varphi_1(y)}{\varphi_1(y+\hbar\nu)},\qquad \varphi_1(y):=\hbar^{Y/\hbar}\,\Gamma(1+Y/\hbar),\quad Y=\langle y+\hbar\rho,\check\alpha\rangle,
\]
Recall that $Y_\nu=Y+\hbar\nu$ and $\Gamma(1+z+k)/\Gamma(1+z)=\prod_{m=1}^k(z+m)$. Thus, along each simple root $\alpha_i$, the function $\varphi$ shoulf be $\varphi_1(\langle y+\hbar\rho,\check\alpha_i\rangle)$.

\emph{(iii) Locality.} The zeros and poles of $d_\nu$ lie on walls, and near a wall of $\gamma$ only the root $\gamma$ should matter; for simple $\gamma$ this is what \S\ref{sec:disc-simple} proves. This suggests the following guess:
\[
\varphi(y)=\prod_{\gamma\in R^+}\varphi_\gamma\bigl(\langle y+\hbar\rho,\check\gamma\rangle\bigr),
\]
which is a product over all positive roots of functions of one variable, with $\varphi_{\alpha_i}=\varphi_1$ by (ii).

\emph{(iv) Weyl symmetry.} The comparison through Levi subgroups, in Step 4 of Lemma~\ref{lem:WD}, requires that: for each simple root $\alpha_i$, with $L=L_i$ and $w=s_i$, the part $\calG_{G/L}$ of the gauge coming from the roots in $R^+\smallsetminus\{\alpha_i\}$ should commute with $w_\Phi$. This requires the invariant under $(y,\nu)\mapsto(s_i\cdot y,s_i\nu)$. Since \footnote{The $\rho$-shift in $Y=y+\hbar\rho$ is what makes the dot action of $W$ act on $y$ in $\langle Y,\check\gamma\rangle$ equivalent to permuting the coroots.}:
\[
\langle s_i\cdot y+\hbar\rho,\check\gamma\rangle=\langle y+\hbar\rho,s_i\check\gamma\rangle,\qquad\langle s_i\nu,\check\gamma\rangle=\langle\nu,s_i\check\gamma\rangle,
\]
and $s_i$ permutes $R^+\smallsetminus\{\alpha_i\}$, this holds when $\varphi_\gamma=\varphi_{s_i\gamma}$ for all $\gamma\in R^+\smallsetminus\{\alpha_i\}$; and for a product as in (iii) this is essentially the only way it can hold. Now every non-simple $\gamma\in R^+$ has some $\alpha_i$ with $\langle\gamma,\check\alpha_i\rangle>0$, and then $s_i\gamma\in R^+$ has smaller height. Thus, by induction on the height, it followa $\varphi_\gamma=\varphi_{\alpha_j}=\varphi_1$ for every $\gamma\in R^+$. Therefore,
\[
\varphi(y)=\prod_{\gamma\in R^+}\hbar^{\langle Y,\check\gamma\rangle/\hbar}\,\Gamma(1+\langle Y,\check\gamma\rangle/\hbar)=\Gamma_G(Y),\qquad
\frac{\varphi(y)}{\varphi(y+\hbar\nu)}=\frac{\Gamma_G(Y)}{\Gamma_G(Y_\nu)}=\delta_\nu .
\]
This is the same induction on the height of roots as in \S\ref{sec:disc-weyl}.

\emph{(v) Rationality and units.} Each factor of $\delta_\nu$ is $\Gamma(1+z)/\Gamma(1+z+k)$ with $z=\langle Y,\check\gamma\rangle/\hbar$, $k=\langle\nu,\check\gamma\rangle\in\bbZ$. It is a product of the linear forms $\langle y,\check\gamma\rangle+(\langle\rho,\check\gamma\rangle+m)\hbar$, or the inverse. Because these are walls, so $\delta_\nu\in\Crs^\times$, and $\delta^L_\nu\in\Cbe^\times$, since only roots of $L$ occur, by Lemma~\ref{lem:gauge}(a). This is also the fact that allows us to use such a correction in Case 2.

So (i)--(v) determine the gauge up to a character $\nu\mapsto\nu(t)$, $t\in T(\bbC)$, which is the constant $t_G$ of Theorem~\ref{thm:DF}. The proof also checks the answer: if the factors for non-simple roots are designed differently, $\calG_{G/L}$ will not commute with $w_\Phi$,  and $\calE_\nu$ will not be integral at the walls of non-simple roots.

\subsection*{The classical limit}
It is also interesting to see what happens to the correction in the classical limit. It is given by the folloing.
\begin{remark}[The case $\hbar=0$]\label{rem:hbar0}
After the extension of \S\ref{app:extension} we may set $\hbar\to0$. The classical limit of the correction is
\[
\delta^L_\nu\longrightarrow\prod_{\beta\in R_L^+}\langle y,\check\beta\rangle^{-\langle\nu,\check\beta\rangle}=\nu(t(y))^{-1},\qquad t(y):=\prod_{\beta\in R_L^+}\check\beta\bigl(\langle y,\check\beta\rangle\bigr).
\]
So classically the correction is the action of a $T$-valued function, which is absorbed in the freedom of normalizing the classical Kostant--Whittaker isomorphism that \cite{macerato2023leviequivariantrestrictionsphericalperverse} fixes by (4.4.13). This is why no correction appears in \cite[proof of Theorem 4.4.5]{macerato2023leviequivariantrestrictionsphericalperverse}. The $\calG_L$ here can be viewed as its quantization. Quantum anomaly appears because there is no canonical splitting of the quantized Whittaker cotangent bundle $D_\hbar(G/U,\psi)$.
\end{remark}

\bibliographystyle{amsalpha}
\bibliography{reference}

@misc{macerato2023leviequivariantrestrictionsphericalperverse,
      title={Levi-Equivariant Restriction of Spherical Perverse Sheaves}, 
      author={Mark Macerato},
      year={2023},
      eprint={2309.07279},
      archivePrefix={arXiv},
      primaryClass={math.RT},
      url={https://arxiv.org/abs/2309.07279}, 
}

@misc{bezrukavnikov2026equivariantsatakecategorykostantwhittaker,
      title={Equivariant Satake category and Kostant-Whittaker reduction}, 
      author={Roman Bezrukavnikov and Michael Finkelberg},
      year={2026},
      eprint={0707.3799},
      archivePrefix={arXiv},
      primaryClass={math.RT},
      url={https://arxiv.org/abs/0707.3799}, 
}

@misc{ginzburg2014differentialoperatorsguaffine,
      title={Differential operators on G/U and the affine Grassmannian}, 
      author={Victor Ginzburg and Simon Riche},
      year={2014},
      eprint={1306.6754},
      archivePrefix={arXiv},
      primaryClass={math.RT},
      url={https://arxiv.org/abs/1306.6754}, 
}

@misc{ginzburg2022differentialoperatorsgugelfandgraev,
      title={Differential operators on G/U and the Gelfand-Graev action}, 
      author={Victor Ginzburg and David Kazhdan},
      year={2022},
      eprint={1804.05295},
      archivePrefix={arXiv},
      primaryClass={math.RT},
      url={https://arxiv.org/abs/1804.05295}, 
}

@article{Kostant1978,
author = {Kostant, Bertram},
journal = {Inventiones mathematicae},
pages = {101-184},
title = {On Whittaker Vectors and Representation Theory.},
url = {http://eudml.org/doc/142586},
volume = {48},
year = {1978},
}

@article{Milicic1997,
    author = {Dragan Milicić and Wolfgang Soergel},
    title = {The composition series of modules induced from Whittaker modules},
    journal ={ Commentarii Mathematici Helvetici} ,
    year = {1997}
}

@article{MCDOWELL1985161,
title = {On modules induced from Whittaker modules},
journal = {Journal of Algebra},
volume = {96},
number = {1},
pages = {161-177},
year = {1985},
issn = {0021-8693},
doi = {https://doi.org/10.1016/0021-8693(85)90044-4},
url = {https://www.sciencedirect.com/science/article/pii/0021869385900444},
author = {Edward McDowell}
}

@ARTICLE{Backelinfunctor,
  author={Backelin, Erik},
  journal={International Mathematics Research Notices}, 
  title={Representation of the category O in Whittaker categories}, 
  year={1997},
  volume={1997},
  number={4},
  pages={153-172},
  doi={10.1155/S1073792897000111}}

@article{KAC197997,
title = {Structure of representations with highest weight of infinite-dimensional Lie algebras},
journal = {Advances in Mathematics},
volume = {34},
number = {1},
pages = {97-108},
year = {1979},
issn = {0001-8708},
doi = {https://doi.org/10.1016/0001-8708(79)90066-5},
url = {https://www.sciencedirect.com/science/article/pii/0001870879900665},
author = {V.G Kac and D.A Kazhdan}
}

@misc{braverman2016instantonmodulispacesmathscr,
      title={Instanton moduli spaces and $\mathscr W$-algebras}, 
      author={Alexander Braverman and Michael Finkelberg and Hiraku Nakajima},
      year={2016},
      eprint={1406.2381},
      archivePrefix={arXiv},
      primaryClass={math.QA},
      url={https://arxiv.org/abs/1406.2381}, 
}

@misc{eberhardt2024proofsjantzenconjecturewhittaker,
      title={Two proofs of a Jantzen Conjecture for Whittaker Modules}, 
      author={Jens Niklas Eberhardt and Anna Romanov},
      year={2024},
      eprint={2407.16330},
      archivePrefix={arXiv},
      primaryClass={math.RT},
      url={https://arxiv.org/abs/2407.16330}, 
}

@book{humphreysrepresentations,
  title={Representations of Semisimple Lie Algebras in the BGG Category O},
  author={Humphreys, J.E.},
  isbn={9780821872529},
  url={https://books.google.com/books?id=8GCP4Ng6risC},
  publisher={American Mathematical Soc.}
}

@misc{mirkovic2018geometriclanglandsdualityrepresentations,
      title={Geometric Langlands duality and representations of algebraic groups over commutative rings}, 
      author={I. Mirkovic and K. Vilonen},
      year={2018},
      eprint={math/0401222},
      archivePrefix={arXiv},
      primaryClass={math.RT},
      url={https://arxiv.org/abs/math/0401222}, 
}

@article{skryabin2002,
  author  = {Skryabin, Serge},
  title   = {A category equivalence},
  note    = {Appendix to A. Premet, Special transverse slices and their enveloping algebras},
  journal = {Advances in Mathematics},
  volume  = {170},
  number  = {1},
  pages   = {1--55},
  year    = {2002}
}

@incollection{DeligneMilne,
  author    = {Deligne, Pierre and Milne, James S.},
  title     = {Tannakian categories},
  booktitle = {Hodge cycles, motives, and {S}himura varieties},
  series    = {Lecture Notes in Mathematics},
  volume    = {900},
  pages     = {101--228},
  publisher = {Springer},
  year      = {1982}
}

@misc{iritani2023gammaclassesquantumcohomology,
      title={Gamma classes and quantum cohomology}, 
      author={Hiroshi Iritani},
      year={2023},
      eprint={2307.15938},
      archivePrefix={arXiv},
      primaryClass={math.AG},
      url={https://arxiv.org/abs/2307.15938}, 
}

@misc{okounkov2017lecturesktheoreticcomputationsenumerative,
      title={Lectures on K-theoretic computations in enumerative geometry}, 
      author={Andrei Okounkov},
      year={2017},
      eprint={1512.07363},
      archivePrefix={arXiv},
      primaryClass={math.AG},
      url={https://arxiv.org/abs/1512.07363}, 
}

\Addresses

\end{document}